\documentclass[fleqn,trsc,nonblindrev]{informs3noheader}

\OneAndAHalfSpacedXI 
\usepackage{endnotes}
\usepackage{csquotes}
\let\footnote=\endnote
\let\enotesize=\normalsize
\def\notesname{Endnotes}%
\def\makeenmark{$^{\theenmark}$}
\def\enoteformat{\rightskip0pt\leftskip0pt\parindent=1.75em
	\leavevmode\llap{\theenmark.\enskip}}
\usepackage{amsmath, threeparttable, float, color, hyperref, epstopdf, multirow, makecell}
\usepackage{graphicx}
\usepackage{multirow}
\usepackage{hhline}
\usepackage{threeparttable}
\usepackage[small, margin=1cm]{caption}
\usepackage{color}
\definecolor{strcolor}{rgb}{0.6, 0.2, 0.6}
\definecolor{commentcolor}{rgb}{0.3125, 0.5, 0.3125}
\definecolor{keycol}{rgb}{0, 0, 1}

\usepackage{subcaption}
\usepackage[short, nocomma]{optidef}
\newcounter{mycounter}
\usepackage{natbib}
 \bibpunct[, ]{(}{)}{,}{a}{}{,}%
 \def\bibfont{\small}%
 \def\newblock{\ }%
\usepackage{caption}
\usepackage{subcaption}
\usepackage{algorithm}
\usepackage{algorithmic}
\usepackage{tikz}
\usepackage{bbm}

\usepackage{listings}
\usepackage{url}

\newcommand {\bea}{\begin{eqnarray}}
	\newcommand {\eea}{\end{eqnarray}}

\def\blot{\quad \mbox{$\vcenter{ \vbox{ \hrule height.4pt
				\hbox{\vrule width.4pt height.9ex \kern.9ex \vrule width.4pt}
				\hrule height.4pt}}$}}

\usepackage{natbib}
\bibpunct[, ]{(}{)}{,}{a}{}{,}%
\def\bibfont{\fontsize{10}{13}\selectfont}%
\def\newblock{\ }%

\TheoremsNumberedThrough     
\EquationsNumberedThrough    

\gdef\AQ#1{}
\gdef\CQ#1{}
\begin{document}
\def\COPYRIGHTHOLDER{INFORMS}%
\def\COPYRIGHTYEAR{2017}%
\def\DOI{\fontsize{7.5}{9.5}\selectfont\sf\bfseries\noindent https://doi.org/10.1287/opre.2017.1714\CQ{Word count = 9740}}

\RUNAUTHOR{Yang et~al.} %

\RUNTITLE{Integrated demand-side management and timetabling}

\TITLE{Integrated demand-side management and timetabling for an urban rail transit line: A Benders decomposition approach}

\ARTICLEAUTHORS{
\AUTHOR{Lixing Yang,\textsuperscript{a,b*}, Yahan Lu,\textsuperscript{a,c*} Jiateng Yin,\textsuperscript{a} Shadi Sharif Azadeh\textsuperscript{c} }
\AFF{$^{a}$ School of Systems Science, Beijing Jiaotong University, Beijing, 100044, China; \\
$^{b}$ Hebei Key Laboratory of Future Urban Intelligent Traffic Management, Beijing Jiaotong University, Beijing 100044, China; \\
$^{c}$ Department of Transport \& Planning, Delft University of Technology, Netherlands}
}

\ABSTRACT{The intelligent upgrading of metropolitan rail transit systems has made it feasible to implement demand-side management policies that integrate multiple operational strategies in practical operations. However, the tight interdependence between supply and demand necessitates a coordinated approach combining demand-side management policies and supply-side resource allocations to enhance the urban rail transit ecosystem. In this study, we propose a mathematical and computational framework that optimizes train timetables, passenger flow control strategies, and trip-shifting plans through the pricing policy. Our framework incorporates an emerging trip-booking approach that transforms waiting at the stations into waiting at home, thereby mitigating station overcrowding. Additionally, it ensures service fairness by maintaining an equitable likelihood of delays across different stations. We formulate the problem as an integer linear programming model, aiming to minimize passengers' waiting time and government subsidies required to offset revenue losses from fare discounts used to encourage trip shifting. To improve the computational efficiency, we develop a Benders decomposition-based algorithm within the branch-and-cut method, which decomposes the model into train timetabling with partial passenger assignment and passenger flow control subproblems. We propose valid inequalities based on our model’s properties to strengthen the linear relaxation bounds at each node of the branch-and-bound tree. Computational results from proof-of-concept and real-world case studies on the Beijing metro show that our solution method outperforms commercial solvers in terms of computational efficiency. We can obtain high-quality solutions, including optimal ones, at the root node with reduced branching requirements thanks to our novel decomposition framework and valid inequalities. Our integrated optimization approach reduces the fleet size for operators by at least 8.33\% and decreases the waiting time of passengers on the tested instances, thereby validating the effectiveness of our proposed methods.
}


\maketitle
	
\section{Introduction}
\label{sec:introdunction}

According to the United Nations (\citeyear{un75}), there will be 9.7 billion people worldwide by the year 2050, and it is anticipated that about 70\% of the world's population by 2050—up from 55\% in 2019—will reside in metropolitan regions. This tendency causes the features of passengers in urban rail transit (URT) systems to change over time, such as their quantitative composition and spatial-temporal distributional characteristics. Take the Beijing rail transit system as an example, it transported 210 million people in 2012 while more than 4.53 billion in 2019. In this context, congestion within URT systems in metropolitan areas has become the norm in operations. \textit{Congestion} refers to the phenomenon caused bythe  high density of passengers on platforms and trains. In addition to making passengers uncomfortable and thereby lowering their satisfaction, congestion can cause delays, disruptions, and disturbances.

Congestion is caused by the mismatch between continuously evolving passenger demand and the relatively stable capacity of the URT system. To address this challenge, URT authorities and operators may consider several strategies. Initially, they could invest in technical and social strategies that involve infrastructure and operational improvements, such as expanding the number of lines and adjusting train timetables to increase train frequency \citep[]{Jin2023}.  For example, despite the headway on Beijing Metro Line 10 being reduced to one minute and 45 seconds during the morning peak hour, there were still many standing passengers inside trains \citep[]{china}. This highlights the significance of \textit{demand-side management methods}.

From the demand perspective, one approach is the implementation of a \textit{passenger flow control strategy}, which controls the boarding rate at stations along the line to maintain spare capacity for downstream stations, thereby optimizing capacity utilization. For instance, in 2019, the Beijing metro regularly applied this strategy at 91 stations on weekdays. Another approach is \textit{congestion pricing}, which differentiates peak and off-peak fares. This strategy has been effectively applied in serveral locations, such as London's congestion charging for road traffic since 1999 \citep[]{London}, Singapore's electronic road pricing in busy areas \citep[]{Singapore}, and the Netherlands where off-peak railway discounts are offered \citep[]{NS}. Congestion pricing essentially aims to achieve \textit{trip shifting}, encouraging passengers who intend to travel during peak periods to shift to less congested times, thereby balancing the demand across time periods. The effectiveness of managing congestion through this strategy has been demonstrated in \cite{BAO2023}. Thereafter, we define \textit{passenger directing} as the combination of the passenger flow control strategy to limit boarding rates and a trip shifting strategy that incentivizes passengers to shift their trips. A more moderate and emerging alternative within URT systems is the \textit{trip booking} strategy, which allows a limited number of passengers to make a reservation for the following day’s travel, bypassing queues outside the station for a passenger flow-controlled entry permit. This reservation system, as explored in our study, enables passengers to reserve a time slot that guarantees platform access and immediate train boarding, rather than requiring them to book specific seats. The Beijing metro applied this strategy in 2020, and by 2021, it reportedly saved passengers a cumulative 40,000 hours in waiting times, with a 88 percent approval rate from users \citep[]{ChinaNews}. Given the fundamental mismatch between supply and demand that results in congestion, integrating these operational and demand-side strategies presents a more effective approach to managing URT operations amid growing congestion worldwide, which is the focus of this paper.

However, the development of passenger flow control and trip booking strategies in practical operations currently stands apart from the production of train timetables and relies entirely on operators' manual experience. That is, the passenger flow control and trip booking plans are set as input when producing train timetables. In particular, the trip-booking strategy, enabled by advancements in intelligent URT systems, is an emerging method still in the proof-of-concept phase in real-world applications. This stage calls for innovative approaches to realize its potential benefits \citep{XIA2024}. Moreover, these demand-side management approaches are mostly managed on a station-by-station basis, which limits the potential for achieving line-wide connectivity \citep[]{Meng2022}.

While existing research has developed mathematical models and solution algorithms for both the passenger flow control problem and the joint optimization of train timetabling and passenger flow control problem, the integration of train timetabling, passenger directing, and trip booking problem is hardly addressed in the literature. However, experimental results from related studies \citep[e.g., ][]{Shi2018, Lu2023} consistently indicate that although the number of detained passengers decreases with the implementation of the timetabling and passenger flow control policies compared to scenarios without them, the phenomenon of oversaturation remains. These findings highlight the necessity to collaboratively optimize additional demand-side management methods. Moreover, the algorithms proposed for the joint optimization models mainly focus on the solution efficiency,  falling short of finding the exact solutions required in operational planning. 

In this paper, we aim to close these gaps by developing an integrated demand-side management and timetabling approach with an exact solution method. Our framework is grounded on three key layers: government, metro corporations, and passengers. The goal is to provide proof-of-concept insights to operators and governments through optimal solutions obtained by the exact solution method. We don't consider game theory approaches or bi-level frameworks with passenger behaviors, as they often lead to models that are challenging to solve exactly within an acceptable timeframe. We also omit the sequential solution approach because it lacks feedback between decisions and often leads to suboptimal solutions. By optimizing timetables, passenger flow control strategies, and providing discounts on ticket prices to encourage passengers to shift trips, the metro corporation maximizes its passenger-oriented resource allocation. Passengers, by shifting their travel plans and making reservations, minimize waiting times. Meanwhile, the government makes up for the metro corporation's loss of fare revenue by providing additional subsidies. 

Specifically, this paper formally addresses \textit{the integrated optimization of trip booking, passenger directing, and train timetabling} (BDTT) problem from a system-optimal perspective, incorporating deterministic passenger demand and a time-varying reservation slot allocation plan. We formulate the BDTT problem as an integer programming (ILP) model that captures the interdependencies among train timetables, passengers with reserved trips, passengers without reservations, and passengers’ trip-shifting plans. The objective is to minimize the weighted sum of passengers' waiting time and governments' additional subsidies provided to incentivize trip shifting. Service fairness among passengers in terms of the possibility of boarding the first train, both with and without reservations, and across different stations is modeled as constraints. To effectively solve the proposed BDTT problem, we introduce a novel decomposition framework for our problem. Our approach is based on Benders decomposition which decomposes the model into a timetabling subproblem combined with partial passenger allocations and a passenger flow control subproblem. We validate our proposed formulation and solution methodology through proof-of-concept and real-world case studies.

The remainder of this paper is organized as follows: Section \ref{sec:literatureReview} provides an overview of relevant literature. Thereafter, Section \ref{sec:problemDescription} presents a detailed description of the studied problem. In Section \ref{sec:model}, we formulate the problem as an integer programming model. In Section \ref{sec:solutionMethods}, we propose the linearization procedure, the tailored decomposition approach and an exact solution method. Section \ref{sec:caseStudy} presents numerical results and managerial insights. Lastly, we conclude in Section \ref{sec:conclusion}.

\section{Literature review}
\label{sec:literatureReview}

In this section, we give an overview on related research. In Section \ref{sec:literatureTimetabling}, we delve into the collaborative optimization of train timetabling and passenger flow control. Section \ref{sec:literatureBooking} describes how the trip reservation and pricing policies are optimized to improve the matching of supply and demand in existing literature. In Section \ref{sec:Contribution}, we compare our method with the reviewed state-of-art methods.

\subsection{Train timetabling problem combined with passenger flow control}
\label{sec:literatureTimetabling}

With the surge in passenger demand in recent years, the train timetbaling problem combined with passenger flow control has become a hot research topic. A common objective is to minimize passengers' waiting time \citep{LI2017, Liang2023, HU2023, YuanPartC}. 

One of the first related models considering the time-dependent passenger demand has been introduced by \citet{Shi2018}. The authors formulated an ILP model to determine the train timetable and the time-dependent passenger flow control strategy, which is solved by a hybrid algorithm combing the local search procedure with CPLEX. Further, \citet{RMLiu2020} proposed a mixed-integer nonlinear programming (MINLP) model for this problem and designed a Lagrangian relaxation-based solution method. Considering uncertain passenger demand, \citet{Lu2022} formulated a two-stage distributionally robust optimization model where the probability of stochastic scenarios is partially known in advance. \citet{Lu2023} formulated three ILP models for this problem, aiming to generate a reliable train timetable and the train-based passenger flow control strategy to cope with the demand uncertainty in reality. The proposed models based on the Light Robustness technique can be solved directly by GUROBI, while the scenario-based stochastic programming model is also addressed by the hybrid algorithm. In summary, due to the complexity of models that consider supply-demand coupling relations, most of these studies develop heuristic algorithms to solve the proposed models for the integrated optimization problem of timetabling and passenger flow control.

\subsection{Trip reservation and congestion pricing}
\label{sec:literatureBooking}

Trip reservation is a demand management strategy frequently studied in the field of airline marketing and road transportation. In the air transportation ares, it usually requires all passengers to buy tickets in advance and thus automatically make reservations \citep{Air1, Air2, Barz2016}. In contrast, there are a limited number of reservations in the URT system that will not allow all demand to be met, where the allocation plan should be dynamic to match the time-varying characteristics of the passenger demand. For a comprehensive review of the trip reservation in the field of road transportation, we recommend \citet{HaiYang1998}. \citet{HaiYang2015} has demonstrated that the traffic congestion can be effectively relieved by accommodating reservation requests to the level that the highway capacity allows. More recently, \citet{HaiYang2023} proposed a novel hybrid framework integrating booking and rationing strategies on road traffic, which is formulated as a linear program to maintain the fairness, efficiency, and flexibility of individual choices. The computational results indicate that with the system-optimal integrated scheme of booking and rationing, the relative travel time reduction of each OD pair can be more than 20\%. The above two efforts have shown that trip reservations and the integration of booking and rationing are highly promising measures to relieve congestion. 

Another demand management measure is congestion pricing, typically by increasing prices during peak hours to reduce demand (i.e., surge pricing) or by decreasing prices in other periods to encourage passengers to shift their departure times \citep{YangHai2011, Xiao2015, He2017, YangHai2023}. For example, \citet{shadi2018} formulated an integrated optimization model for the elastic passenger-centric train timetabling and the pricing problem, which can be solved by a simulated annealing heuristic algorithm. Recently, \citet {HaiYang2020} proposed a rewards program combined with surge pricing, where riders pay an additional amount to a rewards account during peak periods and then use the balance in the rewards account to subsidize off-peak trips. This paper finds that in some cases, all three stakeholders (i.e., passengers, drivers, and platforms) fare better under this reward scheme combined with surge pricing. 

To sum up, the above literature specifically studied the trip booking and congestion pricing problem with the goal of alleviating congestion. The key idea of these two policies is to encourage passengers to shift their travel times, thereby improving the alignment of supply and demand, which is consistent with URT system operations. However, the aforementioned studies mainly focus on road traffic rather than train timetabling and are limited to a single congested bottleneck. 

\subsection{Paper contributions}
\label{sec:Contribution}

Overall, a wide range of mathematical models and solution methodologies have been explored to develop time-dependent passenger flow control strategies and corresponding train timetables in URT systems, as well as congestion pricing methods aimed at encouraging trip shifting in road traffic. However, most of these approaches are limited to addressing train timetabling, passenger flow control, congestion pricing (or trip shifting), and trip booking problems either separately or with only partial integration. In addition, the joint integration of train timetabling and several demand-side management strategies under operational rules that guarantee reserved passenger boarding and regulate non-reserved flow remains relatively underexplored. This integration poses nontrivial structural interactions between demand-side and supply-side decisions. Moreover, the related studies typically rely on heuristic methods to solve the proposed models. While these methods have provided valuable insights, they fall short in guaranteeing solution quality and do not achieve an integrated optimization across potential operational methods, which is essential for maximizing operational efficiency and effectiveness. Table \ref{tab:literature} outlines the characteristics of closely related literature, highlighting the contributions of our study, which are detailed as follows.

\begin{table}[]
\caption{Overview of included aspects in the discussed literature. Abbreviations: HY = Hybrid method combing heuristic and a standard slover; SA = Simulated Annealing; MA = Mathematical Analyses; HE = Heuristic; BD= Benders-decomposition based solution approach.}
\label{tab:literature}
\resizebox{\textwidth}{!}{%
\begin{tabular}{lcccccc}
\hline
Publications & Transportation systems & Trip reservation  & Trip shifting & Train timetabling & Passenger flow control  & Solution methods \\ 
\hline
\cite{Shi2018} & URT & & & \checkmark  &     Time-based         &     HY  \\
\cite{shadi2018} & Railway & &  & \checkmark  &            &  SA  \\  
\cite{Binder2021}& Railway & &  & \checkmark  &             &  SA  \\ 
\cite{Polinder2022} & Railway  & &  & \checkmark  &             & HE  \\
\cite{Leutwiler2022}& Railway  & &  & \checkmark  &             & BD  \\
\cite{Leutwiler2023}  & Railway  & &  & \checkmark            &   & BD  \\
\cite{Lu2023} & URT  & &  & \checkmark   &     Train-based            &  HY  \\
 \cite{BAO2023} & Road  & & At one bottleneck&  &             &  MA  \\
\cite{HaiYang2023} & Road  & Static & &   &           & MA  \\
\cite{HaiYang2020}& On-demand & & One-dimensional &  &              &  MA \\  
\textbf{This paper} & \textbf{URT}     &  \textbf{Time-varying}
&  \textbf{Four-dimensional}  & \checkmark   &  \textbf{Train-based}     & \textbf{BD}  \\
\hline
\end{tabular}%
}
\end{table}

(i) We develop a unified and generalized modeling framework for the complex problem that integrates trip shifting, passenger flow control, and train timetabling in urban rail transit systems with the trip reservation policy. The proposed model captures the interdependencies among these components and is designed with both generality and scalability as demonstrated by the three extensions we formulate. It incorporates a time-varying reservation slot allocation scheme, enforces train capacity constraints, and ensures service fairness. The framework can be extended to address related problems, such as incorporating the peak/off-peak pricing policy and elastic passenger demand from other transport modes. Moreover, we scale the model from a single-line setting to a network-level formulation. These extensions are built on top of our proposed model, rather than constructing a new framework from the ground up. 

(ii) We propose a novel and generalizable decomposition framework tailored to the mathematical structure of passenger-centric timetabling problems. The framework is built upon a Benders decomposition method embedded within a branch-and-cut algorithm. In our approach, partial information about passenger dynamics is incorporated into the timetabling subproblem, which significantly reduces the number of feasibility cuts. To further enhance computational performance, we introduce valid equalities and inequalities that tighten the bounds at each node and reduce the generation of inefficient feasibility cuts. We also implement heuristic acceleration strategies to improve the solution quality at the root node. This solution framework is not limited to the BDTT problem, and it can be extended to a wide range of passenger-oriented timetabling optimization problems. Computational experiments demonstrate the superiority of our solution approach, which incorporates the proposed valid inequalities. This approach is evaluated against both a hybrid algorithm that combines local search with GUROBI and an exact algorithm, as both benchmarks are developed based on the widely adopted decomposition framework for the timetabling problem integrated with passenger flow control.

(iii) We validate the proposed methodologies through both proof-of-concept and real-life case studies. The numerical results demonstrate considerable improvements in operational efficiency, service fairness, fleet size, and congestion levels. For example, our BDDT approach enables a reduction in the required fleet size while still satisfying all passenger demand, and it reduces average passengers' waiting times compared to strategies that depend only on passenger flow control in real-world case studies. These results highlight the practical value of jointly optimizing timetabling, trip shifting, and passenger flow control. On the algorithmic side, the proposed solution framework demonstrates strong computational performance. It is capable of producing high-quality and even optimal solutions directly at the root node. In real-world case studies, our algorithm outperforms GUROBI in terms of computational efficiency. Moreover, the proposed solution method consistently outperforms both the traditional Benders decomposition method and a hybrid algorithm.

\section{Problem description}
\label{sec:problemDescription}

In this section, we give a formal description of the BDTT problem for urban rail transit systems. Then, we introduce the assumptions adopted in our model. A framework of practical applications is provided in Figure \ref{fig:problem} in Appendix~\ref{sec:Application}. In this section, we delve into a more detailed discussion.

\subsection{Line structure and time discretization} 

Consider an oversaturated URT line where passengers arrive during peak periods and frequently experience congestion. The congestion results in them being detained. The set of stations along this line is represented as $\mathcal{S} = \{ 1, 2, \dots, \left| \mathcal{S} \right| \}$, with each station indexed $u$ or $v$. The section that connects stations $u$ and $u+1$ is referred to as section $u$. The set of trains is denoted as $\mathcal{I} = \{ 1, 2, \dots, \left| \mathcal{I} \right| \}$, where train $i$ should depart from the first station and terminate at station $\left| \mathcal{S} \right|$. The capacity of each train is denoted by $C^{max}$. To model the dynamic evolution of trains, the continuous time horizon is discretized into timestamps of duration $\sigma$. The set $\widetilde{\mathcal{T}} = \Big\{ \tilde{t} \mid 1, 2, \dots, \left| \widetilde{\mathcal{T}} \right| \Big\}$ corresponds to these discretized timestamps.

\subsection{Time-dependent passenger demand and time equivalization}

Given the aforementioned notations, the time-dependent Origin-Destination (OD) demand for those arriving at station $u$ at timestamp $\tilde{t}$ and heading to station $v$ is denoted as $D_{uv\tilde{t}}$, and the predetermined reservations for the OD pair from stations $u$ to $v$ at timestamp $\tilde{t}$ is denoted as $\hat{D}_{uv\tilde{t}}$. Note that this problem encompasses four indexes: time, trains, the origin, and the destination of each passenger. To effectively reduce the dimensionality and scale down the complexity of the problem, we adjust the time dimension at the second station and all subsequent stations. Under this adjusted time dimension, a separate index for the arrival and departure times of a train at each station is unnecessary. Instead, the departure time from the first station can be used to derive the arrival and departure times at subsequent stations, thereby reducing the problem's dimensionality. We define this adjusted time dimension as \textit{equivalent time} and denote the studied equivalent time horizon as $\mathcal{T}=\{ t \mid 1,2,\ldots, \left| \mathcal{T} \right|\}$, where $t$ is the index of the \textit{discretized equivalent timestamp}. As for portraying time-dependent passenger demand, its arrival time \( \tilde{t} \) is mapped to the corresponding equivalent timestamp \( t \) through this skewing operation of the time dimension. Besides, the subset of equivalent timestamps during peak hours is denoted as \( \hat{\mathcal{T}} \), which is included within \( \mathcal{T} \), i.e., \( \hat{\mathcal{T}} \subseteq \mathcal{T} \). 

\begin{example}

\textit{To facilitate understanding of the concept of equivalent time, Figure~\ref{fig:time} visually depicts the adjustment of the time dimension in the proposed modeling approach. Consider a URT line involving four stations along a line, with two trains operating from stations A to D. The first train originally departs from stations A, B, and C at timestamps 1, 3, and 5, respectively. After applying the proposed adjustment of the time dimension, the departure times of the first train at all stations are adjusted to align with the first equivalent timestamp. Similarly, the second train leaves each station at the third equivalent timestamp.}

    \begin{figure}[]
     \centering
     \begin{subfigure}{0.49\textwidth}
         \centering
 \includegraphics[width=\textwidth]{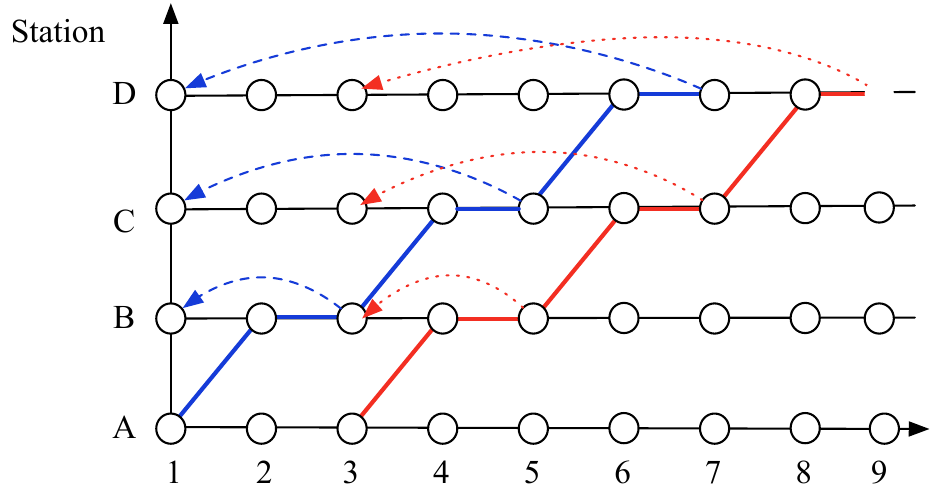}
         \caption{Timetable under the original timestamp}
     \end{subfigure}
     \begin{subfigure}{0.49\textwidth}
         \centering
\includegraphics[width=\textwidth]{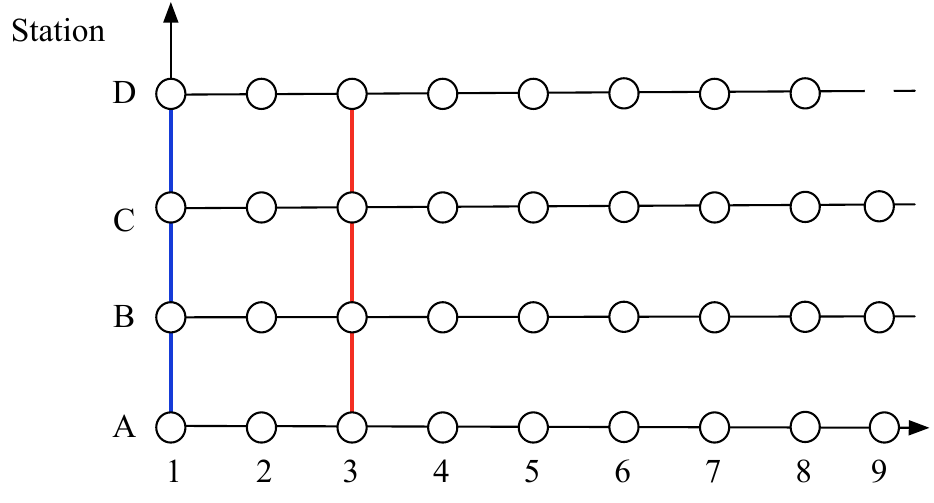}
         \caption{Timetable under the equivalent timestamp}
     \end{subfigure}
        \caption{Illustration of timetables before and after the adjustment of the time dimension.}
        \label{fig:time}
\end{figure}
\end{example}

\subsection{Categorization of passengers and demand-side management strategies }

\begin{figure}
      \centering
\includegraphics[width=\linewidth]{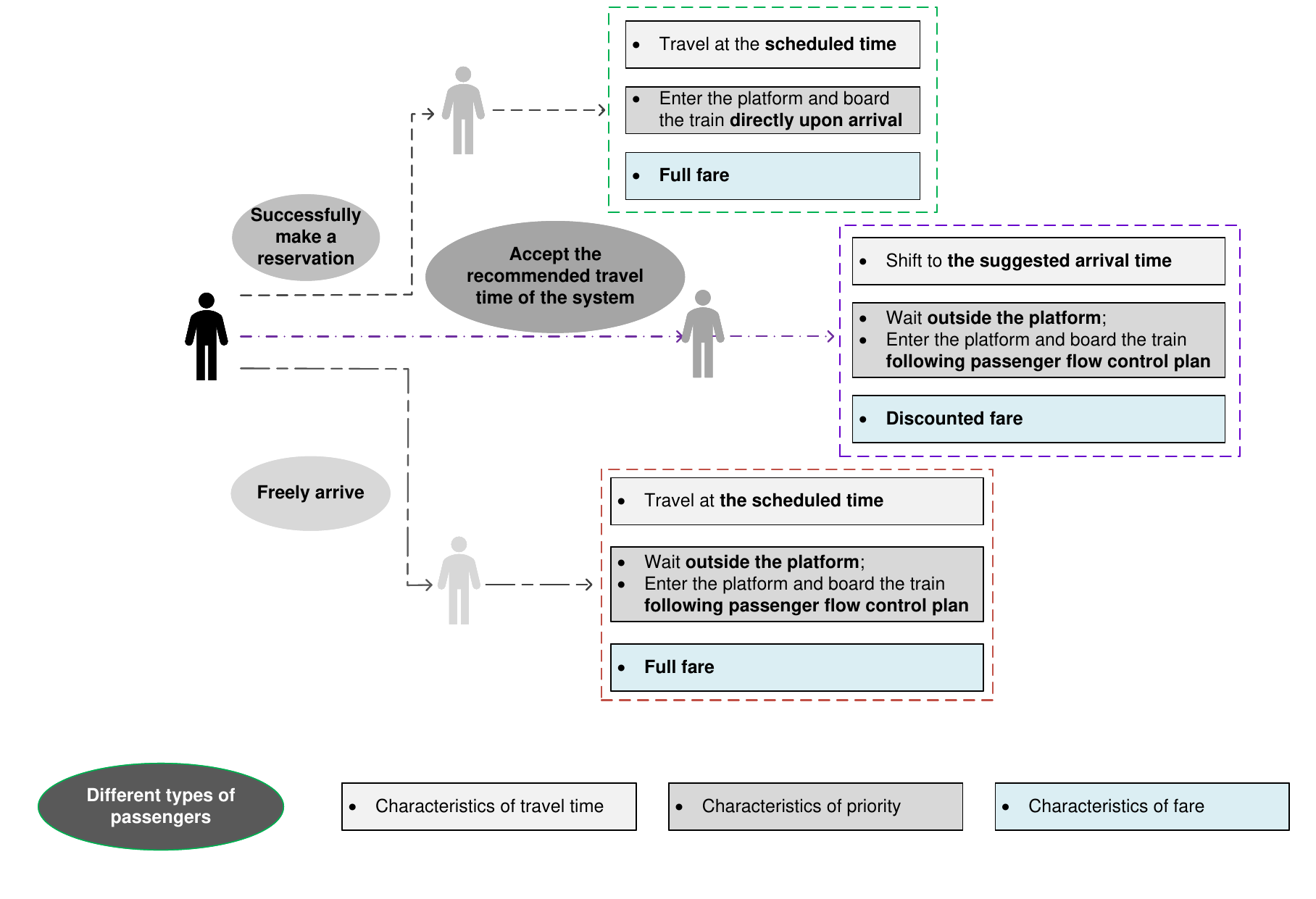}
      \caption{Illustration of passengers' types, travel time, priority, and fare characteristics. \\
      Notes: Passengers with reservations are never detained.}
      \label{fig:characteristics}
  \end{figure}

With advancements in emerging technologies, operators can now utilize intelligent trip reservation systems to implement demand-side management strategies, as adopted by the Beijing metro. These strategies include trip reservations, incentives to encourage off-peak travel, and passenger flow control. Based on these demand-side management methods, passengers can be categorized into three groups: those with reservations, those who attempt to make reservations but are unsuccessful due to limited availability and follow the system's recommended travel times, and those who arrive freely according to their own schedules. The travel time, boarding priority, and fare characteristics for these three passenger types are summarized in Figure~\ref{fig:characteristics}. Next, we introduce these three types of passengers and the implemented demand-side management strategies in detail.

\textbf{(i) Passengers with reservations.} Let $\hat{D}_{uvt}$ denote the number of reserved passengers arriving at station $u$ and traveling to station $v$ at time $t$. These passengers successfully secure reservations for their preferred travel times, granting them the boarding priority. They arrive at their origins at scheduled times, enter to the platform through dedicated reservation gates, and board the first available train. By paying the full ticket fare, they enjoy seamless boarding upon arrival. Passengers with reservations are not subject to passenger flow control measures and are never stranded.

\textbf{(ii) Passengers accepting the system's recommended time.} For users of the intelligent reservation system who attempt to make a reservation but cannot secure a slot due to availability limitations, the system offers a suggested arrival time. If the passenger agrees to travel at this recommended time, they receive an incentive. Upon arrival at stations, these passengers are required to queue outside the platform and wait for permission to enter and board a train. In this paper, we consider fare discounts (denoted as $\phi$) as the incentive, though other benefits such as credits or cashback could also be offered. By adjusting their scheduled travel time, these passengers gain financial savings and potentially a more comfortable journey with fewer in-vehicle passengers. We define the decision variable $\kappa_{uvt't}$ to represent the number of passengers traveling from station $u$ to station $v$ who shift their arrival time from $t$ to $t'$. This variable will be determined in the optimization model from a systematic perspective.

\textbf{(iii) Passengers who arrive freely.} These passengers do not use the reservation system or follow its recommendations, arriving instead at their originally planned travel time. Upon arrival, they must queue outside the platform and wait for permission to enter and board trains according to the passenger flow control plan. Unlike passengers who accept the system's recommended time, freely arriving passengers do not receive any incentives. The passenger flow control plan specifies the number of passengers who are allowed to enter the platform to board the train when a train arrives at a station. This allocation includes both passengers who accept the recommended time and those who arrive freely. We denote the passenger flow control decision variable as $b_{iuv}$, representing the number of passengers traveling to station $v$ who are permitted to enter the platform to board the train when train $i$ arrives at station $u$.

\begin{remark}
In this paper, we propose a reservation system that does not require passengers to book specific seats, but rather allows them to secure entry to platforms and immediate board trains by reserving a time slot. Seat booking is a common demand management strategy in high-speed rail systems with the large headway, such as China's and Europe's (e.g., Germany's ICE trains and the Eurostar from London to Pairs). However, due to the limited seating and high passenger volume during peak hours, metro systems often accommodate standing passengers. For instance, during morning rush hours, the Beijing metro system has more standing than seated passengers. Moreover, in congested metro networks (like the Beijing metro), the headway during peak hours can be as short as 90 seconds. Therefore, it is practical for passengers to simply reserve a time slot, ensuring they can board immediately without the risk of delays and only wait for a maximum of one headway on the platform. In addition, it is worth mentioning that platforms will not be congested, as only passengers with reservations can wait on the platform, and their number is not excessively high.
\end{remark}

\begin{remark}
In real-world operations, the proposed demand-side management strategies, including trip reservations and passenger flow control, have been piloted in the Beijing metro system \citep{BMCT}. However, widespread implementation may still face practical challenges. These strategies require technological infrastructure (e.g., digital platforms for making reservations and physical facilities for managing non-reserved passenger queues), as well as additional operational coordination. In practice, Beijing metro enables reservations through an official WeChat mini-program, which improves accessibility but still entails communication and system management costs. For non-reserved passengers, physical queuing areas have been installed outside station entrances to regulate access during peak hours. 

Currently, these strategies are applied only during peak periods, when travel demand regularly exceeds available capacity. The reservation strategy in peak hours mainly targets regular commuters, whose travel times are predictable and typically involve a single trip in each direction per day. This minimizes the need for multiple reservations. During off-peak hours, the system operates without reservations, as train supply is generally sufficient. This targeted application reflects practical trade-offs and highlights the importance of adaptive, context-specific strategies in real-world rail transit systems.
\end{remark}

\subsection{Interactions among the government, the operator, and passengers}

\begin{figure}
    \centering
\includegraphics[width=0.85\linewidth]{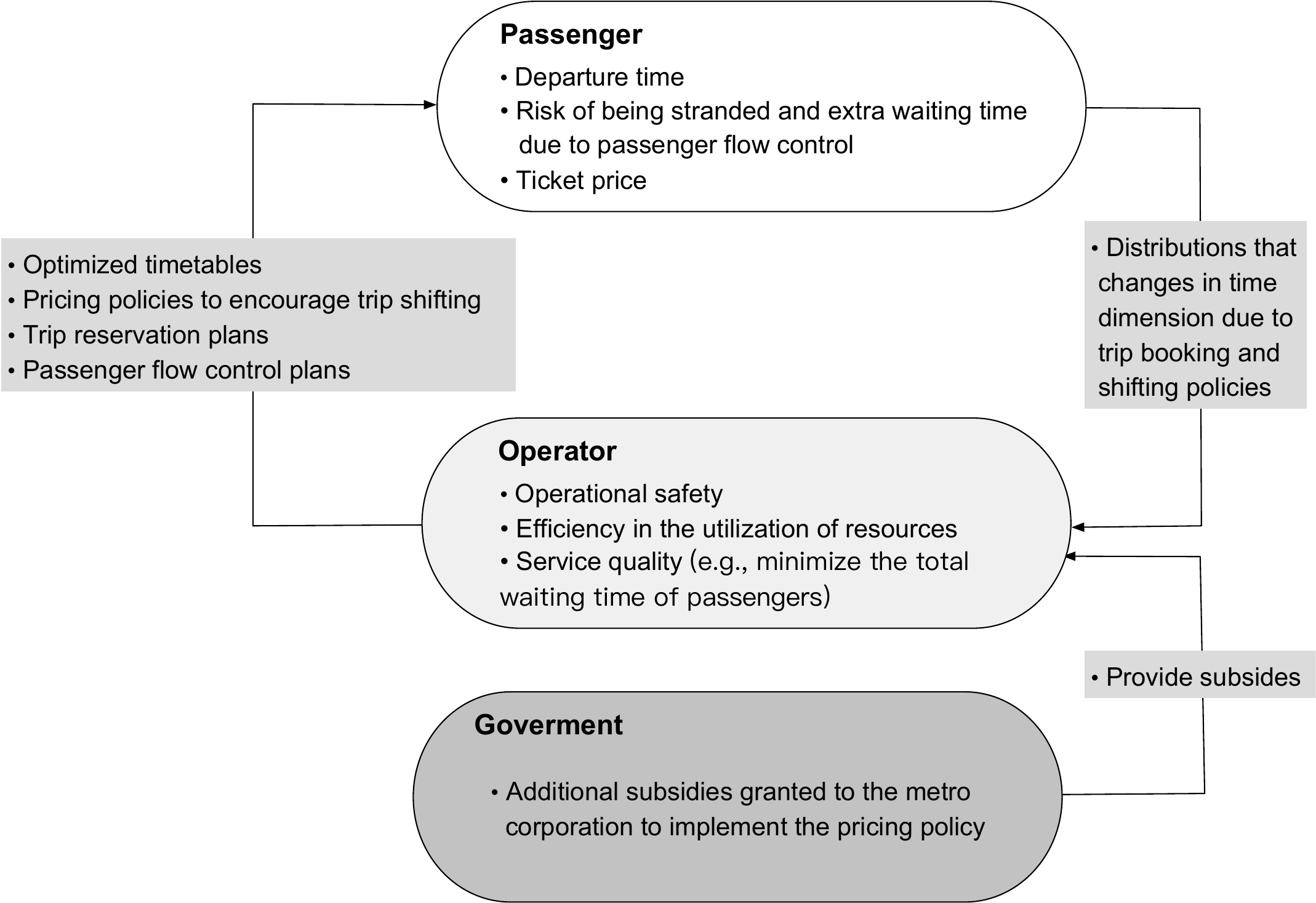}
    \caption{Interactions among the government, passengers, and the operator.}
    \label{fig:interaction}
\end{figure}

In summary, we consider the interactions among three key stakeholders: the government, metro operators, and passengers, as illustrated in Figure~\ref{fig:interaction}. In a megacity with an overcrowded metro system, the government plays a crucial role by providing subsidies to the metro corporation to promote pricing policies that distribute passenger demand more evenly over time, thereby enhancing service quality and improving passenger well-being. While the government ultimately aims to minimize its subsidy expenditures, it does not solve the optimization problem directly. Instead, the optimization problem is formulated and solved from the perspective of operators, who seek to balance service quality (e.g., passengers' waiting time) and subsidy levels. Operators optimize both supply- and demand-side strategies to improve system performance. On the supply side, train timetables are optimized in response to time-dependent passenger demand, particularly during peak periods. On the demand side, the operator employs measures such as pricing incentives for trip shifting, a trip-booking mechanism, and passenger flow control. Although such incentives may reduce fare revenues, these losses are compensated through government subsidies.

The operators' role is to determine optimal operational strategies under varying levels of government subsidy. Based on this, the operator provides the government with a set of feasible and efficient solutions, showing how different subsidy levels influence passenger waiting times and overall service quality. This enables the government to evaluate the cost-effectiveness of potential subsidy policies and make informed decisions, rather than directly minimizing subsidies within the optimization model.

In our study, the interaction between the government and the operators is implicitly modeled such that the government provides subsidies to influence the temporal distribution of passenger demand. Operators then respond by optimizing their timetables based on this adjusted passenger flow distribution. Thus, these two stakeholders are linked through time-varying passenger demand and the trip-shifting policy. The interaction between the operator and passengers is explicitly captured by modeling the dynamics of train operations and passenger movements, with timetables being optimized based on the passenger flows under the trip-shifting policy. Besides, the interaction between the government and passengers is represented through the trip-shifting decision variable and corresponding constraints, which govern how passengers respond to the incentives provided.

Throughout this study, we assume that the objective is to minimize a combination of total passengers' waiting time and government subsidies. Since the government bears the cost of trip-shifting incentives, the metro operators' profit remains unaffected. Under these conditions, the operators are motivated to improve service quality by reducing total waiting times. In the long term, such improvements can enhance the attractiveness of the metro system, potentially leading to increased ridership and better overall performance.

In addition, we make the following assumptions to rigorously formulate the models for the investigated problem.

(\roman{mycounter}\stepcounter{mycounter}) We assume that information about reservation slots is released to the public through the mobile application. The number of reservations is limited and follows a first-come, first-served principle. This assumption is aligned with the practical implementation by the Beijing Metro \citep{BMCT}.

(\roman{mycounter}\stepcounter{mycounter}) We focus on weekday scenarios, as commuting behavior is generally more predictable on these days. Thus, we assume that passengers without a reservation can empirically anticipate the level of congestion they would face and the risk of being stranded if they do not shift their arrival times. 

(\roman{mycounter}\stepcounter{mycounter}) We assume that all reservation slots are selected daily and that reserved passengers are guaranteed to show up on time. Observations from the Beijing metro's trip reservation show that most reserved passengers are commuters during peak hours, and reservation slots are quickly filled within minutes of release.

(\roman{mycounter}\stepcounter{mycounter}) Passenger demand is assumed to be homogeneous and individual behavior is outside the scope of our study. Similar assumption has been widely adopted in
recent previous studies \citep[e.g.,][]{YIN2023, bib0099}. We aim to introduce a proof-of-concept methodology that allows operators and governments to evaluate the effectiveness of integrating demand management strategies with timetabling.

\section{Mathematical formulation}
\label{sec:model}
In this section, we first provide the proposed mathematical model with its associated notations. We then propose three extensions of the proposed model that consider the peak/off-peak pricing strategy, the elastic passenger demand from the other transportation modes, and scale up to the network level.

\subsection{Notations used in the formulation}

To model the BDTT problem, we introduce three families of decision variables, as listed in Table \ref{table:Decision}. Specifically, the integer variable $\kappa_{uvt't}$ denotes the number of non-reserved passengers who shift their arrival time at the origin station from $t$ to $t'$. The second one is the number of non-reserved passengers who are allowed to board train $i$ to reach their destination $v$, which is denoted as $b_{iuv}$. The third set of variables $z_{it}$ aims to model the dynamics of trains. The parameters, dependent variables and abbreviations used throughout this paper are summarized in Table \ref{table-para}. We use the following formulation for the BDTT problem, including upper and lower limit constraints on headway, strict capacity constraints, fairness-preserving constraints, and so on.
\begin{table}[H]
  \centering
  \vskip 0.25cm
  \caption{Decision variables in the BDTT model.} \label{table:Decision}
   \resizebox{\textwidth}{!}{%
  \begin{tabular}{lll}
\hline\noalign{\smallskip}
Symbol & Definition & Type \\
\hline
$\kappa_{uvt't}$ & Number of non-reserved passengers travelling from stations $u$ to $v$ who shift their & Integer \\
&  arrival time from timestamp $t$ to $t'$ & \\
$b_{iuv}$ & Number of non-reserved passengers heading to station $v$ who are allowed to board  & Integer \\
 & train $i$ at station $u$\\
$z_{it}$ & If train $i$ has not departed at timestamp $t$, $z_{it}=1$; otherwise, $z_{it}=0$ & Binary\\
\hline
\end{tabular}%
}
\end{table}

\begin{table}[]
  \centering
  \vskip 0.25cm
  \caption{Parameters, dependent variables, and abbreviations.} \label{table-para}
  \begin{tabular}{ll}
\hline\noalign{\smallskip}
\multicolumn{2}{l}{\bf{Sets}}\\
$\mathcal{I}$& Set of trains, $\mathcal{I}=\{1, 2, \ldots, \left| \mathcal{I} \right| \},$ indexed by $i$, $j$\\
$\mathcal{S}$& Set of stations, $\mathcal{S}=\{1,2, \ldots, \left| \mathcal{S} \right|\},$ indexed by $u$, $v$, $m$\\
$\mathcal{S}_{u+1}$& Set of stations following station $u$, $\mathcal{S}_{u+1}=\{u+1, \ldots, \left| \mathcal{S} \right|\}$\\
$\widetilde{\mathcal{T}}$& Set of discretized timestamps, $\widetilde{\mathcal{T}}=\{1, 2, \ldots, \left| \widetilde{\mathcal{T}} \right|\},$ indexed by $\tilde{t}$\\
$\mathcal{T}$& Set of discretized equivalent timestamps, $\mathcal{T}=\{1,2,\ldots,  \left| \mathcal{T} \right|\},$ indexed by $t$\\
$\hat{\mathcal{T}}$& Set of discretized equivalent timestamps during peaking hours, $\hat{\mathcal{T}} \subseteq \mathcal{T}$\\
\multicolumn{2}{l}{\bf{Parameters}}\\
$\sigma$& Length between two timestamps\\
$s_{u}$ & Train running time on the section between stations $u$ and $u+1$\\
$h^{min}$ & Minimum headway\\
$h^{max}$ & Maximum headway\\
$C^{max}$& Train capacity\\
$\varepsilon_{uv}$ & Ticket price from stations $u$ to $v$\\
$\phi$ & Discount rate of ticket prices\\
$\varrho_{iu}$ & Service fairness factor of train $i$ at station $u$\\
$\imath$& Maximum timestamps that unreserved passengers can shift their trips\\
$D_{uvt}$& Number of unreserved passengers who arrive at station $u$ and head to station $v$ \\
& at timestamp $t$\\
$\hat{D}_{uvt}$& Number of reserved passengers who arrive at station $u$ and head to station $v$ \\
& at timestamp $t$\\
$\omega_t$, $\omega_s$ & Weighting coefficients \\
\multicolumn{2}{l}{\bf{Involved variables}}\\
$x_{it}$ &  Binary variable. $x_{it}=1$ if timestamp $t$ belongs to the headway between train $i-1$ \\
 & and train $i$\\
$d_{i}$&  Departure time of train $i$\\
$h_{i}$&  Headway between trains $i-1$ and $i$, defined as the timestamp from the departure\\
&  time of train \(i-1\) to the timestamp just before the departure of train \(i\)\\
$o_{iu} (\hat{o}_{iu})$& Number of on-board passengers without (with) reservations in train $i$ when it \\
& departs from station $u$ \\
$l_{iu} (\hat{l}_{iu})$& Number of passengers alighting from train $i$ without (with) reservations when it \\
& arrives at station $k$ \\
$w_{iu}$& Number of passengers waiting for train $i$ at station $u$ \\
$w_{iuv}$& Number of passengers waiting for train $i$ at station $k$ and head to station $v$ \\
$r_{iuv}$& Number of passengers detained by train $i$ at station $u$ and head to station $v$\\
$\hat{b}_{iuv}$ & Number of passengers with reservations who are allowed to board train $i$ at \\
& station $u$ and head to station $v$\\
$\hat{p}^{w}_{iut}$ & Number of passengers with reservations who arrive at timestamp $t$ and wait\\
& for train $i$ at station $u$\\
$\hat{p}^{wc}_{iut}$ & The cumulative number of passengers with reservations waiting for train $i$ at\\
&station $u$ and time $t$ \\
$p^{w}_{iut}$ & Number of passengers without reservations who arrive at timestamp $t$ and wait\\
& for train $i$ at station $u$\\
$p^{wc}_{iut}$ & The cumulative number of passengers without reservations waiting for train $i$ at \\
& station $u$ and time $t$ \\
$r_{iuv}$ & Number of passengers without reservations who are detained by train $i$ at station $u$ \\
& and head to station $v$\\
\multicolumn{2}{l}{\bf{Abbreviations}}\\
$F^{t}$ & Total waiting time of passengers\\
$F^{s}$ & Total additional government subsidies that arise from incentives to shift trips\\
\hline
\end{tabular}
\end{table}

\subsection{INLP model for the BDTT problem}

We first introduce the objective function, followed by the constraints. The constraints are formulated to model the waiting times of passengers with and without reservations, train dynamics, interactions between trains and different types of passengers, and the domains of the decision variables.

\textbf{Objective function.} The objective function~\eqref{eq:objective} minimizes the weighted sum of the passengers' waiting time (denoted as $F^{t}$), and additional government subsidies due to encouraging passengers to shift their departure times (denoted as $F^{s}$), where $\omega_t$ and $\omega_s$ represent weighting coefficients, respectively.
\begin{align}
   \min & \quad \omega_t F^{t}+ \omega_s F^{s} \label{eq:objective} \\
  & \quad F^t =\sigma \Big[\sum_{i\in\mathcal{I}}\sum_{u\in\mathcal{S}}\sum_{t\in\mathcal{T}}(\hat{p}^{wc}_{iut}+p^{wc}_{iut}) + \sum\limits_{i \in \mathcal{I}}\sum_{u\in\mathcal{S}}\sum\limits_{t \in \mathcal{T}}(x_{it}\sum\limits_{v \in \mathcal{S}_{u+1}}r_{iuv})\Big], \label{eq:objWaiting} \\
  &\quad F^{s} = \sum\limits_{u\in \mathcal{S}}\sum\limits_{v \in \mathcal{S}_{u+1}}\sum\limits_{t \in \mathcal{T}} \Big[D_{uvt} \varepsilon_{uv} - \varepsilon_{uv}[\sum\limits_{t+1\leq t' \leq \min\{\left|\mathcal{T} \right|, t +\imath\}}(\phi\kappa_{uvtt'}+\kappa_{uvtt})] \Big]. \label{eq:objFare}
\end{align}
Constraint~\eqref{eq:objWaiting} is formulated to calculate the total waiting time of passengers, including both the waiting time and detained time. Specifically, $\hat{p}^{wc}_{iut}$ represents the number of passengers who arrive during the headway between trains $i-1$ and $i$, and wait for train $i$ at station $u$ at timestamp $t$. The variable $p^{wc}_{iut}$ also counts passengers arriving in the same period and waiting at the same station and time, but these passengers do not have reservations. The term $r_{iuv}$ denotes the number of non-reserved passengers who are detained by train $i$ at station $u$ and are heading to station $v$. Constraint~\eqref{eq:objFare} calculates the additional government subsidies, which are defined as the difference in fare revenue before and after applying discounts to incentivize shifts of arrival times. Additionally, we use the terms \textit{additional subsidies} and \textit{lost revenue} interchangeably. The variables involved in constraints~\eqref{eq:objWaiting} and \eqref{eq:objFare} are introduced in detail below.

\begin{figure}
  \centering  \includegraphics[width=0.9\textwidth]{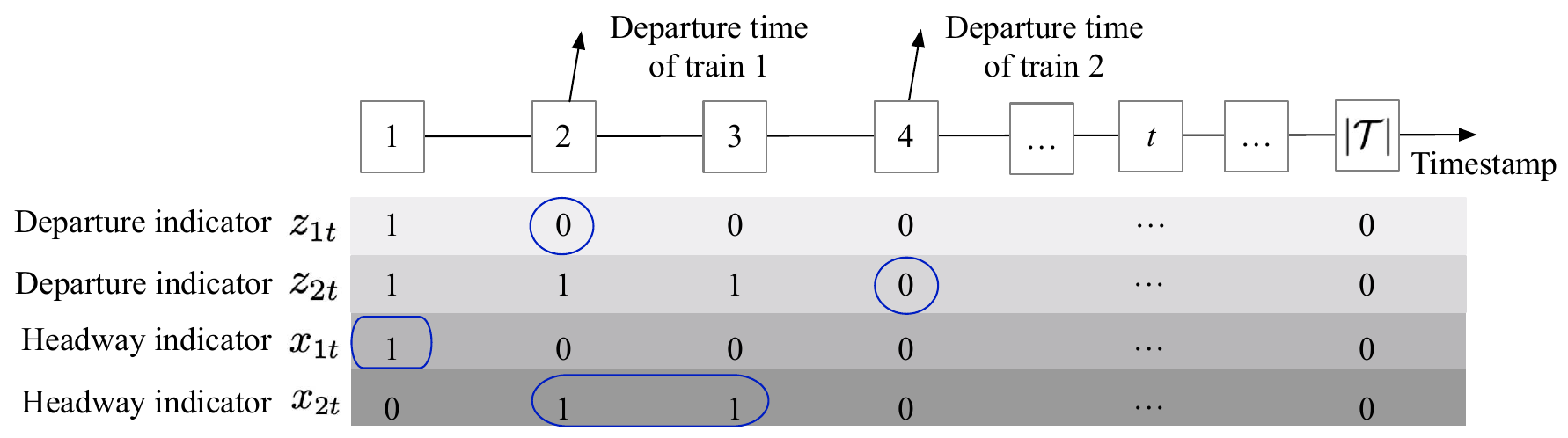}
  \caption{Illustration of departure and headway indicators.}\label{fig:Auxiliary}
\end{figure}

\textbf{Waiting time of passengers with and without reservations.}
To compute the waiting time of passengers, we first propose the binary variable $x_{it}$ to represent the headway indicator, which is coupled with the departure indicator $z_{it}$ in the following constraints. Here, the \textit{headway} between trains $i-1$ and $i$ is defined as the timestamp from the departure time of train \(i-1\) to the timestamp just before the departure of train \(i\). To facilitate understanding, an illustrative example is visualized in Figure~\ref{fig:Auxiliary}. It can be seen that two trains depart at the second and fourth timestamps, respectively. Hence, timestamp 1 corresponds to the first headway, while timestamps 2 and 3 belong to the second headway between trains 1 and 2.

\begin{figure}[h]
     \centering
     \begin{subfigure}{\textwidth}
         \centering
 \includegraphics[width=\textwidth]{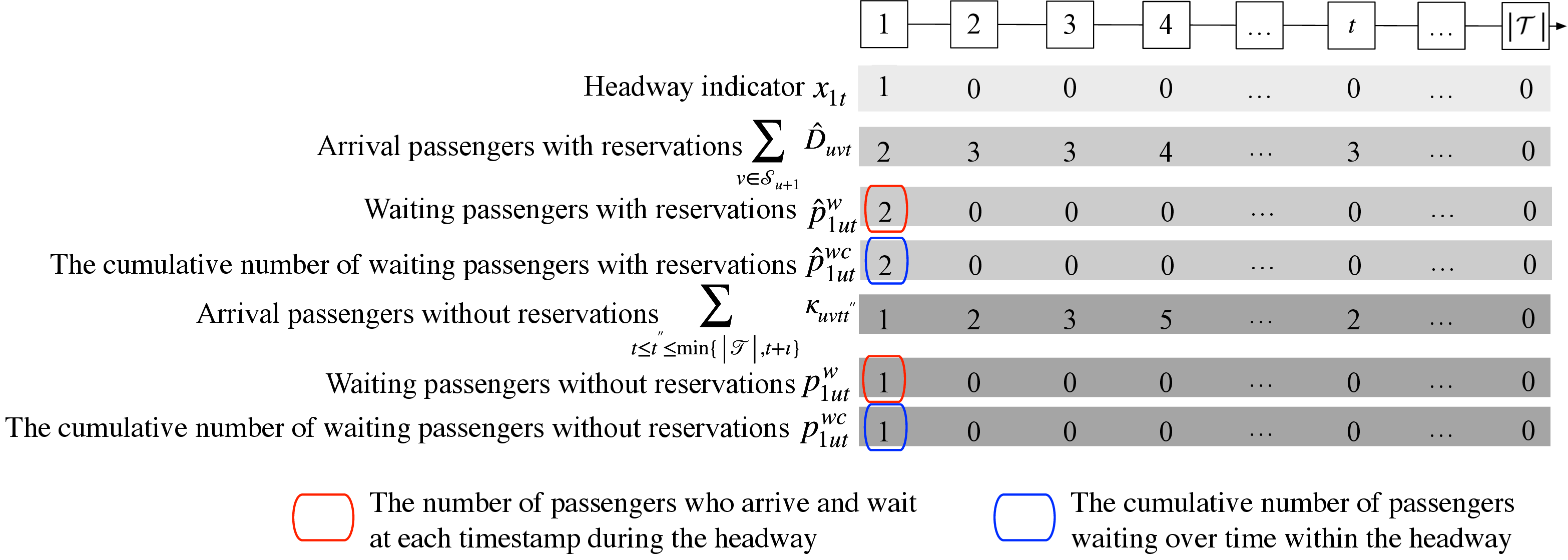}
         \caption{Passengers waiting for the first train}
     \end{subfigure}

\vspace{1em}

     \begin{subfigure}{\textwidth}
         \centering
\includegraphics[width=\textwidth]{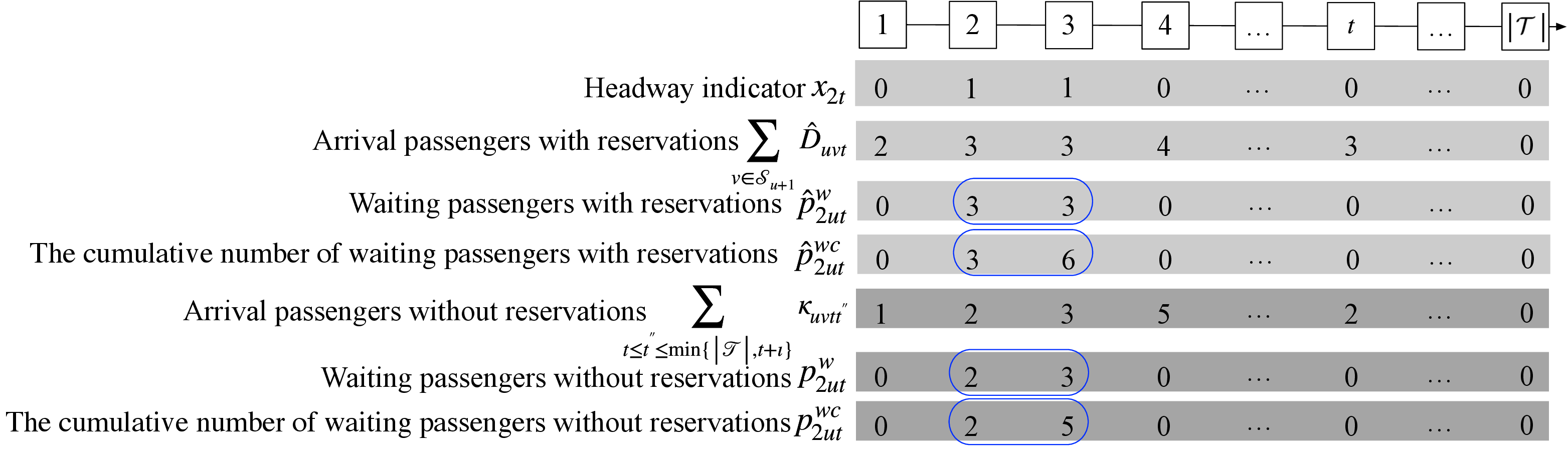}
         \caption{Passengers waiting for the second train}
     \end{subfigure}
        \caption{Number of waiting passengers within each headway.}
        \label{fig:waiting}
\end{figure}

Motivated by \cite{Xia2023}, constraints~\eqref{eq:waiting1} - \eqref{eq:waiting2} are proposed to compute the number of passengers with and without reservations who newly arrive at station $u$ and timestamp $t$ and wait for train $i$, and constraints~\eqref{eq:waiting3} - \eqref{eq:waiting4} track the cumulative number of passengers who arrive during the headway between trains $i-1$ and $i$, and wait for train $i$ at station $u$ at timestamp $t$. Figure \ref{fig:waiting} illustrates the formulations of these dynamics. Train 1 departs at timestamp 2, with two reserved and one non-reserved passengers arriving at timestamp 1. According to constraints \eqref{eq:waiting1}, two reserved passengers wait at timestamp 1, and constraints \eqref{eq:waiting3} reflect a cumulative total of one non-reserved passenger. According to constraints \eqref{eq:waiting2} and \eqref{eq:waiting4}, there is one waiting passenger without reservations. Similarly, as shown in Figure \ref{fig:waiting}(b), train 2 departs at timestamp 4, with timestamps 2 and 3 falling within the second headway. By timestamp 3, a cumulative total of 6 reserved and 5 non-reserved passengers who arrived within this headway are waiting for train 2.

\begin{align}
    \hat{p}^w_{iut} & = x_{it}\sum_{v\in\mathcal{S}_{u+1}}\hat{D}_{uvt} & \forall i \in \mathcal{I}, u \in \mathcal{S},t \in \mathcal{T}, \label{eq:waiting1}\\
    p^w_{iut}& = x_{it}\sum_{v\in\mathcal{S}_{u+1}}\sum_{t\leq t^{''}\leq \min\{ \left|\mathcal{T} \right|, t +\imath \}}\kappa_{uvtt^{''}} & \forall i \in \mathcal{I}, u \in \mathcal{S},t \in \mathcal{T}, \label{eq:waiting2}\\
    \hat{p}^{wc}_{iut} & = x_{it}\sum_{t' \in \mathcal{T}, t'\leq t}\hat{p}^w_{iut'} &  \forall i \in \mathcal{I}, u \in \mathcal{S},t \in \mathcal{T}, \label{eq:waiting3}\\
    p^{wc}_{iut} &= x_{it}\sum_{t' \in \mathcal{T}, t'\leq t}p^w_{iut'} & \forall i \in \mathcal{I}, u \in \mathcal{S},t \in \mathcal{T} \label{eq:waiting4}.
\end{align}
\textbf{Dynamics of trains.} Constraints (\ref{eq:0-1timetableConstraint}) impose restrictions on the binary indicator variables associated with train operations following \cite{Lu2023}. Constraints (\ref{eq:0-1timetablefeasibleConstraint}) ensure that all trains are operated before the end of the study time horizon. Constraints (\ref{eq:0-1DepartureConstraint}) link the binary indicator variable $z_i$ with the real-valued departure time $d_i$ of train $i$. Constraints (\ref{eq:headwayConstraint}) define the headway $h_i$ between two successive trains $i-1$ and $i$. Constraints (\ref{eq:headwayLimitConstraint}) enforce both lower and upper bounds on the headway, denoted by $h^{\min}$ and $h^{\max}$, respectively. Constraints (\ref{eq:0-1headwayConstraint}) introduce a binary variable $x_{it}$, which equals 1 if and only if timestamp $t$ falls within the headway interval of train $i$, and 0 otherwise. Figure \ref{fig:Auxiliary} illustrates the relationship between $z_{it}$ and $x_{it}$. For example, suppose trains 1 and 2 depart at timestamps 2 and 4, respectively. Then, $z_{12} = 1$ and $z_{24} = 1$. The time periods between the departures of train 1 and train 2, i.e., timestamps 2 and 3 fall within the second headway. Therefore, $x_{22} = 1$ and $x_{23} = 1$.
\begin{align}
    z_{i(t+1)} & \leq z_{it} &  \forall i \in \mathcal{I}, t \in \mathcal{T} \backslash \{ \left|\mathcal{T} \right| \}, \label{eq:0-1timetableConstraint}\\
    z_{i\left|\mathcal{T}\right|}  & = 0 & \forall i \in \mathcal{I}, \label{eq:0-1timetablefeasibleConstraint}\\
    d_{i} & = \sum\limits_{t \in \mathcal{T} \backslash \{ 1 \}} [t (z_{i(t-1)} - z_{it})] &\forall i \in \mathcal{I}, \label{eq:0-1DepartureConstraint}\\
    h_{i}  &= \begin{cases}
d_{i} - \sigma & \text{if} \ i=1\\
d_{i} - d_{i-1} & \text{if} \ i \in \mathcal{I} \backslash \{ 1 \}\\
\end{cases}, \label{eq:headwayConstraint}\\
h^{min} & \leq h_{i} \leq h^{max}  & \forall i \in \mathcal{I} \backslash \{ 1 \},  \label{eq:headwayLimitConstraint} \\
x_{it} & = \begin{cases}
z_{it} & \text{if} \ i=1\\
z_{it}-z_{(i-1)t} & \text{if} \ i \in \mathcal{I} \backslash \{1\}\\
\end{cases} & \forall i \in \mathcal{I}, t \in \mathcal{T} \label{eq:0-1headwayConstraint}.
\end{align} 

\textbf{Interactions between trains and different types of passengers.} Due to trip-shifting policies, non-reserved passengers who originally plan to arrive at origin stations during peak hours are allowed to adjust their arrival times. The following constraints ensure that the total number of passengers remains consistent before and after shifting trips.

Constraint (\ref{eq:nobookingshiftConstraint}) ensures that for each timestamp $t \in \hat{\mathcal{T}}$ during peak hours, the number of passengers originally scheduled to arrive at timestamp $t$ from station $u$ to $v$ equals the sum of those who shift their arrival to timestamps $t' \in [\max\{0, t-\imath\}, t ]$. Here, the parameter $\imath$ represents the maximum timestamps that unreserved passengers can shift their trips, and $\kappa_{uvt't}$ denotes the number of such passengers shifting their scheduled arrival time from $t$ to $t'$. These constraints ensure that no passengers are lost or artificially generated due to trip-shifting policies. Constraint (\ref{eq:nobookingnonpeakingConstraint}) ensures that all passengers scheduled to arrive during non-peak hours (i.e., $t \in \mathcal{T} \backslash \hat{\mathcal{T}}$) are assumed to arrive exactly at their original arrival time $t$, as no shifting is permitted.
\begin{align}
    \sum\limits_{\max\{0, t-\imath\} \leq t' \leq t }\kappa_{uvt't} & = D_{uvt} & \forall u \in \mathcal{S}, v \in \mathcal{S}_{u+1}, t \in \hat{\mathcal{T}}, \label{eq:nobookingshiftConstraint}\\
    \kappa_{uvt't} & =
\begin{cases}
 D_{uvt}& \text{if} \ t'=t\\
0 & \text{otherwise} \\
\end{cases} & \forall u \in \mathcal{S}, v \in \mathcal{S}_{u+1}, t \in \mathcal{T} \backslash \hat{\mathcal{T}}. \label{eq:nobookingnonpeakingConstraint}
\end{align}

Thereafter, constraints (\ref{eq:bookingBoardingConstraint}) are formulated to ensure passengers with reservations, who are not controlled by the passenger flow control strategies, can board the first coming train. 
\begin{align}
    \hat{b}_{iuv} &= \sum\limits_{t \in \mathcal{T}} \hat{D}_{uvt} x_{it} &\forall u \in \mathcal{S}, v \in \mathcal{S}_{u+1}.\label{eq:bookingBoardingConstraint}
\end{align}

Constraints (\ref{eq:nobookingFullservedConstraint}) guarantee that all passengers without reservations are served.  Constraints~\eqref{eq:waitingConstraint} compute the number of non-reserved passengers waiting for train $i$ at station $u$ and heading to station $v$, denoted by $w_{iuv}$. For the first train (i.e., $i = 1$), all passengers who have arrived by that time are waiting, as no earlier trains are available. In this case, $w_{1uv}$ is equal to the total number of non-reserved passengers who arrive at station $u$ and head to $v$ before the train departs. For subsequent trains, the waiting passengers are calculated by subtracting the number of passengers who have already boarded previous trains from the cumulative arrivals. Constraint~\eqref{eq:nobookingBoardingLimitConstraint1} then ensures that the number of boarding passengers $b_{iuv}$ does not exceed the number of waiting passengers $w_{iuv}$. In addition, constraints (\ref{eq:nobookingBoardingLimitConstraint1}) also require that at least $\varrho_{iuv}$ percent of passengers without reservations at station $u$ must be served by train $i$ to reach destination $v$. This set of constraints aims to ensure twofold above-mentioned service fairness. 

Constraints~\eqref{eq:nobookingDetainedConstraint1} compute the number of non-reserved passengers who are detained by train $i$ at station $u$, denoted by $r_{iuv}$. These are passengers who wait for train $i$ but are unable to board. This value is calculated as the difference between the number of waiting passengers and those who successfully boarded the train.
\begin{align}
\sum\limits_{i \in \mathcal{I}} b_{iuv} &  = \sum\limits_{t \in \mathcal{T}} D_{uvt} &\forall u \in \mathcal{S}, v \in \mathcal{S}_{u+1},  \label{eq:nobookingFullservedConstraint} \\
    w_{iuv} &=
\begin{cases}
\sum\limits_{t' \in \mathcal{T}}z_{it'}\sum\limits_{t' \leq t \leq \min\{ \left|\mathcal{T} \right|, t' +\imath \} }\kappa_{uvt't} & \text{if} \ i=1\\
\sum\limits_{t' \in \mathcal{T}}z_{it'}\sum\limits_{t' \leq t \leq \min\{ \left|\mathcal{T} \right|, t' +\imath \}  }\kappa_{uvt't} - \sum\limits_{j = 1}^{i-1} b_{juv} & \text{if} \ i \in \mathcal{I} \backslash \{1\}\\
\end{cases} &\forall u \in \mathcal{S}, v \in \mathcal{S}_{u+1},\label{eq:waitingConstraint} \\
\varrho_{iuv} w_{iuv} & \leq b_{iuv} \leq w_{iuv}  & \forall i \in \mathcal{I}, u \in \mathcal{S}, v \in \mathcal{S}_{u+1},  \label{eq:nobookingBoardingLimitConstraint1}\\
r_{iuv} & =  w_{iuv} - b_{iuv} & \forall i \in \mathcal{I}, u \in \mathcal{S}, v \in \mathcal{S}_{u+1}.\label{eq:nobookingDetainedConstraint1}
\end{align}

Lastly, constraints~\eqref{eq:bookingInvechileConstraint}–\eqref{eq:InvehicleLimitConstraint} ensure that the train capacity limit is respected at all stations. Constraints (\ref{eq:bookingInvechileConstraint}) compute the number of in-vehicle passengers with reservations when train $i$ leaves station $u$,  representing those who have boarded at or before $u$ and have not yet reached their destinations. Constraints (\ref{eq:nobookingInvehicleConstraint}) track the number of in-vehicle passengers without reservations, denoted by $o_{iu}$. For the first station, this value equals the number of boarding non-reserved passengers. For intermediate stations, the value is updated by subtracting the number of alighting passengers ($l_{iu}$, defined in constraints \eqref{eq:nobookingAlightingConstraint}) and adding the number of new boarding passengers at station $u$. At the final station, the number of onboard passengers is set to zero. Constraints (\ref{eq:nobookingAlightingConstraint}) are formulated to track the number of alighting passengers without reservations. Constraints~\eqref{eq:InvehicleLimitConstraint} models the coupling relations between passengers with and without reservations, which require the total number of in-vehicle passengers with and without reservations cannot exceed the maximum train capacity $C^{max}$.
\begin{align}
\hat{o}_{iu} & = \sum_{m \leq u, m \in \mathcal{S}}\sum_{v \in \mathcal{S}_{u+1}}\hat{b}_{imv} &\forall i \in \mathcal{I}, u \in \mathcal{S}, \label{eq:bookingInvechileConstraint} \\
o_{iu} & =
\begin{cases}
\sum\limits_{v \in \mathcal{S}_{u+1}}b_{iuv} & \text{if} \ u=1\\
o_{i(u-1)}-l_{iu}+\sum\limits_{v \in \mathcal{S}_{u+1}}b_{iuv} & \text{if} \ u \in \mathcal{S} \backslash \{1, \left|\mathcal{S}\right|\}\\
0 & \text{if} \ u =\left|\mathcal{S}\right|\\
\end{cases}  & \forall i \in \mathcal{I}, \label{eq:nobookingInvehicleConstraint}\\
l_{iu} & =
\begin{cases}
0 & \text{if} \ u=1\\
\sum\limits_{m =1} ^{u-1} b_{imu} & \text{if} \ u \in \mathcal{S} \backslash \{1\}\\
\end{cases}   & \forall i \in \mathcal{I}, \label{eq:nobookingAlightingConstraint}\\
o_{iu}+\hat{o}_{iu} & \leq C^{max} & \forall i \in \mathcal{I}, u \in \mathcal{S}. \label{eq:InvehicleLimitConstraint}
\end{align}

\textbf{Domains of decision variables.} Constraints (\ref{eq:timetabledominConstraint}) and (\ref{eq:shiftdominConstraint}) define variable domains.
\begin{align}
    \mathbf{x}, \mathbf{z} &
 \in \{0, 1\}^{\left|\mathcal{I}\right| \times \left|\mathcal{T}\right|}, \label{eq:timetabledominConstraint} \\
\mathbf{d}, \mathbf{h}, \boldsymbol{\kappa}, \mathbf{b}, \mathbf{w}, \mathbf{o}, \mathbf{l}, \hat{\mathbf{b}}, \hat{\mathbf{o}} & 
\in \mathbb{Z}_{+}. \label{eq:shiftdominConstraint}
\end{align}
Based on the above discussions, we can now formalize the model as follows
\begin{mini!}|l|[0]<b>
{}
{\eqref{eq:objective}} 
{\label{model:orignal}} 
{}
\addConstraint
{\eqref{eq:objWaiting} - \eqref{eq:shiftdominConstraint}.}
{}
\end{mini!}

\subsection{Model extensions}
\label{sec:ModelExtensions}

In this section, to demonstrate the generalizability of our proposed modeling framework, we extend our formulation for incorporating peak/off-peak pricing (i.e., congestion pricing) and elastic passenger demand from other transportation modes. Furthermore, we scale up the formulation to the network level to highlight the scalability of our modeling approach. These extensions are developed on top of our proposed model \eqref{model:orignal}, rather than constructing a new framework from the ground up. These three extensions illustrate how the proposed modelling framework can accommodate broader operational considerations and address more complex integrated demand-side management and train scheduling challenges.

\textit{\textbf{(a) Extension with the off-peak and peak pricing policy.}} In reality, some rail transit systems, such as the metro in Washington DC and London, employ the time-varying pricing policy, known as \textit{Off-Peak and Peak Pricing}. This approach incentivizes passengers to travel during less congested periods, while charging higher fares during peak hours. Our formulation~\eqref{model:orignal} can be extended to incorporate this pricing policy, as detailed below:
\begin{mini!}|s|[2]<b>
{}
{ F^{t} \label{eq:timeVaryingObjective} } 
{} 
{}
\addConstraint
{\sum\limits_{u\in \mathcal{S}}\sum\limits_{v \in \mathcal{S}_{u+1}}\sum\limits_{t \in \mathcal{T}}\nu_{uvt}(\sum\limits_{t \leq t' \leq \min\{\left|\mathcal{T} \right|, t +\imath\}}}
{\kappa_{uvtt'}+\hat{D}_{uvt}) \geq \nonumber}
{}
\addConstraint
{}
{\qquad \qquad \qquad \qquad \qquad \qquad \qquad \varkappa \sum\limits_{u\in \mathcal{S}}\sum\limits_{v \in \mathcal{S}_{u+1}}\sum\limits_{t \in \mathcal{T}} \varepsilon_{uv} (D_{uvt}+\hat{D}_{uvt}), \label{eq:TimeVaryingRevenue}}
{}
\addConstraint
{}
{\ \eqref{eq:objWaiting}, \eqref{eq:waiting1} - \eqref{eq:waiting4}, \eqref{eq:0-1timetableConstraint}- \eqref{eq:shiftdominConstraint}, \nonumber}
{}
\end{mini!}
where $\nu_{uvt}$ is the ticket fare at timestamp $t$ for OD from stations $u$ to $v$ under the off-peak and peak pricing policy, and $\varkappa \in [0, 100]$ (unit: \%) represents the percentage of operator's revenue under the off-peak and peak pricing policy versus the revenue under static ticket fares. Here, the objective function~\eqref{eq:timeVaryingObjective} aims to minimize the total waiting time of passengers. Constraints~\eqref{eq:TimeVaryingRevenue} are formulated to ensure that the operator's fare revenue is not less than $\varkappa$ times the original revenue with the static ticket price.

\noindent \textit{\textbf{(b) Extension with the elastic passenger demand from other transportation modes.}} Considering the entire urban transportation system, which includes various modes, the off-peak and peak pricing policy in the URT network, which provides lower fares during off-peak periods, has the potential to enhance the attractiveness and cost-effectiveness of the URT. Specifically, this policy could encourage passengers who typically use alternative transportation modes to shift to the URT system during low-peak hours. 

To formulate this extended problem, we introduce a parameter $\epsilon_{t}$ to represent the scaling coefficient at timestamp $t$. A decision variable $\Lambda_{uvt}$ is defined to denote the ticket fare at timestamp $t$ for OD from stations $u$ to $v$ under the off-peak and peak pricing policy. Furthermore, we formulate the number of elastic passengers from other transportation modes who head to station $v \in \mathcal{S}_{u+1}$ and shift to take the URT at station $u \in \mathcal{S}$ and time $t \in \mathcal{T}$ as:

$$\epsilon_{t}\frac{\max\limits_{t'' \in \mathcal{T}} \{\Lambda_{uvt''}\}-\Lambda_{uvt}}{\Lambda_{uvt}}D_{uvt}.$$
When $\Lambda_{uvt} = \max\limits_{t'' \in \mathcal{T}}\{\Lambda_{uvt''}\}$, indicating the ticket price at timestamp $t$ for OD from stations $u$ to $v$ is the highest value during the study time horizon, lacks additional appeal. Conversely, the attractiveness coefficient for passengers of other modes is given by $\epsilon_{t}\frac{\max\limits_{t'' \in \mathcal{T}} \{\Lambda_{uvt''}\}-\Lambda_{uvt}}{\Lambda_{uvt}}$.

Based on the above definitions, we can now formulate this problem as follows:
\begin{mini!}|s|[2]<b>
{}
{F^{t} \nonumber} 
{} 
{}
\addConstraint
{\sum\limits_{\max\{0, t-\imath\} \leq t' \leq t }\kappa_{uvt't}}
{= D_{uvt} + \epsilon_{t}\frac{\max\limits_{t'' \in \mathcal{T}} \{\Lambda_{uvt''}\}-\Lambda_{uvt}}{\Lambda_{uvt}}D_{uvt} \label{eq:nobookingnonpeakingshiftConstraint}}
{\ \forall u \in \mathcal{S}, v \in \mathcal{S}_{u+1}, t \in \hat{\mathcal{T}},}
\addConstraint
{\kappa_{uvt't}}
{=
\begin{cases}
 D_{uvt} + \epsilon_{t}\frac{\max\limits_{t'' \in \mathcal{T}} \{\Lambda_{uvt''}\}-\Lambda_{uvt}}{\Lambda_{uvt}}D_{uvt} & \text{if} \ t'=t\\
0 & \text{otherwise} \\
\end{cases}  \label{eq:nobookingnonpeakingExtraConstraint}}
{\  \forall u \in \mathcal{S}, v \in \mathcal{S}_{u+1}, t \in \mathcal{T} \backslash \hat{\mathcal{T}},}
\addConstraint
{\sum\limits_{i \in \mathcal{I}} b_{iuv}}
{= \sum\limits_{t \in \mathcal{T}} D_{uvt} + \epsilon_{t}\frac{\max\limits_{t'' \in \mathcal{T}} \{\Lambda_{uvt''}\}-\Lambda_{uvt}}{\Lambda_{uvt}}D_{uvt} \label{eq:nobookingFullExtraservedConstraint}}
{\ \forall u \in \mathcal{S}, v \in \mathcal{S}_{u+1},}
\addConstraint
{\sum\limits_{u\in \mathcal{S}}\sum\limits_{v \in \mathcal{S}_{u+1}}\sum\limits_{t \in \mathcal{T}}\Lambda_{uvt}(\sum\limits_{t \leq t' \leq \min\{\left|\mathcal{T} \right|, t +\imath\}}}
{\kappa_{uvtt'}+\hat{D}_{uvt}) \nonumber}
{}
\addConstraint
{ \qquad \qquad \quad \geq \varkappa \sum\limits_{u\in \mathcal{S}}\sum\limits_{v \in \mathcal{S}_{u+1}}\sum\limits_{t \in \mathcal{T}} \varepsilon_{uv} (D_{uvt}+\hat{D}_{uvt}), \label{eq:ETimeVaryingRevenue}}
{}
\addConstraint
{\Lambda_{uvt} }
{\in \mathbb{R}_{+} \label{eq:domin_Price}}
{\ \forall u \in \mathcal{S}, v \in \mathcal{S}_{u+1}, t \in \mathcal{T},}
\addConstraint
{}
{\ \eqref{eq:objWaiting}, \eqref{eq:waiting1} - \eqref{eq:waiting4}, \eqref{eq:0-1timetableConstraint}-\eqref{eq:0-1headwayConstraint}, \eqref{eq:bookingBoardingConstraint}- \eqref{eq:shiftdominConstraint}. \nonumber}
{}
\end{mini!}
The objective function aims to minimizes the total waiting time of passengers. Similar to constraints \eqref{eq:nobookingshiftConstraint} and \eqref{eq:nobookingnonpeakingConstraint}, constraints \eqref{eq:nobookingnonpeakingshiftConstraint} and \eqref{eq:nobookingnonpeakingExtraConstraint} are formulated to calculate the number of arrival passengers at each timestamp and station after passengers shifting their departure times. Constraints \eqref{eq:nobookingFullExtraservedConstraint} ensure that the already loyal to the metro or newly attracted passengers are all served. Constraints~\eqref{eq:ETimeVaryingRevenue} guarantee that the operator's revenue under the off-peak and peak pricing policy is not less than $\varkappa$ times the original revenue with the static passenger demand and the ticket price. Lastly, constraints \eqref{eq:domin_Price} give the domain of the decision variable $\Lambda_{uvt}$.

\noindent \textit{\textbf{(c) Extension of scaling up to the network level.}} 

The key ideas for scaling the proposed model for the integrated optimization of demand management and timetabling to the network level are as follows. We still assume that a limited number of reservation slots will be allocated for each OD pair. Passengers who book a reservation can access the platform upon arrival at the origin station. When they need a transfer, they can proceed via dedicated transfer lanes directly to the platform of the line they are transferring to.

For passengers without reservations, we assume that the metro operators implement trip shifting and passenger flow control strategies for them. They need to wait for permission to enter the platform at the origin station according to the passenger flow control plan, and then pay the full fare or a discounted fare depending on whether they follow the travel time recommended by the system or not. For the passenger flow control, we follow the principle commonly employed in the literature related to optimizing passenger flow control at the network level, see, \cite{Lu2022} and \cite{Yuan2022}. Our proposed model~\eqref{model:orignal} for metro lines can be scaled to the network level by increasing the line dimension and transfer constraints. The detailed formulation is presented in Appendix \ref{sec:network}. 

\section{Solution methodology}
\label{sec:solutionMethods}
In this section, we first reformulate the INLP model~\eqref{model:orignal} into a linear version in Section \ref{sec:modelReformulation}. In Section \ref{sec:BD}, we introduce the Benders decomposition approach. In Section~\ref{sec:cutSeparation}, we detail the Benders cut separation. The accelerating strategies are introduced in Section~\ref{sec:acceleration}. Lastly, the overall framework of our solution method is presented in Section~\ref{sec:framework}.

\subsection{Model reformulation}
\label{sec:modelReformulation}

Our solution approach is based on the Benders decomposition (BD) method, which divides the problem into a \textit{relaxed master problem} (RMP) and a \textit{subproblem} (SP). This approach requires access to the SP's dual information, which in turn necessitates that both the RMP and SP are linear and that the SP contains only continuous variables. RMP is obtained by projecting out the decision variables in the SP and contains \textit{Benders cuts}, which include the \textit{optimality cuts} and \textit{feasibility cuts}. The optimality cuts are generated if the SP is feasible. In cases where the SP is infeasible, a separate \textit{feasibility subproblem} (FS) is solved to generate feasibility cuts. The solution of RMP provides a lower bound of the optimal value. In the computational process, an iterative solution procedure is designed where the RMP is firstly solved, the information from the RMP is passed to the SP, and then dual information from the SP (or the FS) is obtained to generate Benders cuts. In our implementation, we construct a branch-and-cut tree for the RMP and solve the SP (or the FS) at each node, generating Benders cuts that are subsequently added to the RMP.  

Note that our original model~\eqref{model:orignal} contains nonlinear terms in constraints~\eqref{eq:objWaiting}, ~\eqref{eq:waiting2}, ~\eqref{eq:waiting3}, ~\eqref{eq:waiting4} and \eqref{eq:waitingConstraint}. The original nonlinear model is first linearized to obtain its equivalent integer linear programming (ILP) form since the BD approach needs to use the dual information of the model. Then, the ILP formulation is relaxed where the passenger-related variables are relaxed to be continuous ones. This relaxed model is decomposed into an RMP where the timetabling and assignments of passengers with reservations are determined, and an SP to optimize the passenger flow control decisions.  Lastly, the optimal timetable obtained from the relaxed model is fixed and used as input to the ILP model to generate the optimal integer demand-side management solutions.

The linearization process is detailed as follows, which covers each step of converting the nonlinear terms into linear expressions. First, we introduce auxiliary variables $q_{iut}$ for all $i \in \mathcal{I}, u \in \mathcal{S}, t \in \mathcal{T}$ to linearize constraints~\eqref{eq:objWaiting}. Specifically, let $q_{iut} = x_{it} \sum\limits_{v \in \mathcal{S}_{u+1}}r_{iuv}$, we have
\begin{align}
\left\{
\begin{aligned}
& q_{iut} \leq M_{u} x_{it}\\
& q_{iut} \leq \sum\limits_{v \in \mathcal{S}_{u+1}}r_{iuv}\\
& q_{iut} \geq \sum\limits_{v \in \mathcal{S}_{u+1}}r_{iuv} - M_{u} (1-x_{it}) \\
& q_{iut} \in [0, M_{u}]
\end{aligned}
\right.
\quad \forall i \in \mathcal{I}, u \in \mathcal{S}, t \in \mathcal{T}. \label{eq:objWaitingLinearTerms}
\end{align}
Therefore, the nonlinear constraints~\eqref{eq:objWaiting} are reformulated as the following linear version:
\begin{align}
    F^t = \sum_{i\in\mathcal{I}}\sum_{u\in\mathcal{S}}\sum_{t\in\mathcal{T}}(\hat{p}^{wc}_{iut}+p^{wc}_{iut}) + \sum\limits_{i \in \mathcal{I}}\sum_{u\in\mathcal{S}}\sum\limits_{t \in \mathcal{T}}q_{iut}. \label{eq:objWaitingLinear}
\end{align}

Besides, the linear version of constraints~\eqref{eq:waiting1} which contain a nonlinear term of multiplication of a 0-1 variable with an integer variable is formulated as follows: 
\begin{align} \label{eq:waiting1Linear}  
\left\{
\begin{array}{ll}
& p^{w}_{iut}\leq M_{ut} x_{it} \\
& p^{w}_{iut} \leq  \sum\limits_{v\in\mathcal{S}_{u+1}}\sum\limits_{t\leq t^{'}\leq \min\{ \left|\mathcal{T} \right|, t +\imath \}} \kappa_{uvtt^{'}} \\
& p^{w}_{iut} \geq \sum\limits_{v\in\mathcal{S}_{u+1}}\sum\limits_{t\leq t^{'}\leq \min\{ \left|\mathcal{T} \right|, t +\imath \}} \kappa_{uvtt^{'}} - M_{ut} (1-x_{it})  \\
& p^{w}_{iut} \in [0, M_{ut}] \\
\end{array} 
\right.
\quad \forall i \in \mathcal{I}, u \in \mathcal{S}, t \in \mathcal{T}.
\end{align}

To derive equivalent linear forms of constraints~\eqref{eq:waiting3} and \eqref{eq:waiting4}, we first introduce auxiliary variables $\theta_{itt'}, \forall i \in \mathcal{I}, t \in \mathcal{T}, t'\leq t$. Let $\theta_{itt'} = x_{it}x_{it'}$, we have
\begin{align}
\left\{
\begin{aligned}
& \theta_{itt'} \leq x_{it}\\
& \theta_{itt'} \leq x_{it'}\\
& \theta_{itt'} \geq x_{it} + x_{it'} - 1\\
& \theta_{itt'} \in [0, 1],
\end{aligned}
\right.
& \quad \forall i \in \mathcal{I}, t, t' \in \mathcal{T}, t'\leq t.
\end{align}
The linear version of constraints~\eqref{eq:waiting3} can be expressed as
\begin{align}\label{eq:waiting3Linear}
    \hat{p}^{wc}_{iut} = \sum\limits_{t' \in \mathcal{T}, t' \leq t}\sum\limits_{v \in \mathcal{S}_{u+1}}\theta_{itt'}\hat{D}_{uvt'} \quad \forall i \in \mathcal{I}, u \in \mathcal{S}, t \in \mathcal{T}.
\end{align}

Further, by introducing auxiliary variables $\mu_{iutt'}, \forall i \in \mathcal{I}, u \in \mathcal{S}, t \in \mathcal{T}, t'\leq t$, let $\mu_{iutt'} = \theta_{itt'} \sum\limits_{v \in \mathcal{S}_{u+1}}\sum\limits_{t'\leq t''\leq \min\{ \left|\mathcal{T} \right|, t^{'} +\imath \}}\kappa_{uvt't''}$, we have
\begin{align}\label{eq:TwofoldLinear}
\left\{
\begin{aligned}
& \mu_{iutt'} \leq M_{ut'}\theta_{itt'} \\
& \mu_{iutt'} \leq \sum\limits_{v \in \mathcal{S}_{u+1}}\sum\limits_{t'\leq t''\leq \min\{ \left|\mathcal{T} \right|, t^{'} +\imath \}}\kappa_{uvt't''} \\
& \mu_{iutt'} \geq \sum\limits_{v \in \mathcal{S}_{u+1}}\sum\limits_{t'\leq t''\leq \min\{ \left|\mathcal{T} \right|, t^{'} +\imath \}}\kappa_{uvt't''} - M_{ut'}(1-\theta_{itt'}) \\
& \mu_{iutt'} \in [0, M_{ut'}]
\end{aligned}
\right.
\quad \forall i \in \mathcal{I}, u \in \mathcal{S}, t \in \mathcal{T}, t'\leq t.
\end{align}

The linear form of constraints~\eqref{eq:waiting4} can be expressed as
\begin{align}\label{eq:waiting4Linear}
p^{wc}_{iut} = \sum\limits_{v \in \mathcal{S}_{u+1}}\sum\limits_{t' \in \mathcal{T}, t'\leq t}\sum\limits_{t'\leq t''\leq \min\{ \left|\mathcal{T} \right|, t^{'} +\imath \}}\mu_{iuvtt't''} \quad \forall i \in \mathcal{I}, u \in \mathcal{S}, t \in \mathcal{T}.
\end{align}

Lastly, to linearize  constraints~\eqref{eq:waitingConstraint}, we define auxiliary variables $\Gamma_{iuvt'}$, let $w_{iuvt'} = z_{it'}\sum\limits_{t' \leq t \leq \min\{ \left|\mathcal{T} \right|, t' +\imath \} }\kappa_{uvt't}$  for all $i \in \mathcal{I}, u \in \mathcal{S}, v \in \mathcal{S}_{u+1}, t' \in \mathcal{T}$, the linearization results of this term in constraint \eqref{eq:waitingConstraint} are shown below:
\begin{align} \label{eq:waitingConstraintLinearTerms}  
\left\{
\begin{array}{ll}
& \Gamma_{iuvt'}\leq M_{uvt'} z_{it'}\\
& \Gamma_{iuvt'}\leq  \sum\limits_{t' \leq t \leq \min\{ \left|\mathcal{T} \right|, t' +\imath \} }\kappa_{uvt't}\\
& \Gamma_{iuvt'} \geq \sum\limits_{t' \leq t \leq \min\{ \left|\mathcal{T} \right|, t' +\imath \} }\kappa_{uvt't} - M_{uvt'} (1-z_{it'}) \\
& \Gamma_{iuvt'} \in [0, M_{uvt'}] \\
\end{array} 
\right.
\quad \forall i \in \mathcal{I}, u \in \mathcal{S}, v \in \mathcal{S}_{u+1}, t' \in \mathcal{T}.
\end{align}

Thus, we have 
\begin{align}\label{eq:waitingConstraintLinear}
w_{iuv}=
\begin{cases}
\sum\limits_{t' \in \mathcal{T}}\Gamma_{iuvt'} & \text{if} \ i=1\\
\sum\limits_{t' \in \mathcal{T}} \Gamma_{iuvt'} - \sum\limits_{j = 1}^{i-1} b_{juv} & \text{if} \ i \in \mathcal{I} \backslash \{1\}\\
\end{cases}
&\quad \forall u \in \mathcal{S}, v \in \mathcal{S}_{u+1}.
\end{align}

The full formulation of the ILP model is presented as follows: 
\begin{mini}
{}
{\omega_t F^{t}+ \omega_s F^{s}} 
{} 
{\label{model:MILP}}
\addConstraint
{}
{\eqref{eq:objFare} - \eqref{eq:waiting1}, \eqref{eq:0-1timetableConstraint} - \eqref{eq:bookingBoardingConstraint}, \eqref{eq:nobookingBoardingLimitConstraint1} - \eqref{eq:shiftdominConstraint}, \eqref{eq:objWaitingLinearTerms}-\eqref{eq:waitingConstraintLinear}.}
{}
\end{mini}

For the sake of clarity, we now present model~\eqref{model:MILP} as follows:
\begin{align*}
\mathcal{O} = 
    \min\{\omega_t F^{t}+ \omega_s F^{s} | f(\mathbf{z}, \boldsymbol{\kappa}, \mathbf{b}) \geq 0, \mathbf{z} \in \{0, 1\}^{\left|\mathcal{I}\right| \times \left|\mathcal{T}\right|}, \boldsymbol{\kappa}, \mathbf{b} \in \mathbb{Z}_{+}\}.
\end{align*}

\subsection{Benders decomposition}
\label{sec:BD}
 
To enable an exact evaluation of the generated timetable during the solution process, we employ Benders decomposition (BD) rather than a heuristic algorithm. In our implementation, we first relax $\boldsymbol{\kappa}$ and $\mathbf{b}$ as continuous variables to facilitate adding the BD cuts, and then, in the final step, the optimized high-quality timetable is used as input to the integer programming model to produce a demand-side management solution that is directly usable by field staff. This approach ensures both computational efficiency and practical applicability. Specifically, the relaxed model $\tilde{\mathcal{O}}$ can be expressed as follows
\begin{align}\label{model:relaxed}
\tilde{\mathcal{O}} = 
    \min\{\omega_t F^{t}+ \omega_s F^{s} | f(\mathbf{z}, \boldsymbol{\kappa}, \mathbf{b}) \geq 0, \mathbf{z} \in \{0, 1\}^{\left|\mathcal{I}\right| \times \left|\mathcal{T}\right|}, \boldsymbol{\kappa}, \mathbf{b} \geq 0\}.
\end{align}

In the literature, a large body of work employs the aforementioned decomposition method, see \citep{DI2022, YIN2023}. However, by omitting constraints that requires all passengers must be served and the strict limitation on capacity during the timetabling process, the generated timetables may result in infeasibilities. In this case, feasibility cuts are generated and incorporated into the RMP \citep{benders1962}. To reduce the high number of iterations required to generate relatively weak feasibility cuts, we introduce a novel decomposition approach that incorporates full information about passengers with reservations and partial information about those without reservations into the RMP. Specifically, model~\eqref{model:relaxed} is decomposed into an RMP, which addresses train timetabling and trip-shifting plans, and an SP with fixed timetables to determine passenger flow control plans.


Furthermore, we embed the BD approach into the branch-and-cut framework. To facilitate efficient cut separation within modern MIP solvers such as GUROBI, we adopt the Modern Benders Decomposition framework proposed by \cite{Fischetti2016}. In this framework, Benders cuts are generated using reduced cost vectors obtained from a modified version of the subproblem, in which \textit{variable-fixing constraints} are introduced. This approach avoids the need for explicit Lagrangian dual derivation, and integrates naturally with the commercial solver's callback mechanism. Specifically, at the beginning, we relax the timetabling-related decision variables (i.e., $z_{it}$) to continuous ones and solve the RMP to the optimum, where no Benders cuts are included. If all the variables $z_{it}$ are integer, we solve the SP and add optimality or feasibility cuts to the RMP. Otherwise, we select a fractional variable $z_{it}$ to run the branch-and-cut process. At each node, we solve the SP and incorporate information from the SP into the RMP. The RMP incorporates only a subset of the Benders cuts, which are added following the procedure proposed in Section~\ref{sec:cutSeparation}. The RMP and the SP can now be formulated as 
\begin{mini!}|s|[2]<b>
{}
{\omega_t \sigma\sum_{i\in\mathcal{I}}\sum_{u\in\mathcal{S}}\sum_{t\in\mathcal{T}}(\hat{p}^{wc}_{iut}+p^{wc}_{iut})+\omega_s F^{s}+ \theta} 
{\label{model:RMP}} 
{}
\addConstraint
{\theta}
{\geq \Omega(\mathbf{z}^{*}_l, \boldsymbol{\kappa}^{*}_l)+\boldsymbol{\xi}^{T}_{l}(\mathbf{z}-\mathbf{z}^{*}_l)+\boldsymbol{\chi}^{T}_{l}(\boldsymbol{\kappa}-\boldsymbol{\kappa}^{*}_l)}
{\ \forall l \in \{1, 2, 3, ..., c_1\},}
\addConstraint
{0 }
{\geq \Psi(\mathbf{z}^{*}_l, \boldsymbol{\kappa}^{*}_l)+\boldsymbol{\xi}^{T}_{l}(\mathbf{z}-\mathbf{z}^{*}_l)+\boldsymbol{\chi}^{T}_{l}(\boldsymbol{\kappa}-\boldsymbol{\kappa}^{*}_l)}
{\ \forall l \in \{1, 2, 3, ..., c_2\},}
\addConstraint
{\hat{o}_{iu} }
{\leq C^{max} \label{eq:RelaxedInvehicleLimitConstraint}}
{\ \forall i \in \mathcal{I}, u \in \mathcal{S},}
\addConstraint
{}
{ \eqref{eq:objFare}, \eqref{eq:waiting1}, \eqref{eq:0-1timetableConstraint} - \eqref{eq:nobookingnonpeakingConstraint}, \eqref{eq:waiting3Linear} - \eqref{eq:waiting4Linear},}
{}
\addConstraint
{\mathbf{z}}
{\in [0, 1]^{\left|\mathcal{I}\right| \times \left|\mathcal{T}\right|},}
{}
\addConstraint
{\boldsymbol{\kappa}, \theta}
{\geq 0,}
{}
\end{mini!}
where $c_1$ and $c_2$ indicates the number of added optimality and feasibility cuts, respectively. $\Omega(\mathbf{z}^*, \boldsymbol{\kappa}^*)$ indicates the objective function of SP under the feasible solution of the RMP, i.e., ($\mathbf{z}^*, \boldsymbol{\kappa}^{*}$). The auxiliary decision variable $\theta$ approximates the objective function of SP. $\Psi(\mathbf{z}, \boldsymbol{\kappa})$ denotes the objective function of the FS. For SP, the feasible solution of the RMP ($\mathbf{z}^*, \boldsymbol{\kappa}^{*}$) is fixed. That is,
\begin{mini!}|s|[2]<b>
{}
{\Big\{ \omega_t\sigma\sum\limits_{i \in \mathcal{I}}\sum_{u\in\mathcal{S}}\sum\limits_{t \in \mathcal{T}}q_{iut} \Big\} }
{\label{model:SP}} 
{ \quad \Omega(\mathbf{z}^*, \boldsymbol{\kappa}^{*}) =}
\addConstraint
{\mathbf{z}}
{= \mathbf{z}^{*}, \label{eq:BDz}}
{}
\addConstraint
{\boldsymbol{\kappa}}
{= \boldsymbol{\kappa}^{*}, \label{eq:BDkappa}}
{}
\addConstraint
{\mathbf{b}}
{\geq 0, \label{eq:BDb}}
{}
\addConstraint
{}
{\eqref{eq:0-1headwayConstraint} - \eqref{eq:nobookingFullservedConstraint}, \eqref{eq:nobookingBoardingLimitConstraint1} - \eqref{eq:nobookingDetainedConstraint1}, \eqref{eq:nobookingInvehicleConstraint} - \eqref{eq:InvehicleLimitConstraint}, \eqref{eq:objWaitingLinearTerms},  \eqref{eq:waitingConstraintLinearTerms} - \eqref{eq:waitingConstraintLinear}, \label{eq:BDsp}}
{}
\end{mini!}
where constraints~\eqref{eq:BDz} and \eqref{eq:BDkappa} are the variable-fixing constraints. 

The FS that are used to generate feasibility cuts  can be expressed as follows:
\begin{mini!}|s|[2]<b>
{}
{\Big\{\sum\limits_{i \in \mathcal{I}}\sum_{u\in\mathcal{S}} \tau_{iu} \Big\} }
{\label{model:FS}} 
{ \quad \Psi(\mathbf{z}^*, \boldsymbol{\kappa}^{*}) =}
\addConstraint
{\mathbf{z}}
{= \mathbf{z}^{*}, }
{}
\addConstraint
{\boldsymbol{\kappa}}
{= \boldsymbol{\kappa}^{*}, }
{}
\addConstraint
{\mathbf{b}}
{\geq 0, }
{}
\addConstraint
{o_{iu} + \hat{o}_{iu}}
{\leq C^{max} + \tau_{iu} }
{\forall i \in \mathcal{I}, u \in \mathcal{S},}
\addConstraint
{}
{\eqref{eq:0-1headwayConstraint} - \eqref{eq:nobookingFullservedConstraint}, \eqref{eq:nobookingBoardingLimitConstraint1} - \eqref{eq:nobookingDetainedConstraint1}, \eqref{eq:nobookingInvehicleConstraint} - \eqref{eq:nobookingAlightingConstraint},}
{}
\addConstraint
{}
{\eqref{eq:objWaitingLinearTerms},  \eqref{eq:waitingConstraintLinearTerms} - \eqref{eq:waitingConstraintLinear}.}
\end{mini!}

\subsection{Benders cut separation}
\label{sec:cutSeparation}

As introduced in Section~\ref{sec:modelReformulation}, we solve the model~\eqref{model:relaxed} using a branch-and-cut method and thus dynamically include the Benders cuts while exploring the branch-and-bound tree. At each node in the branch-and-bound tree, optimality or feasibility cuts for the decision variable $\mathbf{z}$ and $\boldsymbol{\kappa}$ is generated using the dual information of the SP and added into the RMP. In the following, we provide a detailed description of each type of cuts.

First, function $\Omega(\mathbf{z}, \boldsymbol{\kappa})$ can be underestimated by a supporting hyperplane on $(\mathbf{z}^{*}, \boldsymbol{\kappa}^{*})$ because of convexity. If SP~\eqref{model:SP} given solution $(\mathbf{z}^{*}, \boldsymbol{\kappa}^{*})$ is feasible, the \textit{optimality cuts} are generated as follows:
\begin{align}\label{eq:optimalityCuts}
&\theta \geq \Omega(\mathbf{z}^{*}, \boldsymbol{\kappa}^{*})+\boldsymbol{\xi}^{T}(\mathbf{z}-\mathbf{z}^{*})+\boldsymbol{\chi}^{T}(\boldsymbol{\kappa}-\boldsymbol{\kappa}^{*}) &
 \end{align}
where $\Omega(\mathbf{z}^{*}, \boldsymbol{\kappa}^{*})$ is the objective value of the SP, $\boldsymbol{\xi}$ and $\boldsymbol{\chi}$ represent the dual variables of constraints~\eqref{eq:BDz} and \eqref{eq:BDkappa}, respectively.

If SP~\eqref{model:SP} is infeasible, we generate the following \textit{feasibility cuts} by solving FS~\eqref{model:FS} and add them to the RMP:
\begin{align}\label{eq:feasibilityCuts}
0 \geq \Psi(\mathbf{z}^{*}, \boldsymbol{\kappa}^{*})+\boldsymbol{\xi}^{T}(\mathbf{z}-\mathbf{z}^{*})+\boldsymbol{\chi}^{T}(\boldsymbol{\kappa}-\boldsymbol{\kappa}^{*}). 
\end{align}

To strengthen the optimality cuts~\eqref{eq:optimalityCuts}, \cite{Rahmaniani2020} introduces Benders dual decomposition (BDD) approach which generate the \textit{strengthened optimality cuts}. In the BDD method, the local copies of the variables in the RMP are introduced in the SP and then priced out into the objective function. It was proven that strengthened optimality cuts are tighter than the traditional optimality cuts~\eqref{eq:optimalityCuts}. Specifically, when solving our problem with the BDD approach, the variable-fixing constraints~\eqref{eq:BDz} and \eqref{eq:BDkappa} are priced out into the objective function using dual multipliers $\boldsymbol{\xi}$ and $\boldsymbol{\chi}$. By doing so, we obtain the following Lagrangian dual problem
\begin{align}
\max\limits_{\boldsymbol{\xi}, \boldsymbol{\chi}} \min\limits_{\mathbf{z}, \boldsymbol{\kappa}, \mathbf{b}} \Big\{ \omega_t\sigma\sum\limits_{i \in \mathcal{I}}\sum_{u\in\mathcal{S}}\sum\limits_{t \in \mathcal{T}}q_{iut} -\boldsymbol{\xi}^{T}(\mathbf{z}-\mathbf{z}^{*})-\boldsymbol{\chi}^{T}(\boldsymbol{\kappa}-\boldsymbol{\kappa}^{*}): \eqref{eq:BDb}, \eqref{eq:BDsp} \Big\}.
\end{align}

Then, given $\mathbf{z}^{*} \in [0, 1]^{ \left| \mathcal{I} \right| \times \left| \mathcal{T} \right|}$, $\boldsymbol{\kappa}^{*} \in \mathbb{Z}_{+}$, $\boldsymbol{\xi} \in \mathbb{R}$, and $\boldsymbol{\chi} \in \mathbb{R}$, let $(\hat{\mathbf{z}}^{*}, \hat{\boldsymbol{\kappa}}^{*}, \hat{\mathbf{b}}^{*})$ be an optimal solution obtained by solving the following problem
\begin{align}\label{eq:strengthenedoptimality}
\min\limits_{\mathbf{z}, \boldsymbol{\kappa}, \mathbf{b}} \Big\{ \omega_t\sigma\sum\limits_{i \in \mathcal{I}}\sum_{u\in\mathcal{S}}\sum\limits_{t \in \mathcal{T}}q_{iut} -\boldsymbol{\xi^{*}}^{T}(\mathbf{z}-\mathbf{z}^{*})-\boldsymbol{\chi^{*}}^{T}(\boldsymbol{\kappa}-\boldsymbol{\kappa}^{*}): \eqref{eq:BDb}, \eqref{eq:BDsp}, \mathbf{z} \in [0, 1]^{\left|\mathcal{I}\right| \times \left|\mathcal{T}\right|}, \boldsymbol{\kappa} \in \mathbb{Z}_{+},  \mathbf{b} \in \mathbb{R}_{+} \Big\}.
\end{align}

The \textit{strengthened optimality cuts} can be expressed as
\begin{align}\label{eq:strengthenedoptimalityCuts}
&\theta \geq \Omega(\hat{\mathbf{z}}^*, \hat{\boldsymbol{\kappa}}^*)+\boldsymbol{\xi}^{T}(\mathbf{z}-\hat{\mathbf{z}}^*)+\boldsymbol{\chi}^{T}(\boldsymbol{\kappa}-\hat{\boldsymbol{\kappa}}^*). &
 \end{align}
 
A comparison of the performance of six variants with respect the solution algorithm is presented in the following numerical experiments. The variant based on the BD approach solves the SP at each node and then adds an optimality or feasibility cut to the RMP. On the other hand, in the BDD-based solution method, a strengthened optimality or feasibility cut is generated following the solution of the SP at each node and then integrated into the RMP.

\subsection{Strategies for acceleration}
\label{sec:acceleration}

In this section, we first introduce two methodologically accelerating methods to strengthen the linear relaxation bounds at each node in Section~\ref{sec:methodologicallyAccelerated}, subsequently delving into key implementation details which contribute to accelerate computations in Section~\ref{sec:heuristicAcceleration}.

\subsubsection{Accelerating methods}
\label{sec:methodologicallyAccelerated}

Firstly, to further guide the timetabling optimization process in the RMP, we propose the following valid equalities and inequalities that incorporate information from non-reserved passengers, which are added into the RMP~~\eqref{model:RMP} as constraints.

\begin{proposition} \label{proposition}

We define the number of passengers without reservations who head to station $v$ and are ensured to board station $u$ as $\tilde{b}_{iuv}$. Besides, we introduce the number of passengers without reservations who are on board at train $i$, who depart from train $i$ at station $u$, and the total number of boarding passengers without reservations at station $u$ as $\tilde{o}_{iu}$, $\tilde{l}_{iu}$, and $\tilde{b}_{iu}$, respectively. For the RMP~\eqref{model:RMP}, the following equalities and inequalities are valid:
\begin{align}
& \tilde{b}_{iuv} = \varrho_{iuv} \sum\limits_{t' \in \mathcal{T}}x_{it'}\sum\limits_{t' \leq t \leq \min\{ \left|\mathcal{T} \right|, t' +\imath \} }\kappa_{uvt't} & \forall i \in \mathcal{I}, u \in \mathcal{S}, \label{eq:validBoard}\\
&\tilde{o}_{iu}=
   \begin{cases}
    \sum\limits_{v \in \mathcal{S}_{u+1}}\tilde{b}_{iuv} & \text{if} \ u=1\\
    \tilde{o}_{i(u-1)}-\tilde{l}_{iu}+\sum\limits_{v \in \mathcal{S}_{u+1}}\tilde{b}_{iuv} & \text{if} \ u \in \mathcal{S} \backslash \{1, \left|\mathcal{S}\right|\}\\
0 & \text{if} \ u =\left|\mathcal{U}\right|\\
\end{cases}
& \forall i \in \mathcal{I}, \label{eq:validInvehicle}\\
& \hat{o}_{iu} + \tilde{o}_{iu} \leq C^{max},  & \forall i \in \mathcal{I}, u \in \mathcal{S}, \label{eq:validInvehicleTotal}\\
& \tilde{l}_{iu}=
\begin{cases}
0 & \text{if} \ u=1\\
\sum\limits_{m =1} ^{u-1} \tilde{b}_{imu} & \text{if} \ u \in \mathcal{S} \backslash \{1\}\\
\end{cases}
& \forall i \in \mathcal{I},\\
& \tilde{b}_{iu} = \sum\limits_{v \in \mathcal{S}, v \textgreater u}\tilde{b}_{iuv}  & \forall i \in \mathcal{I}, u \in \mathcal{S}. \label{eq:wildeb}
\end{align}
\proof{Proof.}
Recall that we require at least $\varrho$ percent of passengers without reservations to be served by each train to ensure fairness in the SP~\eqref{model:SP}, that is, constraints~\eqref{eq:nobookingBoardingLimitConstraint1}. Therefore, for any $\mathbf{z} \in [0, 1]^{ \left| \mathcal{I} \right| \times \left| \mathcal{T} \right|}$ generated in the RMP, if this timetable is feasible, then it must serve $\varrho$ (unit: percent) of passengers without reservations while satisfying capacity limitation~\eqref{eq:InvehicleLimitConstraint}. Being inspired by this property, we integrate the dynamics of $\varrho$ (unit: percent) of passengers without reservations into the RMP. To be specific, the number of passengers without reservations who head to station $v$ and are guaranteed to board train $i$ at station $u$ is formulated as constraints~\eqref{eq:validBoard}. Thereafter, constraints~\eqref{eq:validInvehicle} model the number of in-vehicle passengers without reservations. Constraints~\eqref{eq:validInvehicleTotal} are formulated to ensure the capacity limitation, which is the key inequalities that enhance the lower bound of the RMP. Hence, Proposition~\ref{proposition} holds.
\Halmos
\endproof

\end{proposition}

Secondly, it is widely recognized that the effectiveness of the branch-and-bound method is significantly influenced by the \textit{Big-$M$} values in formulations, which might return poor bounds and causes large branching
trees. To enhance solution efficiency, we tighten the upper bounds by redefining the \textit{Big-$M$} values in constraints~\eqref{eq:objWaitingLinearTerms}, \eqref{eq:waiting1Linear}, \eqref{eq:TwofoldLinear}, and \eqref{eq:waitingConstraintLinearTerms}. Specifically, considering that the parameter $M_u$ in constraints~\eqref{eq:objWaitingLinearTerms} represents the upper limitation of the detained passengers, its minimum value without cutting out the optimal solution is the one where all the non-reserved passengers are not served during the study time horizon. Thus, we redefine the $M_{u}$ as follows.
\begin{align}
M_{u}    & = \sum\limits_{t \in \mathcal{T}}\sum\limits_{v \in \mathcal{S}_{u+1}}(1-\varrho_{uv})D_{uvt} & u \in \mathcal{S}. 
\end{align}

Similarly, \(M_{ut}\) in constraints~\eqref{eq:waiting1Linear}, \(M_{ut'}\) in constraints~\eqref{eq:TwofoldLinear}, and \(M_{uvt'}\) in constraints~\eqref{eq:waitingConstraintLinearTerms} can be defined as:
\begin{align}
M_{ut}   & = \sum\limits_{v \in \mathcal{S}_{u+1}}\sum\limits_{t\leq t^{'}\leq \min\{ \left|\mathcal{T} \right|, t +\imath \}}D_{uvt'} & \forall u \in \mathcal{S}, t \in \mathcal{T}. \\
M_{ut'}  & = \sum\limits_{v \in \mathcal{S}_{u+1}}\sum\limits_{t'\leq t''\leq \min\{ \left|\mathcal{T} \right|, t^{'} +\imath \}}D_{uvt''}& u \in \mathcal{S}, t' \in \mathcal{T}. \\
M_{uvt'} & = \sum\limits_{t' \leq t \leq \min\{ \left|\mathcal{T} \right|, t' +\imath \} } D_{uvt} & \forall u \in \mathcal{S}, v \in \mathcal{S}_{u+1}, t' \in \mathcal{T}.
\end{align}

\subsubsection{Heuristic accelerating strategy} \label{sec:heuristicAcceleration} Moreover, we accelerate the solution by incorporating the following ideas from \cite{Fischetti2016} and \cite{FischettiMS}.

\textbf{Cut loop stabilization at the root node.} We implement the accelerating method at each cut loop iteration. Specifically, we have two points at the decision variables at each cut loop iteration: the optimal solution $(\mathbf{z}^*,\boldsymbol{\kappa}^*)$ of the current RMP and a stabilizing point $(\mathbf{\widetilde{z}},\boldsymbol{\widetilde{\kappa}})$ that is initialized by solving the following problem:
\begin{equation*}
\max\Big\{\sum\limits_{i\in\mathcal{I}}\sum\limits_{t\in\mathcal{T}}z_{it}+\sum\limits_{u\in\mathcal{S}}\sum\limits_{v\in\mathcal{S}_{u+1}}\sum\limits_{t\in\mathcal{T}}\sum\limits_{t'\leq t\leq \min{\left|\mathcal{T}\right|, t'+ \iota}}\kappa_{uvt't}| (\mathbf{z},\boldsymbol{\kappa}) \in\mathcal{X}\Big\},
\end{equation*}
where $\mathcal{X}$ is the domain of decision variables $(\mathbf{z},\boldsymbol{\kappa})$. At each step, we move $(\mathbf{\widetilde{z}},\boldsymbol{\widetilde{\kappa}})$ towards $(\mathbf{z}^*,\boldsymbol{\kappa}^*)$ by setting $(\mathbf{\widetilde{z}},\boldsymbol{\widetilde{\kappa}}) = (\alpha \mathbf{\widetilde{z}} + (1-\alpha)\mathbf{z}^*, \alpha\boldsymbol{\widetilde{\kappa}} + (1-\alpha)\boldsymbol{\kappa}^* )$ and then apply our optimality cuts to the intermediate point $(\lambda\mathbf{z}^* + (1-\lambda)\mathbf{\widetilde{z}},  \lambda\boldsymbol{\kappa}^* + (1-\lambda) \boldsymbol{\widetilde{\kappa}})$, where parameters $\lambda \in (0, 1]$  and $\alpha\in (0, 1]$. Then, the intermediate point is fed to the SP~\eqref{model:SP}. The optimality cuts~\eqref{eq:optimalityCuts} are updated and statically added to the RMP. After five consecutive iterations which the linear programming (LP) bound does not improve, parameter $\lambda$ is reset to one and the cut loop continues.

\textbf{Tailing off.} We prevent more than $\upsilon$ successive calls of the Benders cut separation functions~\eqref{eq:optimalityCuts} and \eqref{eq:strengthenedoptimalityCuts} ($\hat{\upsilon}$ for the root node) at each fractional node of the branch-and-cut tree. In order to ensure correctness, integer solutions capable of updating the incumbent are always separated.

\textbf{Restart.} Before entering the final branch-and-cut run, we stop the execution right after the root node, add the generated Benders cuts as static cuts to the RMP, update the incumbent, and repeat. This restart mechanism is applied twice before entering the final run in our implementation.

\textbf{Tree search.} We aggressively apply the level of all GUROBI's internal cuts, select the full-strong branching, and use the strongest lower-bound searching strategies. As to heuristics, we apply the relaxation induced neighborhood search heuristic at every node to feed our Benders-cut separator with low-cost integer solutions.

\subsection{Overall framework}
\label{sec:framework}

With these definitions in place, the decomposition algorithm within branch-and-cut framework can be outlined in Algorithm~\ref{alg:algorithm1} in Appendix \ref{sec:overallFramework}, which can be summarized as follows:  (1) At the root node, relax the variable $\mathbf{z}$ and solve the RMP~\eqref{model:RMP} to the optimum. It is worth noting that no optimality cut or feasibility cut has been included. (2) Identify a fractional \( z_{i,t} \notin \{0,1\} \) from the solution at the root node and perform branching on this variable to create two new nodes. (3) At each subsequent node within the branch and bound tree, solve the RMP, and feed its solution into the SP~\eqref{model:SP}. Subsequently, an optimality, strengthened optimality, or feasibility cut is added to the RMP, followed by the evaluation to check if the relaxed Gap, defined as $(UB-LB)/LB \times 100 \%$, falls below the pre-given tolerance $\varepsilon_1 $. Once this criterion is met, the timetabling solution $\mathbf{z}^*_0$ is input into the ILP model (39), which includes the integer decision variables $\boldsymbol{\kappa},\mathbf{b}$. GUROBI is employed to solve this model to optimality. This final integer solution allows on-site staff to implement the optimized passenger flow control plan (i.e., $\mathbf{b}$) directly and without ambiguity, thus aligning perfectly with real-life operational scenarios where passengers are inherently counted as integers. Moreover, in our implementation, Benders cuts are added using GUROBI’s callback framework at RMP solutions. Specifically, feasibility and optimality cuts are generated via the lazy constraint callback at integer solutions, and via the user cut callback at fractional solutions. For each RMP solution:
(1) If the SP is infeasible, a feasibility cut is added;
(2) Otherwise, an optimality cut is generated.
Cuts at integer solutions are added through the lazy constraint callback, while cuts at fractional solutions are added via the user cut callback.

\section{Numerical experiments}
\label{sec:caseStudy}
In this section, two series of numerical experiments are constructed to demonstrate the effectiveness of our proposed approaches. First, in Section \ref{sec:smallCase}, we evaluate the advantages of the proposed models and the performance at the root node of the solution method on a proof-of-concept case study. Thereafter, to gain more insights on the integrated optimization and the full performance of the algorithm, we conduct real-world instances based on the data of Beijing metro in Section~\ref{sec:BeijingMetro}. The model and the solution method are coded in Java in combination with GUROBI 9.5.1. The experiments are run on a personal computer equipped with an Intel i9-14900HX CPU at 2.20 GHz and 64 GB of RAM.

\subsection{Proof-of-concept case study}
\label{sec:smallCase}

The proof-of-concept case study involves a metro line with six stations, as shown in Figure~\ref{fig:smallLine}. It is assumed that the section running times and dwell times at each station are one minute. Furthermore, the minimum and maximum headways are set as two and six minutes, respectively, with each train having a capacity of 600 passengers. For the time horizon, a period of 60 minutes is considered, which is discretized into one-minute timestamps. Time-varying OD passenger demand is randomly generated to simulate realistic situations from peak-hour to off-peak-hour periods. The total time-varying OD passenger demand (i.e., $\hat{D}_{uvt} + D_{uvt}$) is generated for each OD pair $(u,v)$ and time interval $t$, where $\hat{D}_{uvt}$ and $D_{uvt}$ represent the reserved and non-reserved passenger demand, respectively. A reservation ratio is defined to determine $\hat{D}_{uvt}$, and the remaining demand is treated as non-reserved. Sensitivity analyses on the reservation ratio are presented in Section~\ref{sec:insights} to derive managerial insights. The fare on each OD pair is set as 3 RMB. 

To derive insights into the operational strategies integrating booking, directing and timetabling, sensitivity analyses related to various parameter settings are presented in Section \ref{sec:insights}. In Section \ref{sec:rootNode}, a comparison of the lower and upper bounds at the root node among the six variants of the proposed algorithm is conducted to evaluate their effectiveness. 

\begin{figure}[H]
\centering

\tikzset{every picture/.style={line width=0.75pt}} 

\tikzset{every picture/.style={line width=0.75pt}} 

\tikzset{every picture/.style={line width=0.75pt}} 

\begin{tikzpicture}[x=0.75pt,y=0.75pt,yscale=-1,xscale=1]

\draw  [fill={rgb, 255:red, 155; green, 155; blue, 155 }  ,fill opacity=0.28 ] (133,145) -- (169.5,145) -- (169.5,167.75) -- (133,167.75) -- cycle ;
\draw  [fill={rgb, 255:red, 155; green, 155; blue, 155 }  ,fill opacity=0.28 ] (192.5,145) -- (229,145) -- (229,167.75) -- (192.5,167.75) -- cycle ;
\draw  [fill={rgb, 255:red, 155; green, 155; blue, 155 }  ,fill opacity=0.28 ] (252.5,143.5) -- (289,143.5) -- (289,166.25) -- (252.5,166.25) -- cycle ;
\draw  [fill={rgb, 255:red, 155; green, 155; blue, 155 }  ,fill opacity=0.28 ] (312.5,144) -- (349,144) -- (349,166.75) -- (312.5,166.75) -- cycle ;
\draw  [fill={rgb, 255:red, 155; green, 155; blue, 155 }  ,fill opacity=0.28 ] (372,144.5) -- (408.5,144.5) -- (408.5,167.25) -- (372,167.25) -- cycle ;
\draw  [fill={rgb, 255:red, 155; green, 155; blue, 155 }  ,fill opacity=0.28 ] (432.5,144) -- (469,144) -- (469,166.75) -- (432.5,166.75) -- cycle ;
\draw  [color={rgb, 255:red, 155; green, 155; blue, 155 }  ,draw opacity=1 ][fill={rgb, 255:red, 128; green, 128; blue, 128 }  ,fill opacity=1 ] (111.5,129) -- (486.5,129) -- (486.5,133.75) -- (111.5,133.75) -- cycle ;
\draw  [color={rgb, 255:red, 155; green, 155; blue, 155 }  ,draw opacity=1 ][fill={rgb, 255:red, 128; green, 128; blue, 128 }  ,fill opacity=1 ] (111.5,178.5) -- (486.5,178.5) -- (486.5,183.25) -- (111.5,183.25) -- cycle ;

\draw (145,148) node [anchor=north west][inner sep=0.75pt]    {1};
\draw (205,148) node [anchor=north west][inner sep=0.75pt]    {2};
\draw (265,148) node [anchor=north west][inner sep=0.75pt]    {3};
\draw (325,148) node [anchor=north west][inner sep=0.75pt]    {4};
\draw (385,148) node [anchor=north west][inner sep=0.75pt]    {5};
\draw (445,148) node [anchor=north west][inner sep=0.75pt]    {6};

\end{tikzpicture}
\caption{An illustration of the proof-of-concept line.} \label{fig:smallLine}
\end{figure}

\subsubsection{Managerial insights}
\label{sec:insights}

This section assesses the value of integrating directing and timetabling decisions, and provides insights into the trade-off between operational effectiveness and service fairness as well as the interests of passengers and the operator. To do so, we first fix the service fairness parameters and weight coefficients in the objective function, and analyze the resulting service efficiency and the additional government subsidies among various shifting limitations. Further, we construct a set of experiments by varying the booking ratio. Lastly, we use the concept of $\varepsilon$-Constraints to find Pareto solutions, where it is impossible to reduce the additional government subsidies without increasing the waiting time, and vice versa. 

Our managerial insights quantify the benefits of directing and timetabling simultaneously (Insight \ref{insight1}), ensuring service fairness in optimizing passengers' and resources' assignments (Insight \ref{insight2}) and improving operational efficiency through the pricing policy (Insight \ref{insight3}).

\begin{insight}\label{insight1}
\textit{Encouraging 21.14\% of passengers to shift their arrival times by up to 10 minutes through providing a discount on tickets would result in savings of at least 8.33\% in the fleet size and at the expense of a 2.03\% reduction in the additional subsidy.}
\end{insight}

\begin{table}[tbp]
\centering
\caption{Results among various values of maximum allowable shifting time. Abbreviations: STP = Percentage of passengers who shift their trips; AWT-WR = Average waiting time of passengers without reservations; DP = Detained passengers; MW = Maximum number of passengers waiting at stations.}
\label{tab:sensitivityShifting}
\resizebox{\textwidth}{!}{%
\begin{tabular}{rrrrrrr}
\Xhline{1pt}
\# of trains & \begin{tabular}[c]{@{}r@{}} Maximum shifting \\ time (min) \end{tabular} & \begin{tabular}[c]{@{}r@{}}  STP (\%)\end{tabular}&  AWT-WR (min) & \# of DP & \begin{tabular}[c]{@{}r@{}} Additional \\ subsidy (\%) \end{tabular} & MW \\
\Xhline{0.6pt}
\multirow{7}{*}{11} & 0  &Infeasible &$-$  &$-$&$-$  & $-$    \\
                    & 5  &Infeasible &$-$  &$-$ &$-$  & $-$    \\
                    & 10 & 21.14  & 2.48 & 211.00 & 2.03  & 312.00    \\
                    & 15 & 26.25  & 1.91 & 32.00 & 2.52 & 236.00   \\
                    & 20 & 22.48   & 1.93 & 32.00 & 2.16 & 225.00    \\
                    & 25 &22.75 & 1.92 & 32.00 & 2.19 & 236.00 \\
\Xhline{0.6pt}
\multirow{7}{*}{12} & 0  & 0  & 2.53 & 202.00 &0 & 245.00   \\
                    & 5  & 14.22  & 1.93 & 0 & 1.37 & 231.00   \\
                    & 10 & 15.64  & 1.87 & 0 & 1.50 & 227.00   \\
                    & 15 & 16.28  & 1.84 & 0 & 1.56 & 212.00   \\
                    & 20 & 12.76  & 1.93 & 0 & 1.23 & 231.00   \\
                    & 25 & 14.00  & 1.90 & 0 & 1.35 & 214.00   \\
\Xhline{1pt}
\end{tabular}%
}
\end{table}

In Table~\ref{tab:sensitivityShifting}, the results among various maximum values of shifting time are presented, where both the booking ratio (i.e., $\hat{D}_{uvt}/ ( \hat{D}_{uvt}+ D_{uvt} ) \ \forall u, v \in \mathcal{S}, v \geq u, t \in \mathcal{T}$) and the service fairness factor $\varrho_{iu}$ for all $i \in \mathcal{I}, u \in \mathcal{S}$ are set to 50\%, the discount on the ticket price is 20\%, and the weight coefficients are designated as 1 and 5. These results include the percentage of passengers shifting trips (STP), the average waiting time of passengers without reservations (AWT-WR), the number of detained passengers without reservations (DP), the percentage of additional government subsidies, and the maximum congestion experienced at stations during the operation of all trains. SP and AWT-WR are computed as 
\begin{equation*}
\big(\text{Results} / \# \text{ of passenger without reservations} \big) \times 100 ~ (\%),
\end{equation*} 
\begin{equation*}
\big(\text{Additional subsidies} / \text{Revenue when shifting is not allowed} \big) \times 100 ~ (\%).
\end{equation*} 

It can be observed that when the demand management strategy encouraging passengers to shift trips through incentives is not implemented, at least 12 trains are required to serve all demand. However, if passengers are encouraged to shift their travel by up to 10 minutes, it becomes possible to meet all demand with just 11 trains. This finding leads us to conclude that a demand management strategy that promotes passenger shifts through incentives can effectively reduce the needed fleet size.

A second observation is that there is an inverse relationship between the elasticity in arrival time adjustments and the reduction of the number of detained passengers. Specifically, as passengers' shifting behavior is encouraged and the duration of maximum shifting time extends, the number of detained passengers is remarkably reduced. For example, when 11 trains are operated, the number of detained passengers decreases from 211 to 32 as the maximum shifting time changes from 10 to 25 minutes. However, there is no linear connection between the maximum allowable shifting time and the average waiting time. This is because the optimal solution is a trade-off between the waiting time and additional government subsidies. These results indicate that the integration of directing and timetabling policies could increase the effectiveness of the rail transit system. This improvement is demonstrated by a considerable reduction in the waiting time of passengers without doing any serious harm to the operator’s perspective.

\begin{insight}\label{insight2}
\textit{Neglecting the principle of service fairness in the allocation of transit resources increases up to 7.51\% in the percentage of passengers necessitating a shift in their departure times. Besides, it leads to a rise in the total waiting time, observed as 19.74\% while a 1.93\% reduction in the additional government subsidies due to the number of passengers with reservations increases, who do not have the flexibility to shift travel times.}
\end{insight}

\begin{figure}[h]
     \centering
     \begin{subfigure}{0.49\textwidth}
         \centering
 \includegraphics[width=\textwidth]{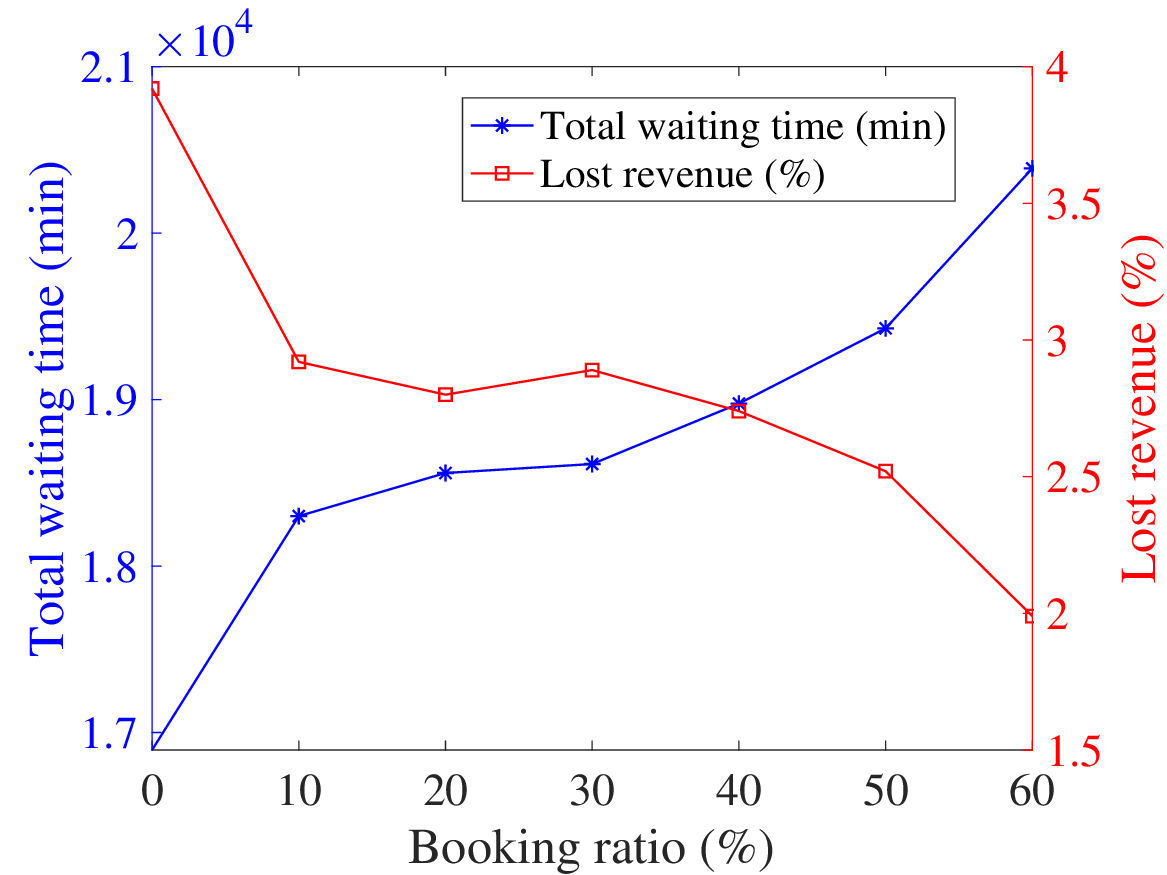}
         \caption{Objectives}
     \end{subfigure}
     \begin{subfigure}{0.49\textwidth}
         \centering
\includegraphics[width=\textwidth]{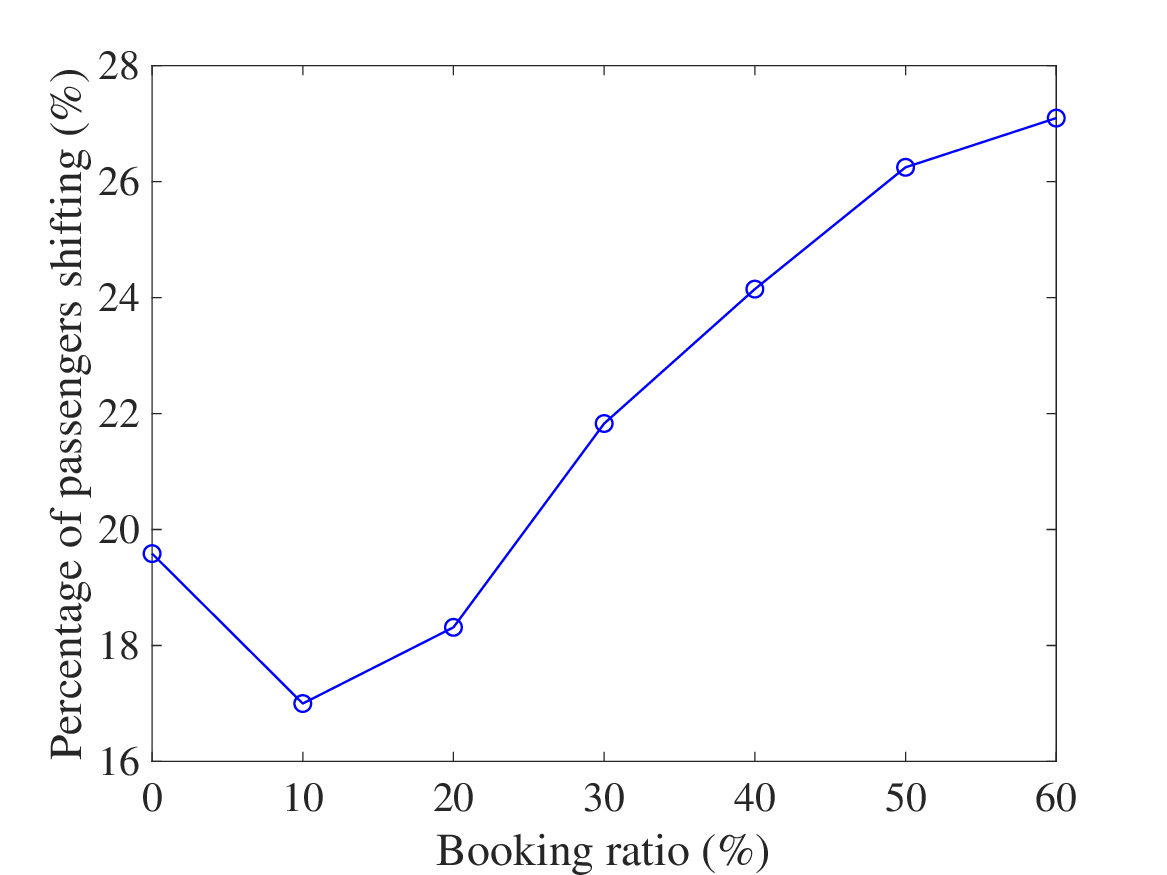}
         \caption{Shift passengers}
     \end{subfigure}
        \caption{Results among various booking ratios.}
        \label{fig:sensitivityFairness}
\end{figure}

In Figure~\ref{fig:sensitivityFairness}, we present the results among various booking ratios, under the operational parameters of 11 trains, the discount of 20\%, a service fairness factor of 50\%, a maximum allowable shifting time of 15 minutes, and weight coefficients established at 1 and 5. Here, the booking ratio refers to the ratio of booked passengers to all passenger demand. An increase in the total waiting time with higher booking ratios is observed, which can be attributed to the lack of temporal flexibility of passengers with reservations. These passengers, constrained by their booking commitments, are unable to realign their departure times to match the optimized timetable in order to reduce waiting times. This finding is further corroborated by the trends depicted in Figure~\ref{fig:sensitivityFairness}(b), where we observe that the percentage of passengers who have to alter their departure times increase as the booking ratio changes from 10\% to 60\%. This alteration arises as a consequence of the occupancy of capacity by reserved passengers adhering to their predetermined departure times, thus forcing a shift among those without reservations to accommodate the system's operational constraints. Moreover, the model is infeasible when the booking ratio is set to 70\%. This result indicates that when too many bookings are announced for passengers, it not only becomes more unfair to those who do not have reservations, but it also prevents some booked passengers from being able to access platforms directly. To sum up, the aforementioned findings show the complexity of achieving a balance between service fairness and system efficiency, highlighting the necessity for strategic planning in the allocation of reservations within the rail transit systems.

\begin{insight}\label{insight3}

\textit{Substantial decreases in the waiting time of passengers (23.78\%) can be achieved at the cost of very limited increases in the lost revenue of the operator (1.00\%).}
\end{insight}

\begin{figure}[]
    \centering
    \includegraphics[height=0.4\textwidth]{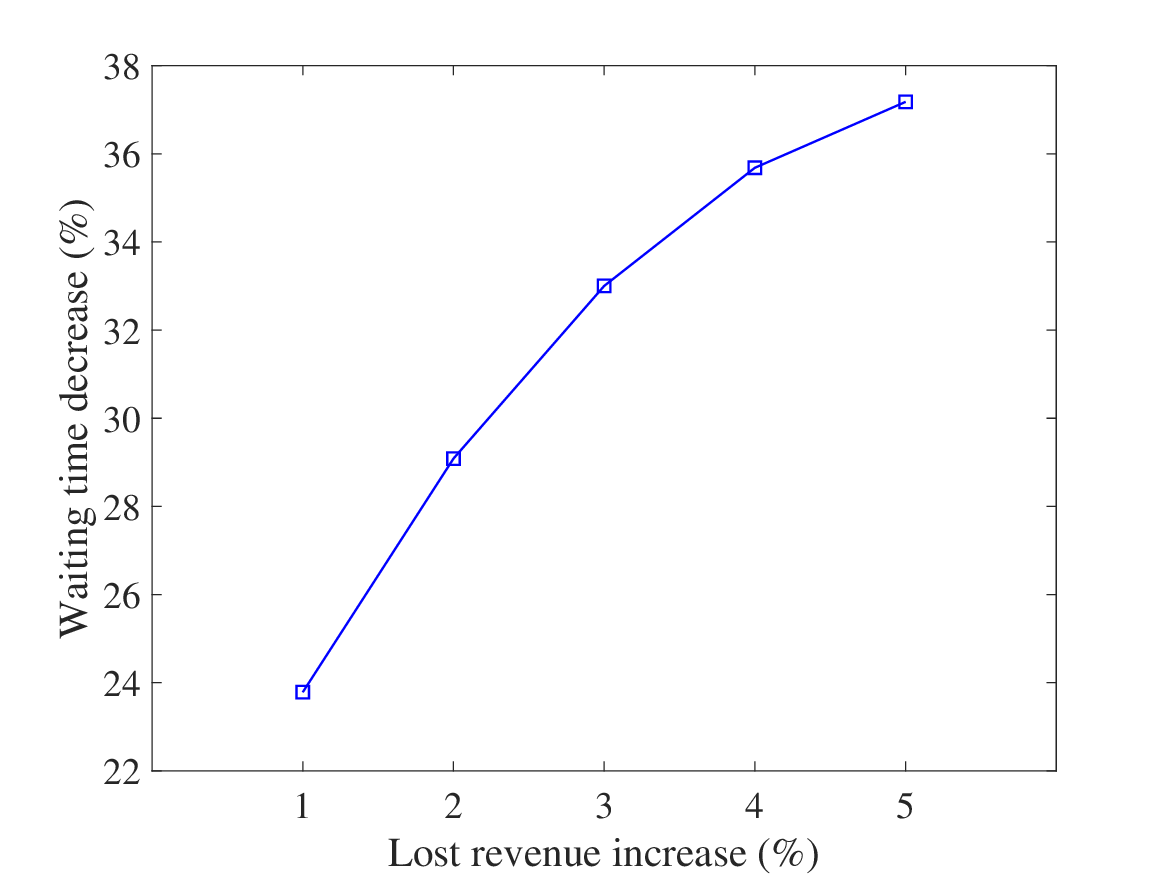}
    \caption{Decrease in waiting time plotted against the increase in the lost revenue, both relative to the solution with the minimum loss of revenue.}
    \label{fig:sensitivityWeight}
\end{figure}

To further analyze the benefit of integrating timetabling, booking, and directing, Figure~\ref{fig:sensitivityWeight} depicts the decrease in waiting time of passengers and the increase in the additional government subsidies for all Pareto solutions, relative to the solution with minimum additional government subsidies. The other motivation behind this analysis is to investigate the trade-off between efficiency and lost revenue which is covered by government subsidies. By starting from the solution that minimizes government subsidies and progressively increasing the upper bound of subsidies, we aim to assess the marginal benefits in terms of reduced passenger waiting time.

To be specific, the objective function is modified as $F^{s}$ to generate the solution with minimum loss of revenue (i.e., the additional government subsidies), resulting in a reduction of $0.51\%$ relative to the original revenue computed as $\sum\limits_{u\in \mathcal{S}}\sum\limits_{v \in \mathcal{S}_{u+1}}\sum\limits_{t \in \mathcal{T}} D_{uvt} \varepsilon_{uv}$. Then, the objective function is modified as $F^{t}$ and a new constraint $F^{s} \leq \varepsilon$ is added. Here, $\varepsilon = (0.51\% + \xi) \times \sum\limits_{u\in \mathcal{S}}\sum\limits_{v \in \mathcal{S}_{u+1}}\sum\limits_{t \in \mathcal{T}} D_{uvt} \varepsilon_{uv}$, where $\xi$ represents the increment of the allowed lost revenue. Finally, we add the allowed lost revenue with a step size of 1\%, and the Pareto optimal points are obtained, as shown in Figure~\ref{fig:sensitivityWeight}. It turns out that a strategic decision to accept a slight 3\% rise in the additional government subsidies can result in a pronounced 33.00\% reduction in waiting time. This result indicates the operational efficiency can be enhanced considerably with a relatively minor increase in the additional government subsidies from the operator’s perspective. All in all, this is a meaningful observation for operating companies, as they can slightly increase the government subsidies to significantly enhance the people's sense of well-being in traveling through public transportation systems. In the longer-term benefit, savings in passengers' waiting time have the potential to translate to higher user satisfaction and increased ridership.

\subsubsection{Performance of all variants of solution methods at the root node.}
\label{sec:rootNode}
To examine the effectiveness of the proposed algorithm, we define the following six variants: 

 (i) \textbf{BD} uses the optimality and feasibility cuts.

(ii) \textbf{TCBD} extends BD by including the cut loop stabilization at the root node and the tailing off strategy.

(iii) \textbf{TTCBD} extends TCBD by including the tree search strategy.

(iv) \textbf{TTSCBD} extends TTCBD by including the strengthened optimality cuts.

(v) \textbf{TRTCBD} extends TTCBD by including the restart strategy.

(vi) \textbf{TRTSCBD} extends TTSCBD by including the restart strategy.

Prior to assessing the full performance of the variants, we first examine the impact that the strategies embedded in the solution methods have on the lower and upper bounds at the root node. In this set of experiments, the booking ratio and the serving fairness factor of passengers without reservations are set as 50\% and 20\%, respectively. Additionally, the discount on the ticket price is set at 20\%. The weight coefficients in the objective function are assigned values of 1 and 10, respectively.
\begin{table}[]
\caption{Lower and upper bounds at the root node for all variants of the solution method.}
\label{tab:rootNode}
\resizebox{\textwidth}{!}{%
\begin{tabular}{crrrrrrr}
\Xhline{1pt}
\begin{tabular}[c]{@{}c@{}}Instance index\\ ( \# of trains,  \# of maximum shifting \\ timestamps, \# of timestamps \\ during peak hours)\end{tabular}  & Root Node          & BD  & TCBD & TTCBD & TTSCBD & TRTCBD & TRTSCBD \\
\Xhline{0.6pt}
\multirow{3}{*}{A (12, 5, 35)}  & Lower Bound   & 20,492.25  & 20,424.86  & 21,408.66& 21,809.93 & 21,414.98    & 21,829.58 \\
                    & Upper Bound  & $-$ & $-$   & 21,847.66 & 22,225.62 & 21,847.66  & 21,847.66\\
                   & Root Node Gap (\%) & $-$ & $-$ & 2.01&1.87 & 1.98 & 0.08\\
\Xhline{0.6pt}
\multirow{3}{*}{B (12, 15, 35)}  & Lower Bound   & 19,931.46  & 19,871.13  & 21,593.46 & 21,524.36 & 21,562.38    & 21,846.26\\
                    & Upper Bound  & $-$ & $-$  & 21,601.12 & 21,601.12 & 21,601.12 &21,853.24\\
                   & Root Node Gap (\%) & $-$ & $-$ & 0.03 & 0.18 & 0.18 &0.03\\
\Xhline{0.6pt}
\multirow{3}{*}{C (15, 5, 35)}  & Lower Bound   & 12,428.30  & 12,730.73  & 17,093.40 & 17,380.99 & 17,452.86    & 17,477.00 \\
                    & Upper Bound  & $-$ & $-$ & 17,501.00 & 17,477.00 & 17,477.00  & 17,477.00\\
                   & Root Node Gap (\%) & $-$ & $-$ & 2.33 & 0.55 &  0.14 &0 \\
\Xhline{0.6pt}
\multirow{3}{*}{D (18, 5, 35)}  & Lower Bound   & 0  & 12,056.10  & 14,573.57 & 14,805.82 & 15,056.99& 15,028.79     \\
                    & Upper Bound  & 11,914.78 & 18,107.00  & 15,334.00 & 15,334.00 & 15,334.00 & 15,334.00  \\
                   & Root Node Gap (\%) & 100.00 & 33.42  & 4.96  &3.44 &  1.81 & 1.99   \\
\Xhline{0.6pt}
\multirow{3}{*}{E (22, 5, 20)}   & Lower Bound   & 13,096.65  & 13,123.82  & 13,155.36 & 13,160.19 & 13,171.69    & 13,183.00 \\
                    & Upper Bound  &  13,183.00 & 13,183.00  & 13,183.00 & 13,199.00 &  13,183.00 & 13,183.00\\
                   & Root Node Gap (\%) & 0.66 & 0.45  & 0.21 & 0.29 & 0.09  & 0 \\
\Xhline{0.6pt}
\multirow{3}{*}{F (22, 5, 30)}  & Lower Bound   & 13,114.52  & 13,102.33  &  13,180.85 &13,182.03 & 13,182.53    & 13,144.14 \\
                    & Upper Bound  & 13,265.00  & 12,236.00  &  13,211.00 & 13,265.00  & 13,227.00  & 13,183.00 \\
                   & Root Node Gap (\%) & 1.13 & 1.01  &  0.23 & 0.63  & 0.34  & 0.29 \\
\Xhline{0.6pt}
\multirow{3}{*}{G (22, 5, 40)}  & Lower Bound   & 13,072.24  & 13,057.69  & 13,173.06  & 13,183.00&  13,141.95   & 13,182.37 \\
                    & Upper Bound  & 13,265.00 & 13,550.00  & 13,265.00 & 13,198.00 & 13,218.00  & 13,198.00 \\
                   & Root Node Gap (\%) & 1.45 & 3.63  & 0.69 & 0.11 & 0.58  &  0.12\\
\Xhline{1pt}
\end{tabular}%
}
\end{table}

In Table~\ref{tab:rootNode}, we present the lower and upper bounds at the root node, as well as the \textit{Root node Gap} at the root node among instances characterized by varying numbers of trains, timestamps, and maximum values for shifting timestamps. The Root node Gap represents the relative difference between lower and upper bounds, calculated using the formula 
\begin{equation*}
\big[(\text{Upper bound} - \text{Lower bound})/ \text{Lower bound} \big] \times 100 ~ (\%).
\end{equation*}
From the results in Table~\ref{tab:rootNode}, the following three observations emerge. First, when comparing TRTCBD to BD, there is a substantial decrease in the root node gap. For example, it can be seen that the root node gap decreases from 100\% under BD to 3.44\% under TTSCBD at Instance D. Second, on average, the strengthened optimality cuts contribute to tightening both the lower bound and the gap at the root node. In particular, in Instances A, B, C, and E, TRTSCBD can find the optimal solution at the root node. Lastly, we observe the benefits of incorporating the tree search strategy, as it improves the solution quality at the root node in all instances when comparing TTCBD with TCBD.

\subsection{Real-world case study}
\label{sec:BeijingMetro}

\begin{figure}[]
\centering
\includegraphics[height=0.4\textwidth]{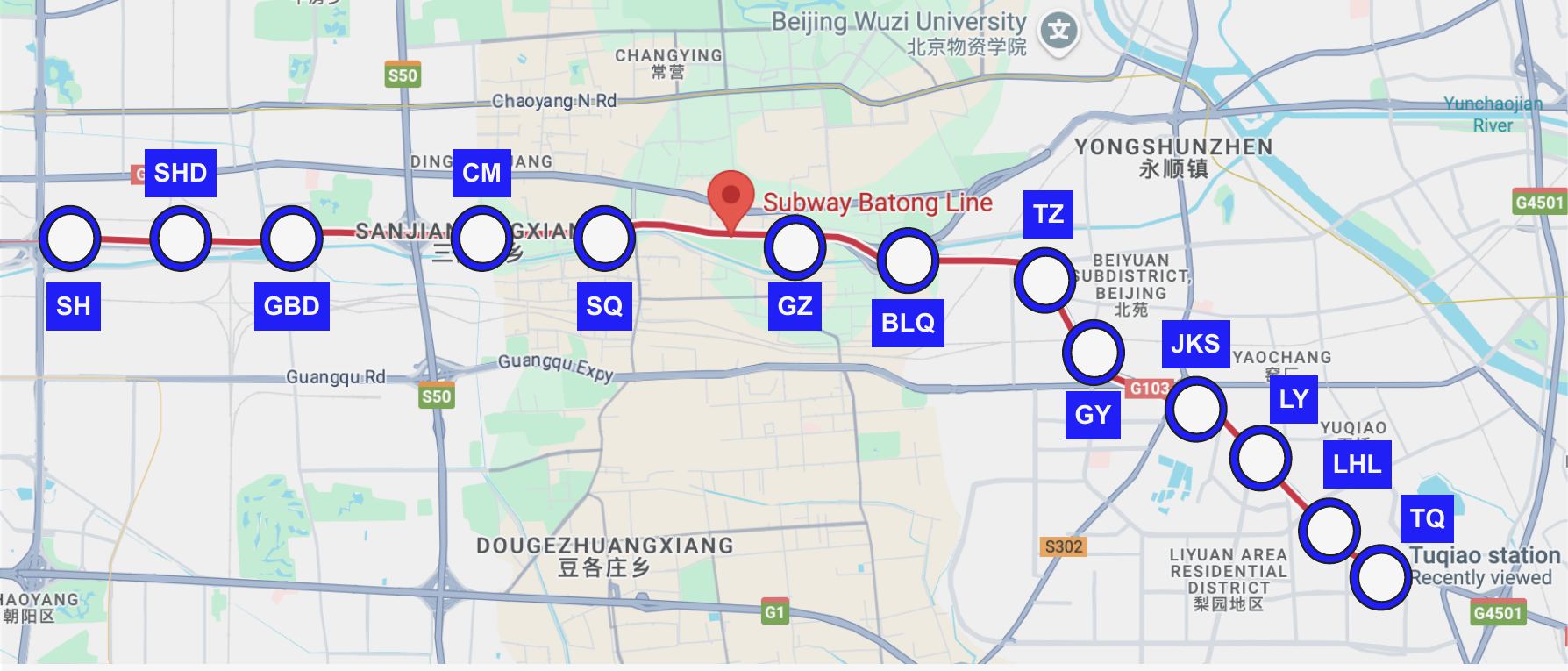}
\caption{An illustration of the Beijing metro Batong line (Source: Google map).}\label{fig 13a}
\end{figure}

The instances used for the experiments are derived from the Beijing metro Batong line, which has 13 stations as depicted in Figure \ref{fig 13a}. This metro line serves as a feeder line, transporting commuters from the suburban district to the city center during morning peak hours. In 2018, nine stations on this line have implemented routine passenger flow control strategies on weekdays to cope with the high volume of passenger demand. Time-dependent passenger demand at each time and each station is derived from historical Automatic Fare Collection (AFC) data. We follow the method described in Section~\ref{sec:smallCase} to preprocess the demand for both reserved and non-reserved passengers. Tables \ref{Table13a} and \ref{price-real} provide detailed information about running times on sections, the dwell time at each station, and the distance-based ticket prices used in practice, respectively. On this basis, we consider five instances with varying study time horizons, as detailed in Table \ref{tab:instance}. To construct these instances, continuous time periods are discretized into one-minute timestamps. We also follow the conclusions derived from the proof-of-concept experiments to set the instance characteristics. For example, when passengers have greater flexibility, such as the ability to shift their trips during extended peaking hours, fewer trains are needed. Besides, the maximum and minimum headways are set to 360 and 120 seconds, respectively. The maximum number of passengers that a train can accommodate is 2,000. In addition, a discount of 20\% is offered to passengers who are willing to shift their travel times.  

\begin{table}[]
\centering
\caption{Dwell and running times on the Beijing metro Batong line.}\label{Table13a}
     \resizebox{\textwidth}{!}{
\begin{tabular}{ccccc}
\hline
Station index & Station name (Abbreviation) & Dwell time (in: second) & Section & Running time (in: second) \\
\hline\noalign{\smallskip}
1&Tuqiao (TQ) & 60 &Tuqiao$\rightarrow$ Linheli & 120 \\
2&Linheli (LHL)&60 & Linheli $\rightarrow$ Liyuan & 120 \\
3&Liyuan (LY) &60 & Liyuan $\rightarrow$ Jiukeshu & 120\\
4& Jiukeshu (JKS)&60 &Jiukeshu$\rightarrow$ Guoyuan& 120 \\
5&Guoyuan (GY)&60 &Guoyuan$\rightarrow$ Tongzhou& 120 \\
6&Tongzhou (TZ)&60 &Tongzhou$\rightarrow$ Baliqiao& 180\\
7&Baliqiao (BLQ)&60 &Baliqiao$\rightarrow$ Guanzhuang& 180 \\
8&Guanzhuang (GZ) &60 &Guanzhuang$\rightarrow$ Shuangqiao& 180 \\
9&Shuangqiao (SQ)&60 &Shuangqiao$\rightarrow$ Chuanmei Uni& 180 \\
10&Chuanmei Uni (CM)&60 &Chuanmei Uni$\rightarrow$ Gaobeidian& 180 \\
11&Gaobeidian (GBD)&60 &Gaobeidian$\rightarrow$ Sihui East& 120 \\
12&Sihui East (SHE)&60 &Sihui East$\rightarrow$ Sihui&180 \\
13&Sihui (SH)&60 & &\\
\hline
\end{tabular}}
\end{table}

\begin{table}[]
\centering
\caption{The distance-based ticket price of Beijing metro Batong Line (unit: RMB) }
\label{price-real}
\begin{tabular}{cccccccccccccc}
\hline
Station & 1 & 2 & 3 & 4 & 5 & 6 & 7 & 8 & 9 & 10 & 11 & 12 & 13 \\ \hline
1       & 0 & 3 & 3 & 3 & 3 & 3 & 4 & 4 & 4 & 4  & 5  & 5  & 5  \\
2       & 0 & 0 & 3 & 3 & 3 & 3 & 4 & 4 & 4 & 4  & 5  & 5  & 5  \\
3       & 0 & 0 & 0 & 3 & 3 & 3 & 3 & 4 & 4 & 4  & 5  & 5  & 5  \\
4       & 0 & 0 & 0 & 0 & 3 & 3 & 3 & 3 & 4 & 4  & 4  & 5  & 5  \\
5       & 0 & 0 & 0 & 0 & 0 & 3 & 3 & 3 & 3 & 4  & 4  & 5  & 5  \\
6       & 0 & 0 & 0 & 0 & 0 & 0 & 3 & 3 & 3 & 4  & 4  & 4  & 5  \\
7       & 0 & 0 & 0 & 0 & 0 & 0 & 0 & 3 & 3 & 3  & 4  & 4  & 4  \\
8       & 0 & 0 & 0 & 0 & 0 & 0 & 0 & 0 & 3 & 3  & 3  & 4  & 4  \\
9       & 0 & 0 & 0 & 0 & 0 & 0 & 0 & 0 & 0 & 3  & 3  & 3  & 4  \\
10      & 0 & 0 & 0 & 0 & 0 & 0 & 0 & 0 & 0 & 0  & 3  & 3  & 3  \\
11      & 0 & 0 & 0 & 0 & 0 & 0 & 0 & 0 & 0 & 0  & 0  & 3  & 3  \\
12      & 0 & 0 & 0 & 0 & 0 & 0 & 0 & 0 & 0 & 0  & 0  & 0  & 3  \\
13      & 0 & 0 & 0 & 0 & 0 & 0 & 0 & 0 & 0 & 0  & 0  & 0  & 0  \\ \hline
\end{tabular}
\end{table} 

\begin{table}[]
    \centering
\caption{Characteristics of the instances used in numerical experiments.}
    \label{tab:instance}
    \begin{tabular}{crccr}
\Xhline{1pt}
        Instance & $\left| \mathcal{T} \right|$ & The whole time period & \begin{tabular}[c]{@{}c@{}}The period during \\ peaking hours\end{tabular}  & \# of trains  \\
\Xhline{0.6pt}
        H        & 60  & 6:30 - 7:30 & 6:35 - 7:30  & 6  \\
        I        & 75  & 6:30 - 7:45 & 7:00 - 7:45  & 14 \\
        J        & 90  & 6:30 - 8:00 & 7:00 - 8:00  & 18  \\  
        K        & 100  & 6:50 - 8:30 & 7:00 - 8:30  & 24 \\ 
        L        & 120 & 6:30 - 8:30 & 7:00 - 8:30 & 33 \\
        M        & 100  & 6:50 - 8:30 & 6:50 - 8:30 & 23 \\
        N        & 120 & 6:30 - 8:30 & 6:50 - 8:30& 32 \\
        O        & 120 & 6:30 - 8:30 & 6:40 - 8:30 & 31 \\
    \hline
\end{tabular}
\end{table}

Based on the prepared data, four sets of numerical experiments are conducted to evaluate the performance of the proposed approaches. First, Section \ref{sec:solutionMethod} compares the six variants of our BD algorithm with GUROBI to validate the solution efficiency of our algorithm and identify the most effective variant for subsequent experiments. Second, Section \ref{sec:performanceAlgorithms} evaluates the performance gains of our algorithm over a heuristic approach and the traditional BD method. Third, Section \ref{sec:QoS} evaluates the benefits of integrating timetabling, trip shifting, and passenger flow control. Finally, Section \ref{sec:sensitivity} presents a sensitivity analysis of key parameters.

\subsubsection{Performance comparison between GUROBI and six variants of the BD algorithm.}
\label{sec:solutionMethod}

In this section, we quantify the advantages of our TTCBD compared to GUROBI and the other five variants, and derive insights from an algorithmic perspective (Insight \ref{insight4}).

\begin{insight}\label{insight4}

\textit{TTCBD is the most effective variant of our proposed Benders decomposition algorithm, outperforming GUROBI and the other five variants in large-scale instances (e.g., Instances J, K, and L).}

\end{insight}

We first evaluate the full performance of the six variants of the proposed BD algorithm based on the data of Instances H, I, J, K, L, as shown in Table \ref{tab:instance}. The coefficient values $\omega_s$ and $\omega_t$ are set to 1 and 10, whose impacts will be discussed in Section \ref{sec:sensitivity}. In Table~\ref{tab:fullPerformance}, we present the objective function values (abbreviated as Obj.) and the computational times (abbreviated as Time) obtained from GUROBI and six variants of solution methods within a Gap limit of 5\%. It can be seen that the cut loop, tailing off, and tree search strategies are able to reduce the computational time considerably. For example, GUROBI necessitates approximately 5,430 seconds to reach a near-optimal solution with a Gap of 5\% for instance L. When employing TTCBD, however, a better solution is found within 2,200 seconds. We can also observe that including the strengthened optimality cuts worsens the solution efficiency. This can be explained by the fact that an integer programming model~\eqref{eq:strengthenedoptimality} has to be solved to generate these cuts, which introduces additional computational burden.
\begin{table}[]
\caption{Performance for GUROBI and the six variants of algorithms with a Gap limit of 5\%}
\label{tab:fullPerformance}
\resizebox{\textwidth}{!}{%
\begin{tabular}{crrrrrrrr}
\Xhline{1pt}
\begin{tabular}[c]{@{}c@{}}Instance \\ index\end{tabular}  &    & GUROBI & BD & TCBD & TTCBD & TTSCBD & TRTCBD & TRTSCBD \\
\hline
\multirow{2}{*}{H} & Obj. & 13,229,00  & 13,229,00  & 13,229,00 & 13,229,00   & 13,229,00 & 13,229,00 & 13,229,00 \\
& Time (sec.)  & 1.48 & 4.07 & 3.84  & 4.42 & 6.24  & 9.43 &  15.05 \\     
\multirow{2}{*}{I} & Obj. & 41,307.00  & 41,307.00  &  41,307.00&  41,307.00  &  41,307.00  &  41,307.00     &  41,307.00 \\
& Time (sec.)  & 56.37 & 168.78 & 139.08 & 124.31 & 399.64 & 223.68 & 1044.36  \\
\multirow{2}{*}{J} & Obj. & 43,042.00  & 43,128.00 & 43,098.00 & 43,042.00 &  43,463.00 &  43,042.00    & 43,080.00 \\
& Time (sec.)  & 233.55 & 303.83 & 292.16  & \textbf{175.20}& 524.77 & 404.97 & 1316.53   \\ 
\multirow{2}{*}{K} & Obj. & 59,543.00  & 59,377.00 & 59,377.00 &  59,377.00 & 59,377.00  &   59,377.00   &  59,377.00\\
& Time (sec.)  & 910.43 & 2,185.19  & 22,39.00 & \textbf{895.33} & 1,240.36 & 1,678.77 &  3,491.24  \\ 
\multirow{2}{*}{L} & Obj. & 72,010.00  & 71,945.00 & 71,945.00 & 71,945.00 & 71,945.00  &  72,081.00    & 71,945.00 \\
& Time (sec.)  & 5,428.54 & 4,176.90 & 4,006.74  & \textbf{2,192.77} & 5,859.72 & 7,524.54 &  1 5,800.98 \\ 
\Xhline{1pt}
\end{tabular}%
}
\end{table}

\begin{figure}[]
     \centering
 \begin{subfigure}{\textwidth}
         \centering
\includegraphics[width=0.6\textwidth]{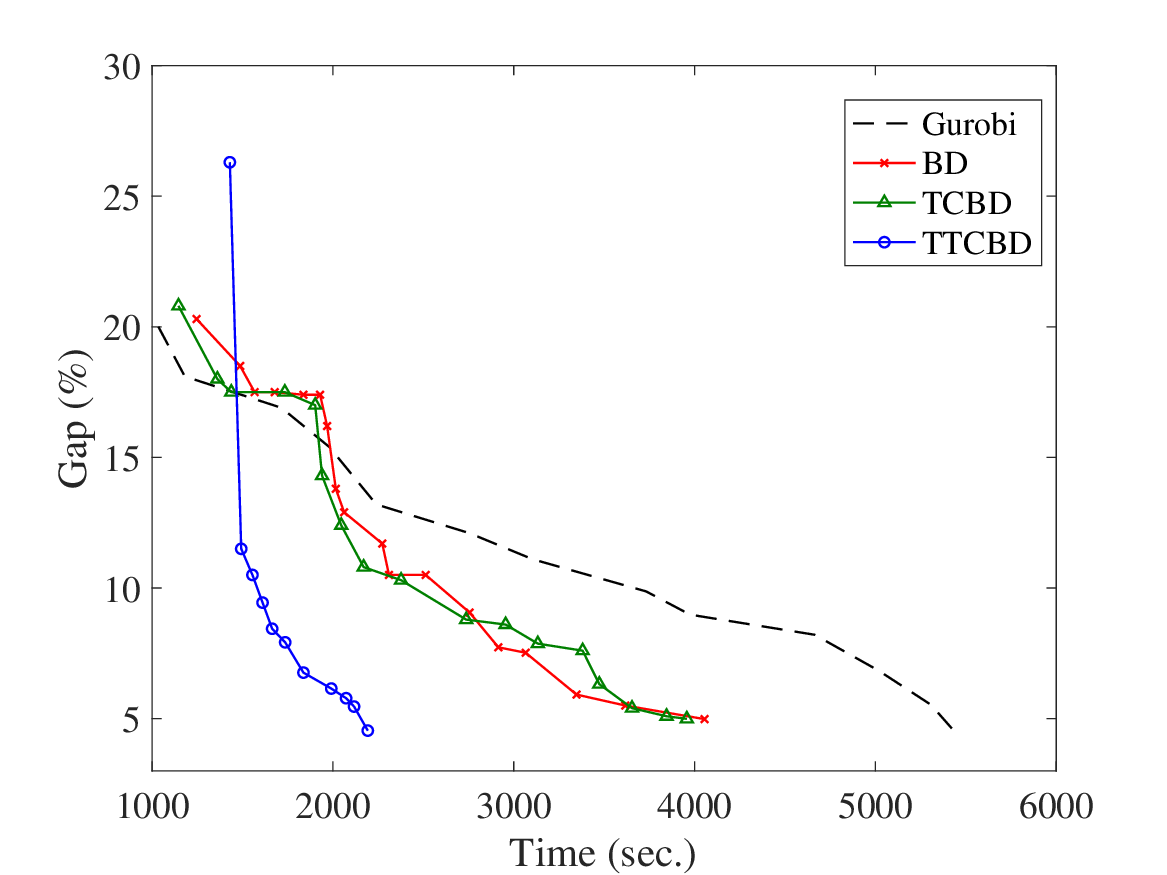}
         \caption{Gap}
     \end{subfigure}
     
     \begin{subfigure}{\textwidth}
         \centering
 \includegraphics[width=0.6\textwidth]{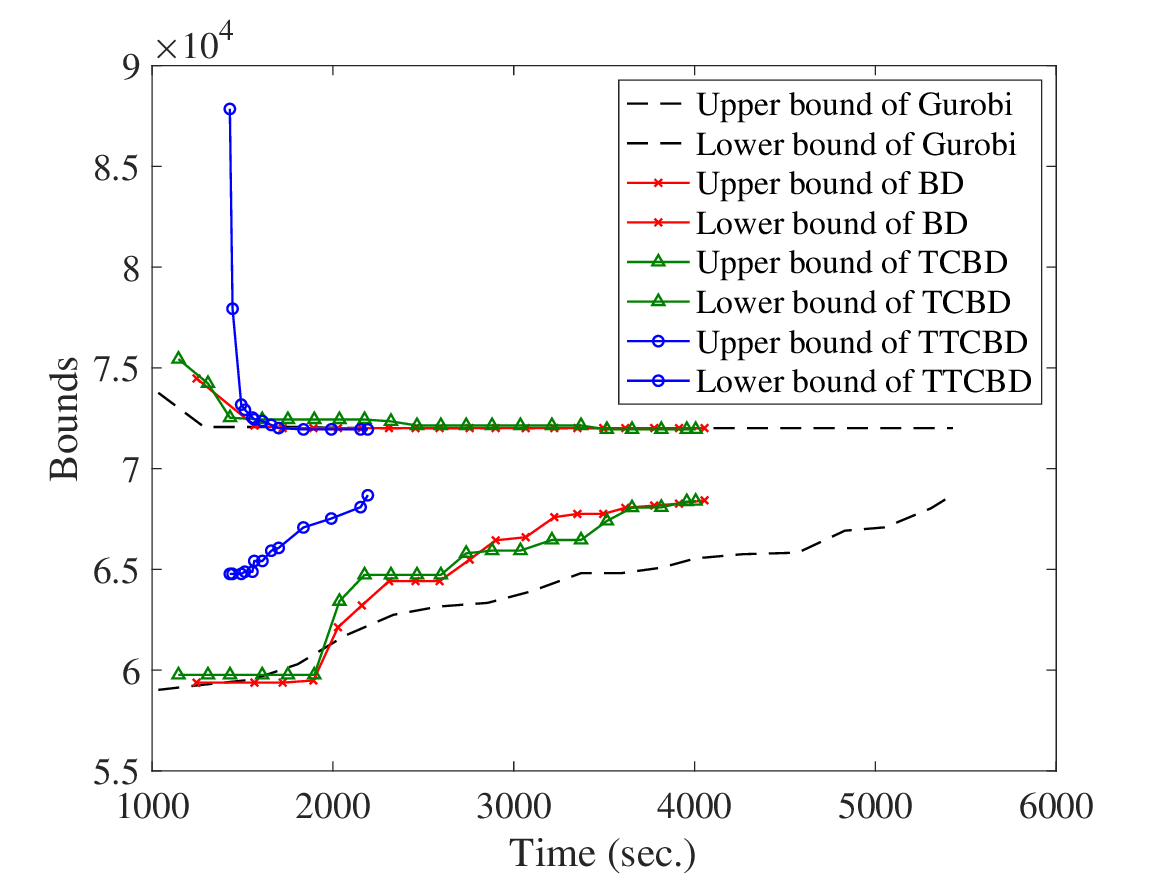}
         \caption{Upper and lower bounds}
     \end{subfigure}
        \caption{Convergence trends of the upper and lower bounds and the Gap with respect to Instance L with a Gap limit of 5\%.}
        \label{fig:convergence}
\end{figure}

Figure~\ref{fig:convergence} presents an overall comparison among GUROBI, BD, TCBD, and TTCBD in terms of the convergence trends of upper and lower bounds and the optimality gap over time with respect to the Instance L. Looking to Gap displayed in Figure~\ref{fig:convergence}(a), it is evident that TTCBD outperforms the other algorithms without the tree search strategy. Specifically, TTCBD finds a solution with a Gap of less than 5\% within 2,200 seconds. In comparison, both TCBD and BD reach solutions with similar quality after 4,000 seconds, whereas GUROBI does not converge to a comparable solution within 5,000 seconds. Figure~\ref{fig:convergence}(b) provides detailed results of upper and lower bounds over time for these solution methods. Notably, the initial lower bound obtained by TTCBD is highest. The significant improvement in performance of TTCBD can be attributed to the incorporation of a tree search strategy, which facilitates strong branching and in-depth exploring for the lower bound. The second observation is that GUROBI, BD, and TCBD can quickly find a good upper bound, while the lower bound slowly increases. All in all, this experiment demonstrates that TTCBD outperforms the other alternatives. Therefore, we employ TTCBD for the subsequent experiments.

\subsubsection{ Performance comparison of different solution approaches.}
\label{sec:performanceAlgorithms}

In this section, we compare the performance of our proposed TTCBD algorithm with several benchmarks in terms of generating optimal solutions: (i) a hybrid algorithm combining a local search method with GUROBI (referred to as Hybrid-LS), and (ii) the traditional Benders decomposition algorithm introduced in Appendix \ref{sec:traditional}, where the RMP only determines the timetables (referred to as Traditional BD). The traditional Benders decomposition algorithm also incorporates cut loop stabilization at the root node, the tailing-off strategy, and the tree search strategy, consistent with TTCBD, to ensure a fair comparison. In addition, the pseudo-code of the hybrid algorithm is presented in Appendix \ref{sec:local}. In this hybrid solution framework, we decompose the problem into a timetabling subproblem and a passenger-related subproblem, whose key decomposition idea is the same as the decomposition method proposed in the traditional Benders decomposition algorithm in Appendix \ref{sec:traditional}. The local search algorithm is used to solve the timetabling subproblem with timetabling-related variables and constraints, and GUROBI is employed to solve the passenger-related subproblem with the fixed timetable generated by the timetabling subproblem. The solution quality of the generated timetables is evaluated by the passenger-related subproblem.

In this experiment, the termination criteria for the two exact solution methods (i.e., TTCBD and traditional BD) are set as follows: (i) the optimality gap reaches zero; or (ii) the computational time exceeds 7,200 seconds. The hybrid algorithm (i.e., Hybrid-LS) terminates when the computational time exceeds 7,200 seconds. The goals of the performance comparison between our proposed TTCBD with the aforementioned benchmarks are to address the following questions:

(i) To what extent does TTCBD outperform Hybrid-LS, considering that the latter may get stuck in local optima?

(ii) To what extent does TTCBD outperform the traditional BD, where the RMP only determines the timetables? In other words, what is the value of our novel decomposition framework, which incorporates partial passenger-related variables into the RMP and integrates the valid equalities and inequalities \eqref{eq:validBoard}–\eqref{eq:wildeb}?

Table \ref{tab:algorithm} presents the detailed results of our TTCBD, Traditional BD, and Hybrid-LS algorithms, reporting the upper and lower bounds, optimality gaps, and computational times for each solution method. Since Hybrid-LS cannot generate a lower bound, we use the best lower bound obtained from the other two algorithms to compute the optimality gap for this hybrid approach. We first focus on analyzing the performance of TTCBD and Hybrid-LS, and then compare TTCBD with Traditional BD. We observe that even for the smallest instance among the four, i.e., Instance J, Hybrid-LS can only find a solution with a 10.45\% optimality gap. This can be attributed to two main reasons. First, as a heuristic method, it tends to converge to local optima and lacks mechanisms to escape them. Second, this algorithm simply generates candidate timetables to satisfy the timetabling subproblem without incorporating any passenger information into the search process. These results highlight the advantages of using exact methods to solve the investigated problem.

\begin{table}[]
\caption{Performance comparison of our TTCBD, Traditional BD, and Hybrid-LS.}
\label{tab:algorithm}
\begin{tabular}{ccrrrr}
\Xhline{1pt}
Instance           & Solution method & Upper bound & Lower bound & Optimality gap (\%) & Time (s) \\
\Xhline{0.6pt}
\multirow{3}{*}{J} & TTCBD &   43,042.00                                                     & 43,042.00                                                        &  0                                                              & 282.02                                                                            \\
                   & Traditional BD                                             &    53,773.00                                                    &  0                                                      &                100.00                                                                                                                   &  7,200.00        \\

                   & Hybrid-LS                                               & 47,540.00                                                       &  -                                                      &  10.45                                                                                                                                & 7,200.00  \\
\Xhline{0.6pt}  
\multirow{3}{*}{M} & TTCBD                      &   61,345.00                             &    61,345.00                                                    &       0                                                 &        2,240.07                                                                                                                                   \\
                   & Traditional BD                                             &   82,004.00                                                     &   -                                                     & 100.00    & 7,200.00              
                  \\
                   & Hybrid-LS    &  68,255.00                                                    & -  &     11.26                                                   &  7,200.00                                                                                                                                \\
\Xhline{0.6pt} 
\multirow{3}{*}{N} & TTCBD &   73,615.00  &  73,615.00  & 0   &  4,681.99 \\
  & Traditional BD  &   102,721.000   &  0   &   100.00   &  7,200.00        \\
                   & Hybrid-LS                                               &    84,234.00                                                                                                          &   -                                                             &         14.43                                                          &  7,200.00 \\
\Xhline{0.6pt} 
\multirow{3}{*}{O} & TTCBD &  75,384.00   &  73,768.01  & 1.61   & 7,200.00  \\
  & Traditional BD  &   109,204.00   &  0   & 100.00     &  7,200.00        \\
                  
                   & Hybrid-LS                                               &    90,208.00                                                                                                          &     0                                                           &       19.66                                                            & 7,200.00  \\
\Xhline{0.6pt} 
\end{tabular}
\end{table}

\begin{figure}[]
     \centering
 \begin{subfigure}{0.49\textwidth}
         \centering
\includegraphics[width=\textwidth]{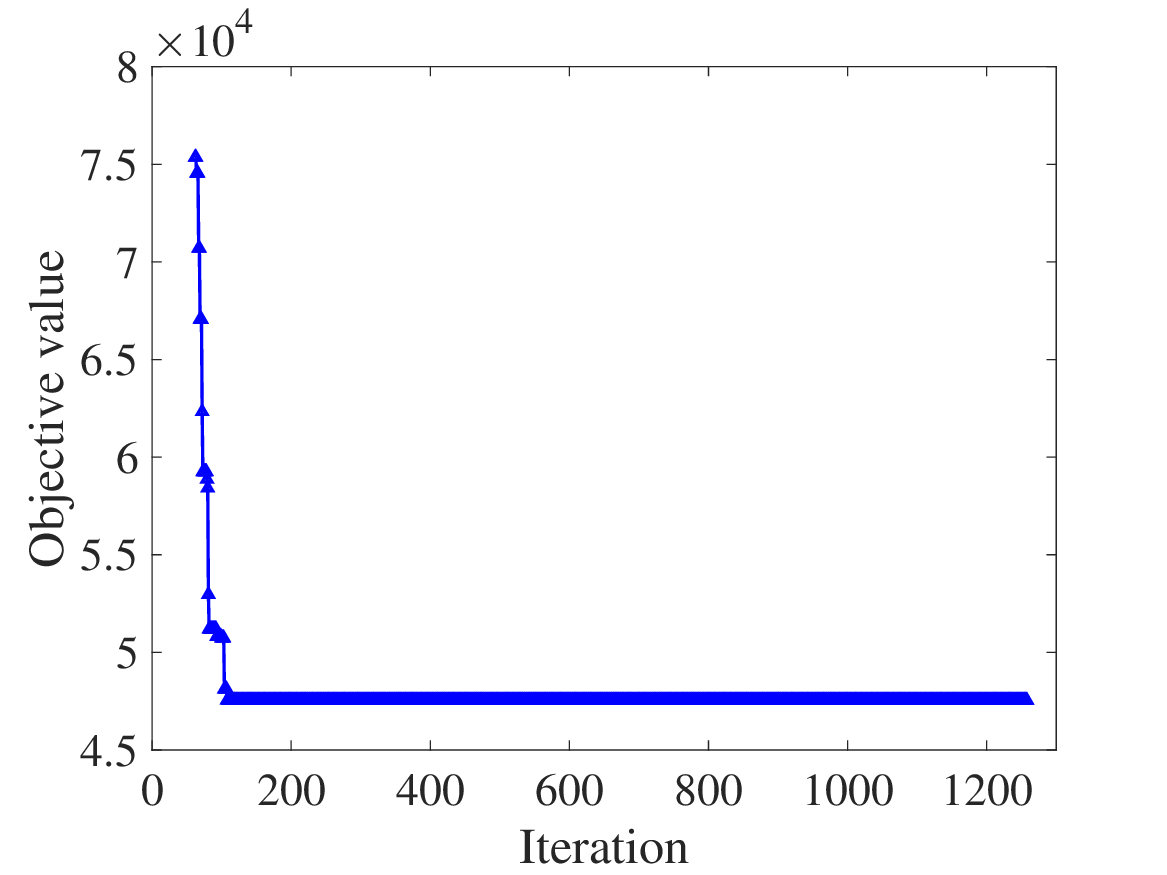}
         \caption{Instance J}
     \end{subfigure}
     \begin{subfigure}{0.49\textwidth}
         \centering
 \includegraphics[width=\textwidth]{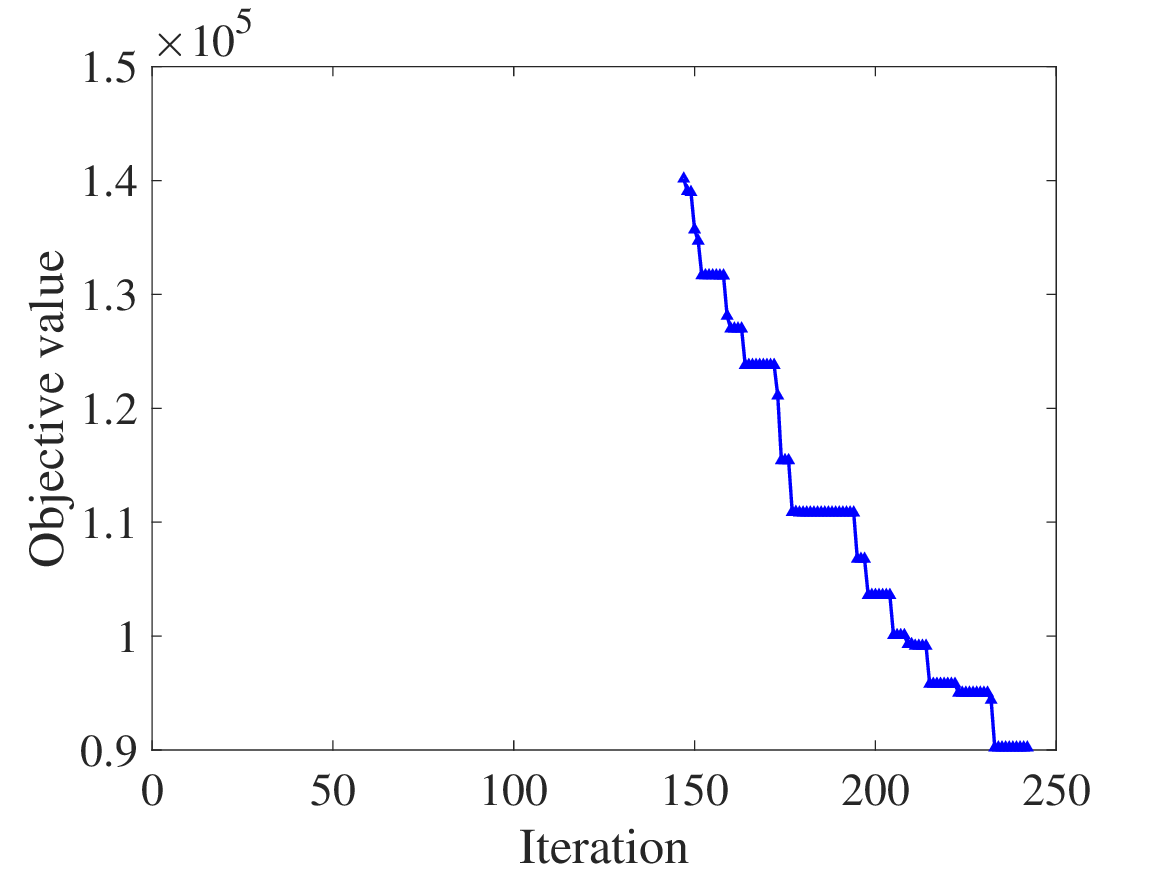}
         \caption{Instance O}
     \end{subfigure}
        \caption{Convergence trend of the hybrid algorithm on the objective value.}
        \label{fig:convergence_hybird}
\end{figure}

Figure \ref{fig:convergence_hybird} shows the convergence trend of Hybrid-LS on the objective value for Instances J and O. The first feasible solution emerges at 63th iteration and 146th iteration, respectively, suggesting that Hybrid-LS initially encountered difficulties in finding a timetable that makes the passenger-related subproblem feasible. These findings indicate that the algorithm cannot easily generate a timetable that satisfies all constraints in the passenger-related subproblem, such as ensuring reserved passengers board the first arriving train, partially accommodating unreserved passengers, respecting train capacity limits, and serving all passengers within the study time horizon. A possible reason is the complexity and strict limitations of the passenger-related subproblem, combined with the lack of feedback from this subproblem to the timetabling subproblem. In other words, no passenger information is integrated into the timetabling process to guide the search. These results indicate the necessity of effectively incorporating passenger information into the solution approach for our integrated optimization problem of demand-side management and timetabling, in order to obtain high-quality timetables.

A second observation is that the final solution is identified by the 109th iteration for Instance J, after which the objective value remains unchanged for the subsequent 1148 iterations, indicating that the algorithm has converged to a local optimum. For the larger-scale Instance O, the final solution is found at the 232nd iteration and remains unchanged for the following iterations, until the computational time limit is reached and the algorithm terminates. These findings highlight the necessity of incorporating passenger information directly into the timetabling process. Such integration helps guide the search more effectively, reduces the time required to identify feasible solutions, and may assist in escaping local optima.

\begin{insight}\label{insight5}

\textit{TTCBD outperforms both the traditional Benders decomposition algorithm and the hybrid algorithm that combines a local search method with GUROBI, in terms of solution efficiency and quality. This advantage is attributed to our novel decomposition framework and the valid equalities and inequalities, which effectively incorporate passenger information into the timetabling process.}
\end{insight}

We now turn to comparing the performance of our TTCBD and Traditional BD and derive insights from an algorithmic perspective. An interesting observation from the results in Table \ref{tab:fullPerformance} and Table \ref{tab:algorithm} is that when TTCBD terminates with a 5\% gap limit for Instance J, the corresponding objective value is exactly the optimal one reported in Table \ref{tab:algorithm}. This indicates that the upper bound obtained under the 5\% optimality gap termination rule is already optimal, even though the lower bound has not yet fully converged. In practical operations, an optimality gap limit of 5\% is often sufficient to ensure a high-quality solution. In our tested instances, this setting performs even better, as one of the solutions obtained under this termination criterion turns out to be optimal. These results further support the use of a 5\% optimality gap as a practical and effective termination criterion in real-world applications.

The second observation from Table \ref{tab:algorithm} is that Traditional BD fails to improve the lower bound within 7,200 seconds, while TTCBD consistently achieves solutions with an optimality gap of at most 1.61\% across all instances under the same time limit. This is because Traditional BD decomposes the problem into a timetabling master problem and a passenger-related subproblem without strong integration between the two. In the Traditional BD framework, many feasibility cuts are added, which are weak. In contrast, TTCBD employs a novel decomposition framework, enhanced with valid equalities and inequalities, that tightly couples passenger demand with the timetabling decision process. Our developed framework contributes to reducing the number of feasibility cuts and lifting the lower bound, thus improving the solution efficiency.

\begin{figure}[]
     \centering
 \begin{subfigure}{0.49\textwidth}
         \centering
\includegraphics[width=\textwidth]{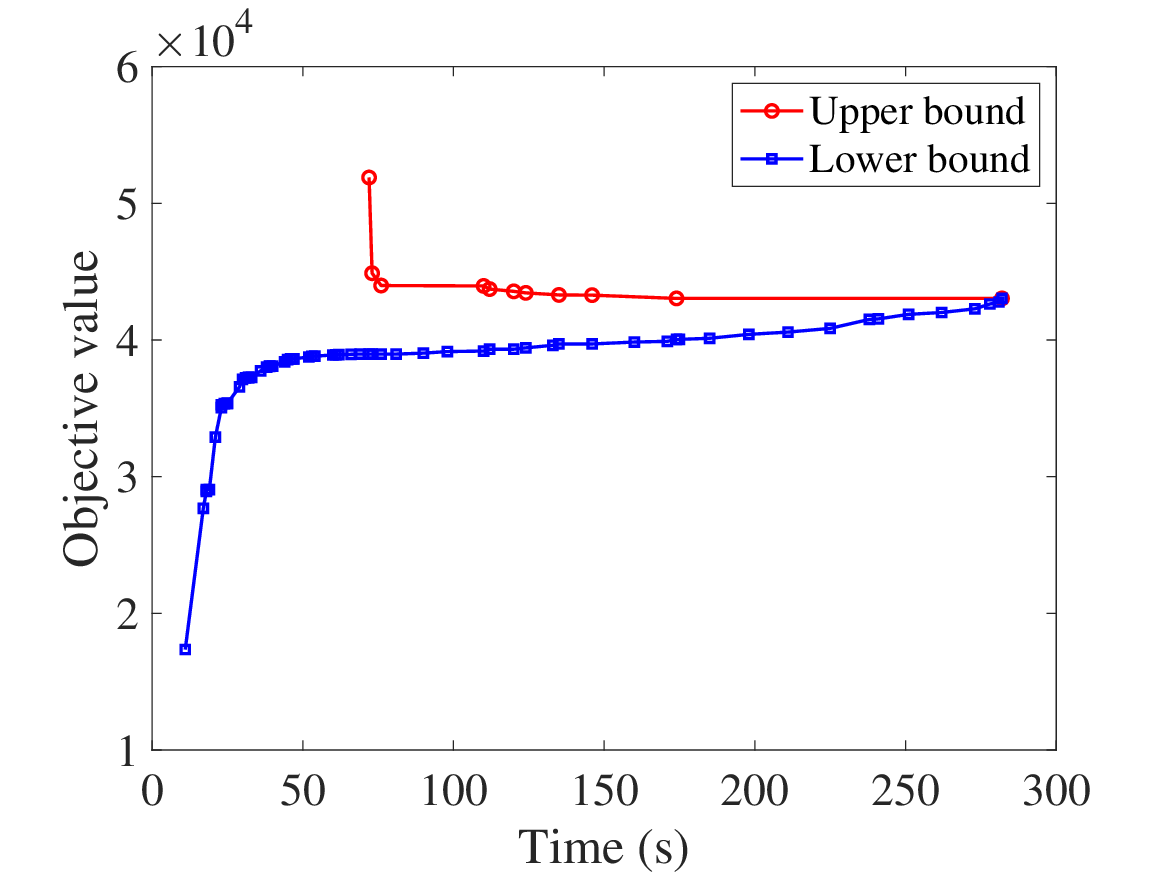}
         \caption{Instance J}
     \end{subfigure}
     \begin{subfigure}{0.49\textwidth}
         \centering
 \includegraphics[width=\textwidth]{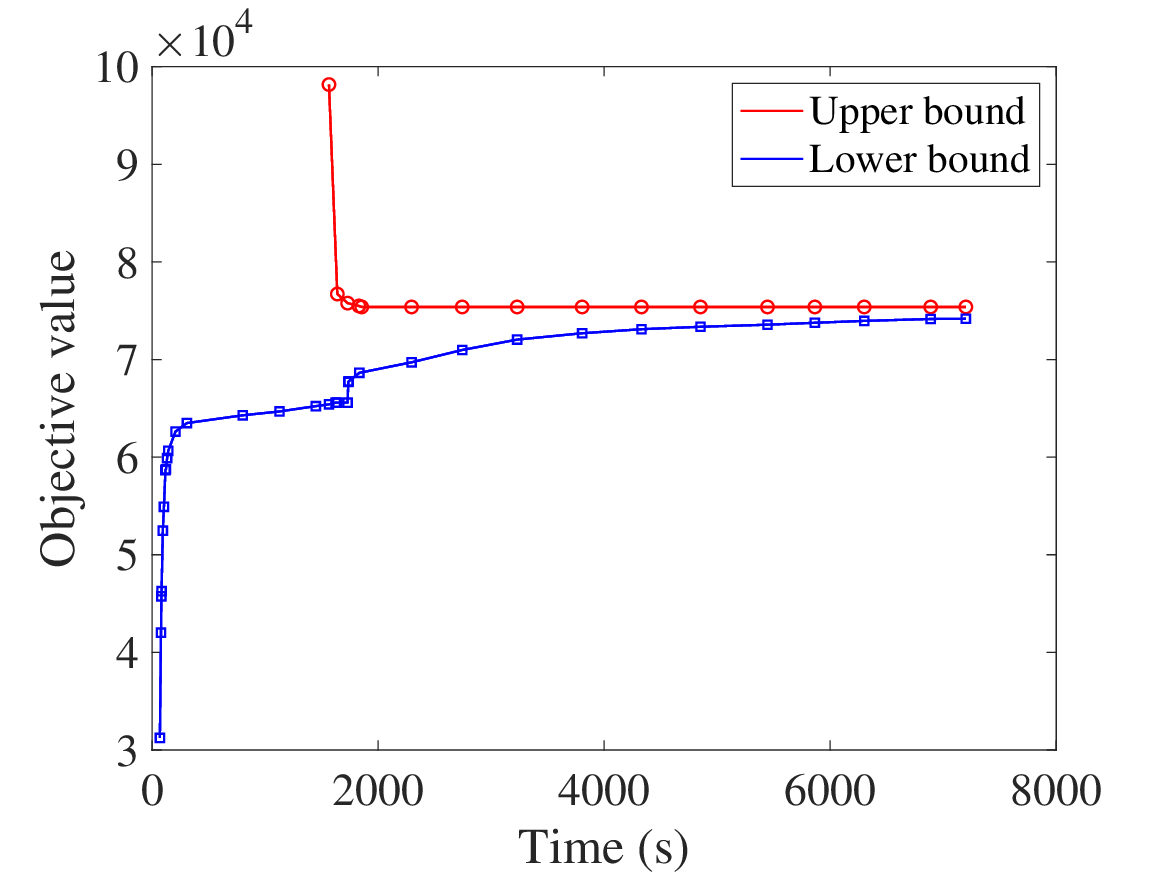}
         \caption{Instance O}
     \end{subfigure}
        \caption{Convergence trend of TTCBD with respect to the upper and lower bounds.}
        \label{fig:convergence_TTCBD}
\end{figure}

Figure \ref{fig:convergence_TTCBD} illustrates the convergence trends of TTCBD for instances J and O. In both small-scale (Instance J) and large-scale (Instance O) cases, the lower bound improves rapidly in the early stages. The upper bound also drops sharply and stabilizes early, which narrows the optimality gap quickly. These trends highlight the practical advantage of our TTCBD. To summarize, the results demonstrate the effectiveness of our decomposition framework armed with valid equalities and inequalities in accelerating convergence and improving solution quality.

\begin{figure}
    \centering
    \includegraphics[width=0.65\linewidth]{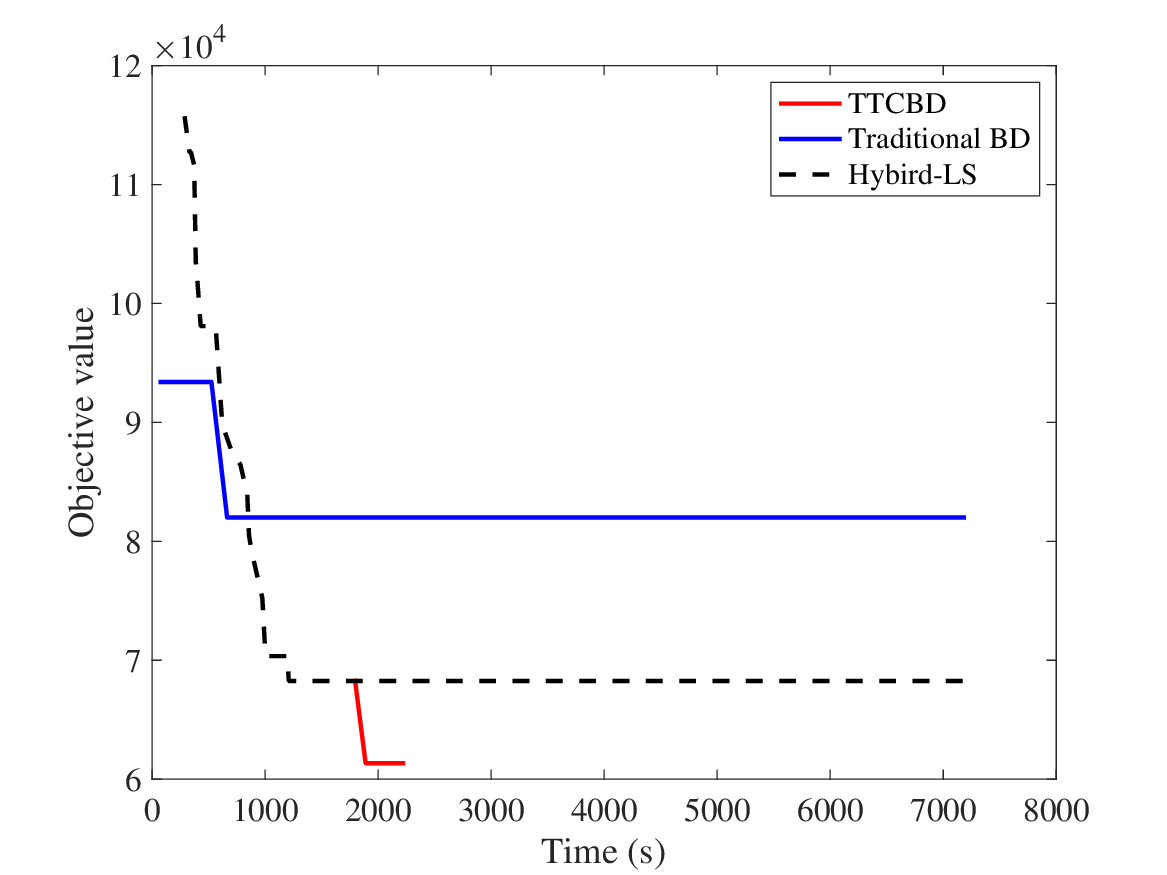}
    \caption{Comparison of upper bounds generated by TTCBD, Traditional BD, and Local search for Instance M.}
    \label{fig:upperBouned}
\end{figure}

Lastly, we examine the performance of the three solution methods in terms of upper bound convergence. Figure \ref{fig:upperBouned} presents the comparison for Instance M. We observe that Traditional BD finds an upper bound the fastest but is unable to improve it after 664 seconds of computation, with a final upper bound of 82,004. Hybrid-LS improves the upper bound to 68,255 after 1,209 seconds. Our TTCBD finds the upper bound more slowly but delivers the best quality. It identifies an upper bound of 68,437 after 1,794 seconds and further refines it to 61,345 after 1,939 seconds. Considering both the computation time efficiency and the best upper bound, we conclude that our developed TTCBD outperforms the other methods.

\subsubsection{Benefits of integrating booking, directing and timetabling.}
\label{sec:QoS}

We next perform experiments to assess the benefits of the integrated BDDT approach. The instance covering the full time period from 6:30 to 8:00, with a peak period from 7:00 to 8:00, is solved while allowing departure time shifts of up to 20 minutes. The weighting coefficients $\omega_t$ and $\omega_s$ are set to 1 and 10.

Based on these settings, we compare our BDDT approach against two alternative approaches. 

(i) \textit{PFC: Optimizing the passenger flow control strategy without allowing for trip shifting under the optimized train timetable through BDDT.}

(ii) \textit{SPE: Simulating the evolution of passengers under the same optimized train timetable.}

Upon inputting the optimized timetable obtained through the BDTT approach into the PFC and SPE models, both models are proven to be infeasible. This result illustrates that, in the absence of trip shifting, the existing number of trains are incapable of satisfying all passenger demand while ensuring service fairness constraints. The results are in line with those of the previous experiment: encouraging passengers to shift their departure times would save the number of operated trains. Moreover, two additional trains are introduced and the timetable is reoptimized. When applied to the PFC and SPE approaches, the results show that PFC is feasible while SPE remains infeasible. The reason is that the PFC method can control the number of boarding passengers at upstream stations, thus preserving capacity for those boarding at later stations. Essentially, this finding highlights the necessity of implementing the passenger flow control strategy to optimize the allocation of capacity resources, thereby efficiently serving all passengers from a system-optimal perspective.

Further, to quantify the effectiveness of the proposed BDDT approach, we construct the other set of experiments. We relax the capacity constraints~\eqref{eq:InvehicleLimitConstraint} in the PFC  model, which is denoted as \textit{relaxed PFC}. In Table~\ref{tab:BDDT}, we present the results under the BDDT and relaxed PFC approaches, stating the waiting time of passengers, the additional government subsidies, and the number of overloading situations, the maximum number of in-vehicle passengers, and the relative differences of these indicators. The relative differences (denoted as $Dev$) is calculated using the formula $\big[(\text{BDDT} - \text{PFC})/ \text{PFC} \big] \times 100 ~ (\%).$ The results show that encouraging passengers to shift their departure time leads to a reduction of 16.53\% in the waiting time, and 45.55\% decrease in the maximum number of in-vehicle passengers. These findings emphasize the importance of implementing both the trip-shifting policy and the passenger flow control strategy at the same time, as opposed to only implementing the latter.
\begin{table}
    \centering
\caption{Performance comparison between the BDDT and relaxed PFC approaches.}
    \label{tab:BDDT}
    \begin{tabular}{crrr}
\Xhline{1pt}
        Approach  & Waiting time (min) & Lost revenue  & \begin{tabular}[c]{@{}r@{}} Maximum number of \\ in-vehicle passengers \end{tabular}  \\
\Xhline{0.6pt}
        Relaxed PFC  & 97,152.00 & 0 & 3,481 \\
        BDDT & 81,097.00 & 5,518.40 &2,000\\
        $Dev (\%)$ & -16.53 & 100.00  &45.55\\
\Xhline{1pt}
\end{tabular}
\end{table}

\subsubsection{Sensitivity analysis of key parameters.}
\label{sec:sensitivity}

In the context of multi-objective optimization, the weights assigned to each term in the objective function significantly influence the resulting solutions. In our problem, the objective function \eqref{eq:objective} includes passengers’ waiting time ($F^t$), which reflects the interests of metro operators, and government-provided subsidies ($F^s$). As previously discussed, metro operators determine optimal operational strategies under different levels of government subsidy, thereby offering the government a set of feasible and efficient solutions and providing advice. To support this decision-making process, it is important for metro operators to solve the problem under various settings of the weighting coefficients. The following discussion analyzes the trade-off between these two objectives by conducting a series of experiments on Instance J, using different combinations of weighting coefficients $\omega_t$ and $\omega_s$. Specifically, the weighting coefficients $(\omega_t, \omega_s)$ are set to {(20, 1), (5, 1), (2, 1), (1, 1), (1, 2), (1, 5), (1, 20)}. The maximum allowable shift is set to 20 minutes, and all experiments are solved to optimality using our TTCBD solution method. 

\begin{figure}
    \centering
    \includegraphics[width=0.6\linewidth]{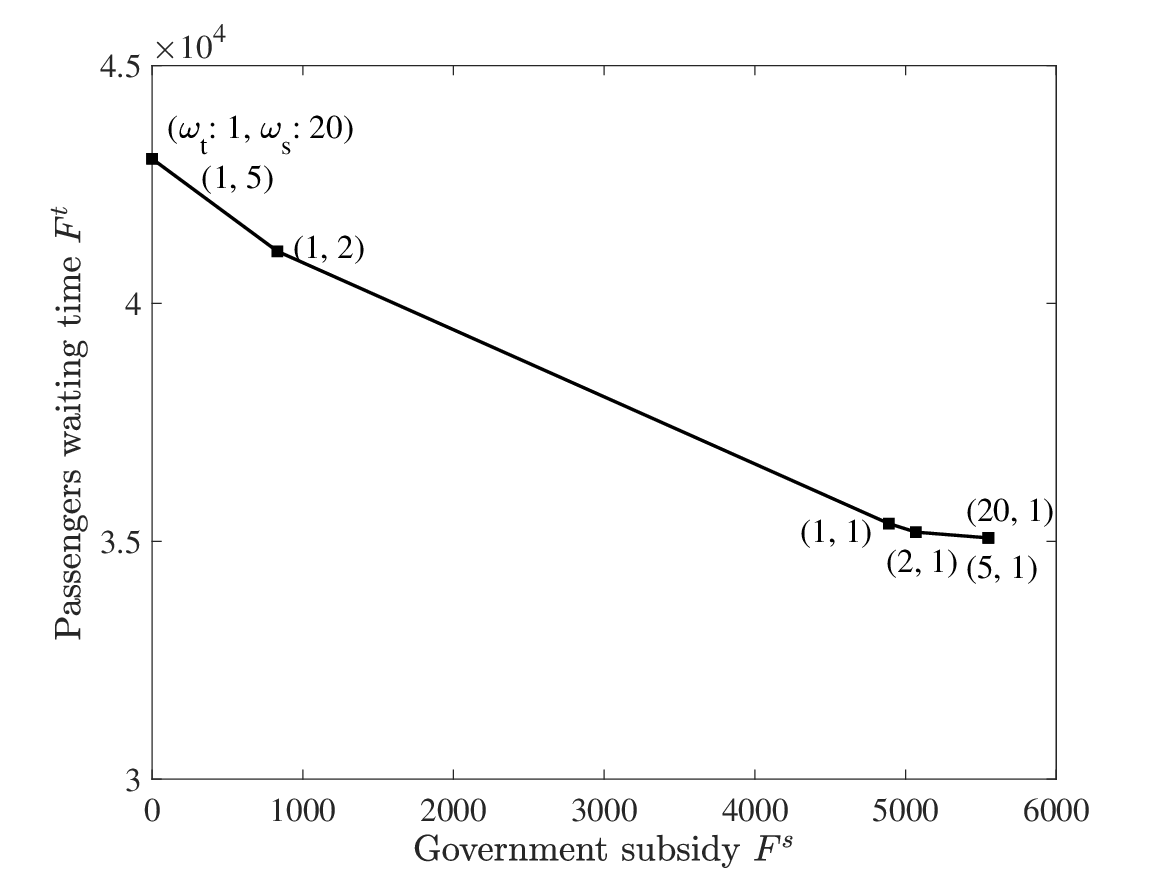}
    \caption{Pareto frontier.}
    \label{fig:Pareto}
\end{figure}

Figure \ref{fig:Pareto} shows the resulting Pareto frontier. These results provide a practical decision-support tool for metro operators when coordinating with the government. As the ratio $\omega_t/\omega_s$ increases, i.e., when the model prioritizes reducing waiting time more heavily than minimizing subsidies, the system achieves better service quality, reflected in reduced passenger waiting times. For instance, when $(\omega_t, \omega_s) = (1, 20)$, the total passenger waiting time is approximately 44,000 minutes with minimal government subsidy. However, when $\omega_s$ decreases from 20 to 2, which corresponds to an increased emphasis on passenger waiting time, the waiting time is substantially reduced to around 35,000 minutes, although this comes at the cost of a higher government subsidy. This shows how shifting model priorities can substantially improve service quality when more government support is allowed. 

A second observation is that the solutions for $(\omega_t, \omega_s) = (1, 20)$ and $(1, 5)$ are identical. This observation suggests that emphasizing government subsidy (via higher $\omega_s$) beyond a certain threshold does not further improve service quality, indicating a saturation point. This insight implies that prioritizing subsidies is a viable strategy for metro operators, as it can enable them to reach optimal service quality with appropriate support.

\section{Conclusions and future research}
\label{sec:conclusion}

In this section, we summarize the main findings of our study and outline promising directions for future research. 

\subsection{Conclusions}
In this paper, we studied the integrated optimization of booking, directing and timetabling on oversaturated urban rail transit lines. Specifically, the main focus of this study consists of determining effective passenger directing approaches, which encompass passenger flow control and trip shifting strategies, as well as to schedule train timetables. To this end, we developed an integrated INLP model to minimize the passengers' waiting time and the additional government subsidies. To improve the computational efficiency, we first linearized the above model by introducing auxiliary variables and big-$M$ constraints, and derived the most appropriate values of big-$M$ to obtain an ILP with a tighter lower bound. Thereafter, we proposed a Benders-decomposition-based approach in a branch-and-cut framework, which decomposes the integrated ILP model into a timetabling problem and a passenger assignment problem. To further enhance the solution efficiency, we proposed a novel decomposition method that incorporates partial passenger information into the timetabling problem to guide the optimization direction of the passenger assignment problem. We also integrated accelerated methods in terms of mathematical approaches and specific implementations.

To verify the effectiveness of the proposed approaches, two series of numerical experiments derived from a proof-of-concept line and the Beijing metro Batong line are implemented. The first series of experiments suggest that encouraging 21.14\% of passengers to shift their arrival times by up to 10 minutes saves at least 8.33\% of the number of 
operated
trains. Compared to the minimum additional government subsidies to serve all passengers using the limited available trains, a 23.78\% improvement in operational efficiency can be reached at the expense of a 1.00\% increase in the additional government subsidies. The second experiment compares our integrated BDDT approach and six variants of solution methods. The results indicate that the variant embedded with the cut loop, tailing off, and tree search strategies (TTCBD) performs the best in terms of execution time, especially for medium-scale and large-scale instances. In addition, we compare the performance of multiple solution approaches, namely, our TTCBD, Traditional BD, and Hybrid-LS on four test instances with various problem sizes. The results indicate that TTCBD surpassed the other two approaches in terms of both the optimality gap and the solution efficiency. 

\subsection{Future research}

This study focuses on line-level optimization under deterministic passenger demand. While the proposed modeling framework is general and scalable, and the designed algorithm is effective for solving line-level problems, several promising directions remain open for future research.

We position this study as a \textit{proof-of-concept} work that introduces an integrated framework for train timetabling and demand-side management in line-level public transport planning. It offers foundational methodologies and insights that serve as a basis for future research. The framework opens up multiple avenues for extending both the proposed model and the solution method to more complex and realistic settings. However, the proposed algorithm is somewhat tailored to line-level problems and cannot be directly applied to network-scale settings. To address network-level problems, future research can focus on the following directions.

First, extending the model to the network level is of significant interest but presents challenges. A full network-level formulation would require explicitly modeling the line dimension, capturing the interactions between lines, and accounting for transfer passenger flows across space and time. Such coupling introduces intricate interdependencies, especially when the arrival time of a transfer passenger on one line depends on the timetable of another line. 

Second, the proposed exact algorithm, which is highly effective for solving line-level problems, cannot be directly employed to the network-level setting. The network extension introduces tight couplings between lines, trains, passengers, and timetables, as well as with various demand-side management decisions. These couplings result in strong nonlinearity, increased dimensionality, and large-scale mixed-integer programming structures that would make our current solution approach computationally intractable. Addressing these challenges would require the development of new heuristic or decomposition-based algorithms capable of efficiently breaking down and solving the coupled subproblems.

Finally, network-level optimization requires the information of dynamic passenger flow across the entire system, which is difficult to obtain in practice. To address this, future studies may incorporate demand uncertainty and forecasting algorithms to enhance the robustness of the resulting timetables and demand-side management strategies. In particular, a ``{Scenario Predict-then-Optimize}'' framework that combines advanced demand forecasting technologies with stochastic programming could be explored to support more adaptive and resilient operational decisions.

\ACKNOWLEDGMENT{%
This work was supported by the National Natural Science Foundation of China (Nos. 72288101, 72322022).
}

\begin{APPENDICES}
\newpage

\section{A framework for practical applications of the proposed approaches}
\label{sec:Application}
 
In this section, we outline a framework for practical applications of our proposed integrated demand-side management and timetabling approach for an urban transit system.

\begin{figure}[h]
    \centering
    \includegraphics[height=18cm]{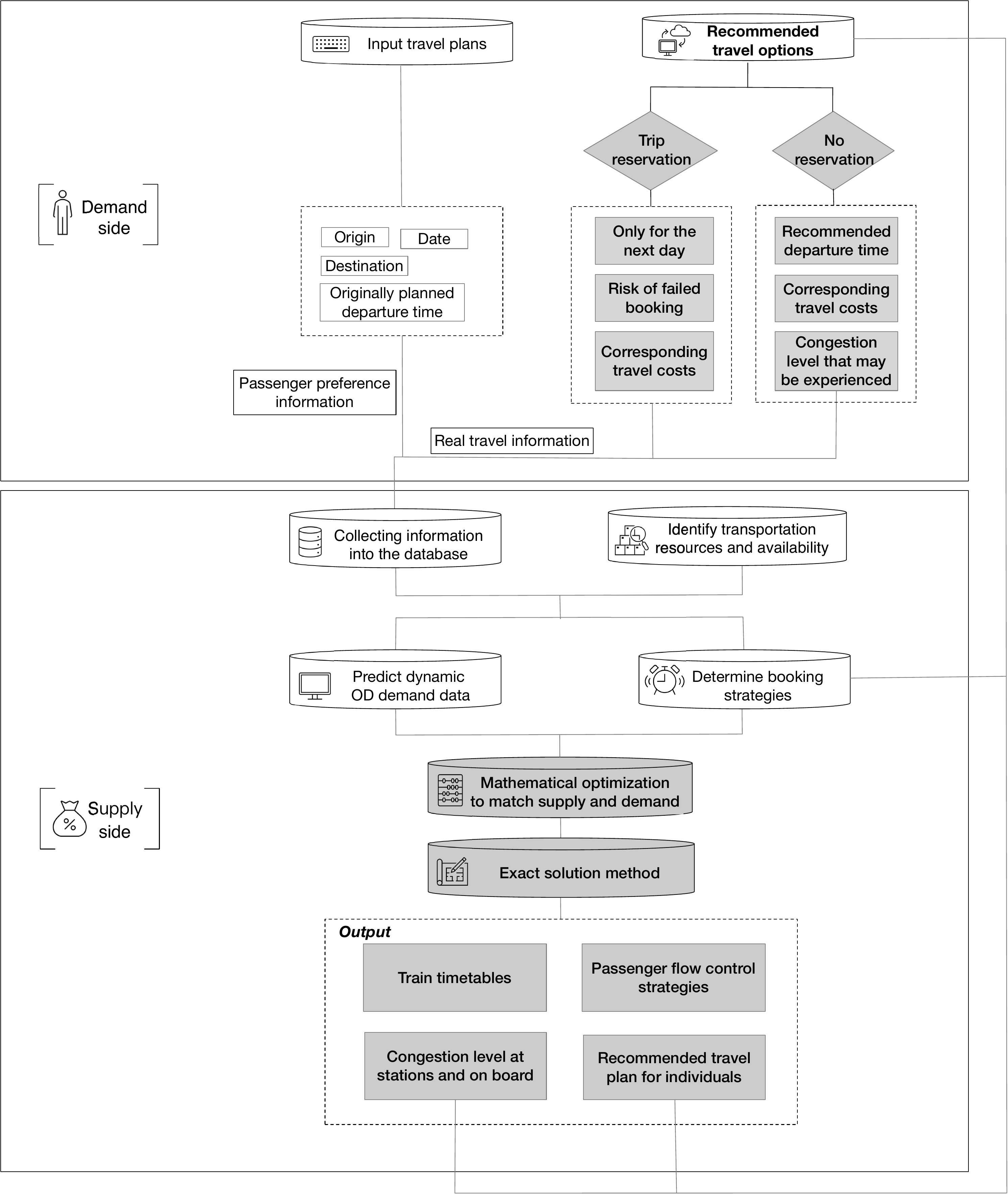}
    \caption{A framework for practical applications.}
    \label{fig:problem}
\end{figure}

\newpage

\section{Overview of the Benders decomposition algorithm for the Problem~\eqref{model:MILP}}
\label{sec:overallFramework}
In this section, we present the pseudocode of the proposed Benders decomposition algorithm for the Problem~\eqref{model:MILP} in Algorithm \ref{alg:algorithm1}.

\begin{algorithm}
\caption{The Benders decomposition algorithm for the Problem~\eqref{model:MILP}}
\label{alg:algorithm1}
\begin{algorithmic}[1]
\REQUIRE The set of stations $\mathcal{S}$, the set of trains $\mathcal{I}$, the set of timestamps $\mathcal{T}$, the time-varying OD demand $\mathbf{D}$, the reservations $\hat{\mathbf{D}}$, the tolerance $0 \leq \varepsilon_1 \leq 100$
\STATE Set the lower bound ($LB$) as $-\infty$, and set the upper bound ($UB$) as $+\infty$. Create an empty list of nodes and insert the root node into that list.
\WHILE{$(UB - LB) / LB \times 100\% > \varepsilon_1$}
    \IF{the node list is empty}
        \STATE \textbf{break}
    \ELSE
        \IF{the objective value at the current node is greater than or equal to $UB$}
            \STATE Fathom the current node and return to the beginning of the loop;
        \ELSE
            \STATE Select and branch on a non-binary variable. Remove the current node and append the resulting two branch nodes to the list of pending nodes.
            \STATE Select a pending node from the node list; Let $(\tilde{\mathbf{z}}, \tilde{\boldsymbol{\kappa}})$ be the solution of the RMP~\eqref{model:RMP} at this node, obtaining the optimal solution. Update $LB$.
            \STATE Solve the SP~\eqref{model:SP} under the solution $(\tilde{\mathbf{z}}, \tilde{\boldsymbol{\kappa}})$.
            \IF{SP is feasible}
                \STATE Obtain the optimal solution $\tilde{\mathbf{b}}$ of SP.
                \STATE Add the optimality cuts~\eqref{eq:optimalityCuts} or the strengthened optimality cuts~\eqref{eq:strengthenedoptimalityCuts} to the RMP.
                \STATE Calculate $\tilde{UB} = \omega_t F^t(\tilde{\mathbf{z}}, \tilde{\boldsymbol{\kappa}}, \tilde{\mathbf{b}}) + \omega_s F^s(\tilde{\mathbf{z}}, \tilde{\boldsymbol{\kappa}}, \tilde{\mathbf{b}})$.
                \IF{$\tilde{UB} < UB$}
                    \STATE Set $UB = \tilde{UB}$.
                    \STATE Update $(\mathbf{z}^*_0, \boldsymbol{\kappa}^*_0, \mathbf{b}^*_0)$ to be $(\tilde{\mathbf{z}}, \tilde{\boldsymbol{\kappa}}, \tilde{\mathbf{b}})$.
                \ENDIF
            \ELSE
                \STATE Add feasibility cuts~\eqref{eq:feasibilityCuts} to the RMP.
            \ENDIF
        \ENDIF
    \ENDIF
\ENDWHILE
\STATE Use $\mathbf{z}^*_0$ to solve Problem~\eqref{model:MILP} using GUROBI, obtaining the optimal integer solution $(\mathbf{z}^*, \boldsymbol{\kappa}^*, \mathbf{b}^*)$.
\RETURN $(\mathbf{z}^*, \boldsymbol{\kappa}^*, \mathbf{b}^*)$
\end{algorithmic}
\end{algorithm}

\newpage

\section{The traditional Benders decomposition framework}
\label{sec:traditional}

In this section, we introduce the traditional Benders decomposition framework commonly used in the related literature \citep[e.g.,][]{DI2022}, where the timetabling problem is the RMP and passenger assignments are all determined in SP. The RMP and the SP can now be formulated as
\begin{mini!}|s|[2]<b>
{\mathbf{z}, \theta}
{\theta} 
{} 
{}
\addConstraint
{\theta}
{ \geq \Omega(\mathbf{z}^{*}_{l})+\boldsymbol{\xi}^{T}_l(\mathbf{z}-\mathbf{z}^{*}_{l})}
{\quad \forall l \in \{1, 2, 3, ..., c_1\},}
\addConstraint
{0 }
{\geq \Psi(\mathbf{z}^{*}_{l})+\boldsymbol{\xi}^{T}_l(\mathbf{z}-\mathbf{z}^{*}_{l})}
{\quad \forall l \in \{1, 2, 3, ..., c_2\},}
\addConstraint
{}
{ \eqref{eq:0-1timetableConstraint} - \eqref{eq:0-1headwayConstraint},}
{}
\addConstraint
{\mathbf{z}}
{\in [0, 1]^{\left|\mathcal{I}\right| \times \left|\mathcal{T}\right|},}
{}
\addConstraint
{\theta}
{\geq 0,}
{}
\end{mini!}
where $c_1$ and $c_2$ indicate the number of added optimality and feasibility cuts, respectively. $\Omega(\mathbf{z}^*)$ indicates the objective function of SP under the solution of the RMP, i.e., $z^{*}$. The auxiliary decision variable $\theta$ approximates the objective function of SP. For SP, the solution $\mathbf{z}^*$ is fixed as the solution generated by the RMP. We can now develop the SP as follows:
\begin{mini!}|s|[2]<b>
{\boldsymbol{\kappa}, \mathbf{b}}
{ \omega_t F^{t}+ \omega_s F^{s} }
{} 
{ \quad \Omega(\mathbf{z}^*) =}
\addConstraint
{\mathbf{z}}
{= \mathbf{z}^{*},}
{}
\addConstraint
{\boldsymbol{\kappa}, \mathbf{b}}
{\geq 0, }
{}
\addConstraint
{}
{\eqref{eq:objFare} - \eqref{eq:waiting1}, \eqref{eq:nobookingshiftConstraint} - \eqref{eq:bookingBoardingConstraint}, \eqref{eq:nobookingBoardingLimitConstraint1} - \eqref{eq:shiftdominConstraint}, \eqref{eq:objWaitingLinearTerms}-\eqref{eq:waitingConstraintLinear}.}
{}
\end{mini!}

\newpage

\section{The pseudo code of the hybrid algorithm}
\label{sec:local}

In this section, the pseudo code of the hybrid algorithm is presented in Algorithm \ref{algorithm:hybrid}.

\begin{algorithm}[h]
\caption{The hybrid algorithm combining local search and GUROBI}
\label{algorithm:hybrid}
\begin{algorithmic}[1]
\REQUIRE $data$ (\# of trains, \# of timestamps, maximum and minimum headway limitations)
\REQUIRE $timeLimit$ (s)
\REQUIRE $initTemp,\,\alpha$ (cooling rate)
\ENSURE  $bestSol$ (vector of departure times)

\STATE Build GUROBI model
\STATE Obtain initial feasible schedule $\gets$ bestSol
\STATE curSol $\gets$ bestSol, curObj $\gets$ Obj(curSol), bestObj $\gets$ curObj
\STATE temp $\gets$ initTemp, start $\gets$ now()

\WHILE{now() -- start $<$ timeLimit}
  \STATE Randomly choose train $i$, set newTime = departure$_i\pm1$
  \IF{headway bounds violated}
    \STATE \textbf{continue}
  \ENDIF
  \STATE Fix all $z_{i,t}$ in the passenger-related subproblem according to candidate
  \STATE Solve passenger-related subproblem $\to$ candObj; if infeasible, \textbf{continue}
  \STATE $\Delta\gets$ candObj – curObj
  \IF{$\Delta<0$ \OR rand() $<\exp(-\Delta/\text{temp})$}
    \STATE curSol $\gets$ candidate, curObj $\gets$ candObj
    \IF{candObj $<$ bestObj}
      \STATE bestSol $\gets$ candidate, bestObj $\gets$ candObj
    \ENDIF
  \ENDIF
  \STATE temp $\gets\alpha\times \text{temp}$
\ENDWHILE

\RETURN bestSol
\end{algorithmic}
\end{algorithm}

\newpage
\section{Formulation for the integrated demand-side management and timetabling for urban transit networks}
\label{sec:network}

We define the set $\mathcal{L}$ as the collection of all lines in the urban transit network, the set $\mathcal{S}_l$ as the set of all stations on line $l$, and the set $\mathcal{I}_l$ as the set of all trains on line $l$. Let the set $\mathcal{H}_{u}$ denote all destination stations of passengers arriving at the origin station $u \in \mathcal{S}_l$ on line $l \in \mathcal{L}$. All variables in Table \ref{table-para} and Table \ref{table:Decision} add the line dimension to scale up to the network level. 

$\bullet$ \textbf{Timetabling-related constraints.}
\begin{align}
    z_{li(t+1)} & \leq z_{lit} &  \forall l \in \mathcal{L}, i \in \mathcal{I}_l, t \in \mathcal{T} \backslash \{ \left|\mathcal{T} \right| \}, \\
    z_{li\left|\mathcal{T}\right|}  & = 0 & \forall l \in \mathcal{L}, i \in \mathcal{I}_l, \\
    d_{li} & = \sum\limits_{t \in \mathcal{T} \backslash \{ 1 \}} [t (z_{li(t-1)} - z_{lit})] &\forall l \in \mathcal{L}, i \in \mathcal{I}_l, \\
    h_{li}  &= \begin{cases}
d_{li} - \sigma & \text{if} \ i=1\\
d_{li} - d_{l(i-1)} & \text{if} \ i \in \mathcal{I}_l \backslash \{ 1 \}\\
\end{cases} & \forall l \in \mathcal{L}, \\
h^{min} & \leq h_{li} \leq h^{max}  & \forall l \in \mathcal{L}, i \in \mathcal{I}_l \backslash \{ 1 \},  \\
x_{lit} & = \begin{cases}
z_{lit} & \text{if} \ i=1\\
z_{lit}-z_{l(i-1)t} & \text{if} \ i \in \mathcal{I}_l \backslash \{1\}\\
\end{cases} & \forall l \in \mathcal{L}, i \in \mathcal{I}_l, t \in \mathcal{T}.
\end{align} 

$\bullet$ \textbf{Constraints related to trip shifting.}

\begin{align}
    \sum\limits_{\max\{0, t-\imath\} \leq t' \leq t }\kappa_{luvt't} & = D_{luvt} & \forall l \in \mathcal{L}, u \in \mathcal{S}_l, v \in \mathcal{H}_{u}, t \in \hat{\mathcal{T}},\\
    \kappa_{luvt't} & =
\begin{cases}
 D_{luvt}& \text{if} \ t'=t\\
0 & \text{otherwise} \\
\end{cases} & \forall l \in \mathcal{L}, u \in \mathcal{S}_l,  v \in \mathcal{H}_{u}, t \in \mathcal{T} \backslash \hat{\mathcal{T}}.
\end{align}

$\bullet$ \textbf{Constraints related to serving passengers with reservations.}

\begin{align}
    \sum\limits_{i \in \mathcal{I}} \hat{b}_{liuv} &= \sum\limits_{t \in \mathcal{T}} \hat{D}_{uvt} &\forall l \in \mathcal{L}, u \in \mathcal{S}_l,  v \in \mathcal{H}_{u}.
\end{align}

$\bullet$ \textbf{Constraints related to serving all passengers without reservations.}

\begin{align}
\sum\limits_{i \in \mathcal{I}} b_{liuv} &  = \sum\limits_{t \in \mathcal{T}} D_{luvt} &\forall l \in \mathcal{L}, u \in \mathcal{S}_l, v \in \mathcal{H}_{u}.
\end{align} 

$\bullet$ \textbf{Dynamics of passengers with reservations.}

We define variable $\hat{e}_{liu}$ to denote the number of passengers with reservations who alight train $i$ at station $u$ on line $l$. The set $\mathcal{M}_u$ represents all origin stations of passengers whose destination is station $u$.
\begin{align}
\hat{o}_{liu} & =
\begin{cases}
\sum\limits_{ v \in \mathcal{H}_{u}}\hat{b}_{liuv} & \text{if} \ u=1\\
\hat{o}_{li(u-1)}-\hat{e}_{liu}+\sum\limits_{ v \in \mathcal{H}_{u}}\hat{b}_{liuv} & \text{if} \ u \in \mathcal{S}_l \backslash \{1, \left|\mathcal{S}_l\right|\}\\
0 & \text{if} \ u =\left|\mathcal{S}_l\right|\\
\end{cases}  & \forall l \in \mathcal{L}, i \in \mathcal{I}_l,\\
\hat{e}_{liu} & =
\begin{cases}
0 & \text{if} \ u=1\\
\sum\limits_{ m \in \mathcal{M}_{u},} \hat{b}_{limu} & \text{if} \ u \in \mathcal{S}_l \backslash \{1\}\\
\end{cases}   & \forall l \in \mathcal{L}, i \in \mathcal{I}_l.
\end{align}

 $\bullet$ \textbf{Dynamics of passengers without reservations.}

We now let $e_{liu}$ to denote the number of passengers without reservations who alight train $i$ at station $u$ on line $l$.
\begin{align}
o_{liu} & =
\begin{cases}
\sum\limits_{v \in \mathcal{S}_l, v \geq u}b_{liuv} & \text{if} \ u=1\\
o_{li(u-1)}-e_{liu}+\sum\limits_{v \in \mathcal{S}_l, v \geq u}b_{liuv} & \text{if} \ u \in \mathcal{S}_l \backslash \{1, \left|\mathcal{S}_l\right|\}\\
0 & \text{if} \ u =\left|\mathcal{S}_l\right|\\
\end{cases}  & \forall l \in \mathcal{L}, i \in \mathcal{I}_l,\\
e_{liu} & =
\begin{cases}
0 & \text{if} \ u=1\\
\sum\limits_{m \in \mathbf{M}_{u}} b_{limu} & \text{if} \ u \in \mathcal{S}_l \backslash \{1\}\\
\end{cases}   & \forall l \in \mathcal{L}, i \in \mathcal{I}_l.
\end{align}

$\bullet$ \textbf{Constraints related to waiting passengers with/without reservations and transferring.}

The number of waiting passengers at station $u$ on line $l$ comes from five sources: (1) newly arriving passengers with reservations whose origin is station $u$; (2) newly arriving passengers without reservations whose origin is station $u$; (3) transferring passengers who had reservations at their origin stations; (4) transferring passengers who did not have reservations at their origin stations; and (5) passengers stranded by previous trains.

We define the parameter $\alpha_{liuv}^{l'}$ is defined to indicate the ratio of passengers who alight train $i$ at station $u$ on line $l$ and aim to transfer to station $v$ on line $l'$. The parameter $\delta_{lul'}$ is introduced to indicate whether station $u$ on line $l$ is a transfer station connecting line $l'$. The parameter $L_l$ denotes the set of all possible transfer lines for line l. $\mathcal{U}_l$ denotes the set of transfer station on line $l$. $\mathcal{V}_{lu}$ denotes the set of all od pairs that will transfer at transfer station $u$ on line $l$. The variables $\hat{w}_{liuv}$ and $w_{liuv}$ are introduced as the number of passengers at station $u$ on line $l$ waiting for train $i$, with and without reservations, respectively, where station $u$ is their origin station and station $v$ is their destination. Thereafter, we introduce the variable $g_{liu}^{l'i'}$, representing the number of passengers alighting train $i'$ on line $l'$ and transferring to train $i$ at station $u$ on line $l$. The variable $r_{liuv}$ is defined to denote the number of stranded passengers. The variable $f_{liu}$ denotes the number of transfer passengers.

\begin{align}
    \hat{w}_{liuv} &=
\begin{cases}
\sum\limits_{t \in \mathcal{T}}z_{lit}\hat{D}_{luvt} & \text{if} \ i=1\\
\sum\limits_{t \in \mathcal{T}}z_{lit}\hat{D}_{luvt} - \sum\limits_{j \in I_{l}, j \leq i-1} \hat{b}_{juv} & \text{if} \  i \in \mathcal{I}_{l} \backslash \{1\}\\
\end{cases} &\forall l \in \mathcal{L}, u \in \mathcal{S}_l, v \in \mathcal{H}_{u}.
\end{align}
\begin{align}
   w_{liuv} &=
\begin{cases}
\sum\limits_{t' \in \mathcal{T}}z_{it'}\sum\limits_{t' \leq t \leq \min\{ \left|\mathcal{T} \right|, t' +\imath \} }\kappa_{uvt't} & \text{if} \ i=1\\
\sum\limits_{t' \in \mathcal{T}}z_{it'}\sum\limits_{t' \leq t \leq \min\{ \left|\mathcal{T} \right|, t' +\imath \}  }\kappa_{uvt't} - \sum\limits_{j \in I_{l}, j \leq i-1} b_{ljuv} & \text{if} \ i \in \mathcal{I}_l \backslash \{1\}\\
\end{cases} &\forall u \in \mathcal{S}_l, v \in \mathcal{H}_{u}.\\
f_{liu} &= \sum_{l'\in L_l}\sum_{(u',v')\in V_{lu}}g^{l',i'}_{lim} (\hat{b}_{l'i'u'v'} + b_{l'i'u'v'}), &\forall l\in\mathcal{L}, u \in \mathcal{U}_l, i \in \mathcal{I}_{l}.
\end{align}

The number of stranded passengers can be expressed as 
\begin{align}
    r_{liuv} &= w_{liuv} - b_{liuv} & \forall l \in \mathcal{L}, i \in \mathcal{I}_{l}, u \in \mathcal{S}_{l}, v \in \mathcal{H}_{u}.
\end{align}

\textbf{$\bullet$ Constraints related to passenger flow control.}
\begin{align}
\varrho_{liuv} w_{liuv} & \leq b_{liuv} \leq w_{liuv}  & \forall l \in \mathcal{L}, i \in \mathcal{I}_{l}, u \in \mathcal{S}_{l}, v \in \mathcal{H}_{u}.
\end{align}

\textbf{$\bullet$ Constraints related to train capacity.}
\begin{align}
o_{liu} + \hat{o}_{liu} + f_{liu} &\leq C^{max}, &\forall  \in \mathcal{L}, i \in \mathcal{I}_{l}, u \in \mathcal{S}_{l}.
\end{align}

\textbf{$\bullet$ Constraints related to transfer.}
\begin{align}
d_{li} - (d_{l',i'} +\varphi^{l'}_{l,i}) \geq M( q_{liu}^{l'i'} - 1) && \forall l \in \mathcal{L}, i \in \mathcal{I}_{l}, u \in \mathcal{U}_{l}, l'\in L_l, i'\in\mathcal{U}_{l'}.\\
d_{li} - (d_{l',i'} +\varphi^{l'}_{l,i}) \leq M q_{liu}^{l'i'}&& \forall l \in \mathcal{L}, i \in \mathcal{I}_{l}, u \in \mathcal{U}_{l}, l'\in L_l, i'\in\mathcal{U}_{l'}.\\
  g_{liu}^{l'i'} =
\begin{cases}
q_{liu}^{l'i'}  &\text{if} \ i=1\\
q_{liu}^{l'i'} - q_{l(i-1)u}^{l'i'}  &\text{if} \ i \in \mathcal{I}_l \backslash \{1\}\\
\end{cases} &&\forall l \in \mathcal{L}, i \in \mathcal{I}_{l}, u \in \mathcal{U}_{l}, l'\in L_l.
\end{align}

\textbf{$\bullet$ Objective functions.}

\begin{align}
   \min & \quad \omega_t F^{t}+ \omega_s F^{s}  \\
   \quad F^t &=\sigma \Big[\sum_{l\in\mathcal{L}}\sum_{i\in\mathcal{I}_l}\sum_{u\in\mathcal{S}_l}\sum_{t\in\mathcal{T}}(\hat{p}^{wc}_{iut}+p^{wc}_{iut}) + \sum_{l\in\mathcal{L}}\sum\limits_{i \in \mathcal{I}_l}\sum_{u\in\mathcal{S}_l}\sum\limits_{t \in \mathcal{T}}(x_{lit}\sum\limits_{v \in \mathcal{S}_{u+1}}r_{liuv}) \\ &+\sum_{l\in\mathcal{L}}\sum_{l'\in L_l}\sum_{i\in\mathcal{I}_l}\sum_{i'\in\mathcal{I}_l'}\sum_{u\in\mathcal{U}_l} g_{liu}^{l'i'} f_{liu} (d_{li}-d_{l'i'} + \varphi^{l'}_{l,i})   \Big],  \\
  \quad F^{s} &= \sum_{l\in\mathcal{L}}\sum\limits_{u\in \mathcal{S}_l}\sum\limits_{v \in \mathcal{H}_{u}}\sum\limits_{t \in \mathcal{T}} \Big[D_{luvt} \varepsilon_{uv} - \varepsilon_{uv}[\sum\limits_{t+1\leq t' \leq \min\{\left|\mathcal{T} \right|, t +\imath\}}(\phi\kappa_{luvtt'}+\kappa_{luvtt})] \Big]. 
\end{align}

\begin{align}
    \hat{p}^w_{liut} & = x_{lit}\sum_{v\in\mathcal{H}_{u}}\hat{D}_{luvt} & \forall  l\in \mathcal{L}, i \in \mathcal{I}_{l}, u \in \mathcal{S}_{l} ,t \in \mathcal{T}, \\
    p^w_{liut}& = x_{lit}\sum_{v\in\mathcal{H}_{v}}\sum_{t\leq t^{''}\leq \min\{ \left|\mathcal{T} \right|, t +\imath \}}\kappa_{luvtt^{''}} & \forall  l\in \mathcal{L}, i \in \mathcal{I}_{l}, u \in \mathcal{S}_{l} ,t \in \mathcal{T}, \\
    \hat{p}^{wc}_{liut} & = x_{lit}\sum_{t' \in \mathcal{T}, t'\leq t}\hat{p}^w_{liut'} &  \forall  l\in \mathcal{L}, i \in \mathcal{I}_{l}, u \in \mathcal{S}_{l} ,t \in \mathcal{T}, \\
    p^{wc}_{liut} &= x_{lit}\sum_{t' \in \mathcal{T}, t'\leq t}p^w_{liut'} & \forall  l\in \mathcal{L}, i \in \mathcal{I}_{l}, u \in \mathcal{S}_{l} ,t \in \mathcal{T}.
\end{align}

\end{APPENDICES}


\begin{thebibliography}{43}
\providecommand{\natexlab}[1]{#1}
\providecommand{\url}[1]{\texttt{#1}}
\providecommand{\urlprefix}{URL }

\bibitem[{Bao et~al.(2023)Bao, Yang, Gao and Xu}]{BAO2023}
\bibinfo{author}{Bao, Y.}, \bibinfo{author}{Yang, H.}, \bibinfo{author}{Gao,
  Z.}, \bibinfo{author}{Xu, H.}, \bibinfo{year}{2023}.
\newblock \bibinfo{title}{How do pre-event activities alleviate congestion and
  increase attendees’ travel utility and the venue's profit during a special
  event?}
\newblock \bibinfo{journal}{Transportation Research Part B: Methodological}
  \bibinfo{volume}{173}, \bibinfo{pages}{332--353}.
\bibitem[{Barz and Gartner(2016)}]{Barz2016}
\bibinfo{author}{Barz, C.}, \bibinfo{author}{Gartner, D.},
  \bibinfo{year}{2016}.
\newblock \bibinfo{title}{Air cargo network revenue management}.
\newblock \bibinfo{journal}{Transportation Science} \bibinfo{volume}{50},
  \bibinfo{pages}{1206--1222}.
\bibitem[{{Beijing Municipal Commission of Transport}(2020)}]{BMCT}
\bibinfo{author}{{Beijing Municipal Commission of Transport}}, \bibinfo{year}{2020}.
\newblock \bibinfo{title}{{The station entry reservation trial will be launched at two Beijing Metro stations starting from March 6}}.
\newblock \bibinfo{howpublished}{\url{https://jtw.beijing.gov.cn/xxgk/dtxx/202003/t20200305_1679196.html}}.
\newblock \bibinfo{note}{[Accessed April 12, 2025]}.

\bibitem[{Benders(1962)}]{benders1962}
\bibinfo{author}{Benders, J.}, \bibinfo{year}{1962}.
\newblock \bibinfo{title}{Partitioning procedures for solving mixed-variables
  programming problems}.
\newblock \bibinfo{journal}{Numerische Mathematik} \bibinfo{volume}{4},
  \bibinfo{pages}{238–252}.
\bibitem[{Binder et~al.(2021)Binder, Maknoon, {Sharif Azadeh} and
  Bierlaire}]{Binder2021}
\bibinfo{author}{Binder, S.}, \bibinfo{author}{Maknoon, M.},
  \bibinfo{author}{{Sharif Azadeh}, S.}, \bibinfo{author}{Bierlaire, M.},
  \bibinfo{year}{2021}.
\newblock \bibinfo{title}{Passenger-centric timetable rescheduling: A user
  equilibrium approach}.
\newblock \bibinfo{journal}{Transportation Research Part C: Emerging
  Technologies} \bibinfo{volume}{132}, \bibinfo{pages}{103368}.
\bibitem[{{China News}(2020)}]{ChinaNews}
\bibinfo{author}{{China News}}, \bibinfo{year}{2020}.
\newblock \bibinfo{title}{{The trial of station entry reservation will be
  launched at Caofang Station of Beijing Subway Line 6 starting from April
  29th}}.
\newblock
  \bibinfo{howpublished}{\url{https://m.chinanews.com/wap/detail/chs/zw/9169668.shtml}}.
\newblock \bibinfo{note}{[Accessed April 12, 2024]}.
\bibitem[{Copeland and McKenney(1988)}]{Air1}
\bibinfo{author}{Copeland, D.G.}, \bibinfo{author}{McKenney, J.L.},
  \bibinfo{year}{1988}.
\newblock \bibinfo{title}{Airline reservations systems: Lessons from history}.
\newblock \bibinfo{journal}{MIS Quarterly} \bibinfo{volume}{12},
  \bibinfo{pages}{353--370}.
\bibitem[{Croella et~al.(2024)Croella, Luteberget, Mannino and
  Ventura}]{Carlo2024}
\bibinfo{author}{Croella, A.L.}, \bibinfo{author}{Luteberget, B.},
  \bibinfo{author}{Mannino, C.}, \bibinfo{author}{Ventura, P.},
  \bibinfo{year}{2024}.
\newblock \bibinfo{title}{A maxsat approach for solving a new dynamic
  discretization discovery model for train rescheduling problems}.
\newblock \bibinfo{journal}{Computers \& Operations Research}
  \bibinfo{volume}{167}, \bibinfo{pages}{106679}.
\bibitem[{Di et~al.(2022)Di, Yang, Shi, Zhou, Yang and Gao}]{DI2022}
\bibinfo{author}{Di, Z.}, \bibinfo{author}{Yang, L.}, \bibinfo{author}{Shi,
  J.}, \bibinfo{author}{Zhou, H.}, \bibinfo{author}{Yang, K.},
  \bibinfo{author}{Gao, Z.}, \bibinfo{year}{2022}.
\newblock \bibinfo{title}{Joint optimization of carriage arrangement and flow
  control in a metro-based underground logistics system}.
\newblock \bibinfo{journal}{Transportation Research Part B: Methodological}
  \bibinfo{volume}{159}, \bibinfo{pages}{1--23}.
\bibitem[{Ding et~al.(2023)Ding, Yang, Qin and Xu}]{YangHai2023}
\bibinfo{author}{Ding, H.}, \bibinfo{author}{Yang, H.}, \bibinfo{author}{Qin,
  X.}, \bibinfo{author}{Xu, H.}, \bibinfo{year}{2023}.
\newblock \bibinfo{title}{Credit charge-cum-reward scheme for green multi-modal
  mobility}.
\newblock \bibinfo{journal}{Transportation Research Part B: Methodological}
  \bibinfo{volume}{178}, \bibinfo{pages}{102852}.
\bibitem[{Fischetti et~al.(2016)Fischetti, Ljubi\'{c} and
  Sinnl}]{Fischetti2016}
\bibinfo{author}{Fischetti, M.}, \bibinfo{author}{Ljubi\'{c}, I.},
  \bibinfo{author}{Sinnl, M.}, \bibinfo{year}{2016}.
\newblock \bibinfo{title}{Benders decomposition without separability: A
  computational study for capacitated facility location problems}.
\newblock \bibinfo{journal}{European Journal of Operational Research}
  \bibinfo{volume}{253}, \bibinfo{pages}{557--569}.
\bibitem[{Fischetti et~al.(2017)Fischetti, Ljubi\'{c} and Sinnl}]{FischettiMS}
\bibinfo{author}{Fischetti, M.}, \bibinfo{author}{Ljubi\'{c}, I.},
  \bibinfo{author}{Sinnl, M.}, \bibinfo{year}{2017}.
\newblock \bibinfo{title}{Redesigning Benders decomposition for large-scale
  facility location}.
\newblock \bibinfo{journal}{Management Science} \bibinfo{volume}{63},
  \bibinfo{pages}{2146--2162}.
\bibitem[{He et~al.(2017)He, Chen, Xiong, Zhu and Zhang}]{He2017}
\bibinfo{author}{He, X.}, \bibinfo{author}{Chen, X.M.}, \bibinfo{author}{Xiong,
  C.}, \bibinfo{author}{Zhu, Z.}, \bibinfo{author}{Zhang, L.},
  \bibinfo{year}{2017}.
\newblock \bibinfo{title}{Optimal time-varying pricing for toll roads under
  multiple objectives: A simulation-based optimization approach}.
\newblock \bibinfo{journal}{Transportation Science} \bibinfo{volume}{51},
  \bibinfo{pages}{412--426}.
\bibitem[{Hu et~al.(2023)Hu, Li, Wang, Zhang, Wei and Yang}]{HU2023}
\bibinfo{author}{Hu, Y.}, \bibinfo{author}{Li, S.}, \bibinfo{author}{Wang, Y.},
  \bibinfo{author}{Zhang, H.}, \bibinfo{author}{Wei, Y.},
  \bibinfo{author}{Yang, L.}, \bibinfo{year}{2023}.
\newblock \bibinfo{title}{Robust metro train scheduling integrated with
  skip-stop pattern and passenger flow control strategy under uncertain
  passenger demands}.
\newblock \bibinfo{journal}{Computers \& Operations Research}
  \bibinfo{volume}{151}, \bibinfo{pages}{106116}.
\bibitem[{Leutwiler and Corman(2022)}]{Leutwiler2022}
\bibinfo{author}{Leutwiler, F.}, \bibinfo{author}{Corman, F.},
  \bibinfo{year}{2022}.
\newblock \bibinfo{title}{A logic-based benders decomposition for microscopic
  railway timetable planning}.
\newblock \bibinfo{journal}{European Journal of Operational Research}
  \bibinfo{volume}{303}, \bibinfo{pages}{525--540}.
\bibitem[{Leutwiler and Corman(2023)}]{Leutwiler2023}
\bibinfo{author}{Leutwiler, F.}, \bibinfo{author}{Corman, F.},
  \bibinfo{year}{2023}.
\newblock \bibinfo{title}{Set covering heuristics in a benders decomposition
  for railway timetabling}.
\newblock \bibinfo{journal}{Computers \& Operations Research}
  \bibinfo{volume}{159}, \bibinfo{pages}{106339}.
\bibitem[{Li et~al.(2017)Li, Dessouky, Yang and Gao}]{LI2017}
\bibinfo{author}{Li, S.}, \bibinfo{author}{Dessouky, M.M.},
  \bibinfo{author}{Yang, L.}, \bibinfo{author}{Gao, Z.}, \bibinfo{year}{2017}.
\newblock \bibinfo{title}{Joint optimal train regulation and passenger flow
  control strategy for high-frequency metro lines}.
\newblock \bibinfo{journal}{Transportation Research Part B: Methodological}
  \bibinfo{volume}{99}, \bibinfo{pages}{113--137}.
\bibitem[{Li et~al.(2023)Li, Yang and Ke}]{HaiYang2023}
\bibinfo{author}{Li, X.}, \bibinfo{author}{Yang, H.}, \bibinfo{author}{Ke, J.},
  \bibinfo{year}{2023}.
\newblock \bibinfo{title}{Booking cum rationing strategy for equitable travel
  demand management in road networks}.
\newblock \bibinfo{journal}{Transportation Research Part B: Methodological}
  \bibinfo{volume}{167}, \bibinfo{pages}{261--274}.
\bibitem[{Liang et~al.(2023)Liang, Lyu, Teo and Gao}]{Liang2023}
\bibinfo{author}{Liang, J.}, \bibinfo{author}{Lyu, G.}, \bibinfo{author}{Teo,
  C.P.}, \bibinfo{author}{Gao, Z.}, \bibinfo{year}{2023}.
\newblock \bibinfo{title}{Online passenger flow control in metro lines}.
\newblock \bibinfo{journal}{Operations Research} \bibinfo{volume}{71},
  \bibinfo{pages}{768--775}.
\bibitem[{Liu et~al.(2020)Liu, Li and Yang}]{RMLiu2020}
\bibinfo{author}{Liu, R.}, \bibinfo{author}{Li, S.}, \bibinfo{author}{Yang,
  L.}, \bibinfo{year}{2020}.
\newblock \bibinfo{title}{Collaborative optimization for metro train scheduling
  and train connections combined with passenger flow control strategy}.
\newblock \bibinfo{journal}{Omega} \bibinfo{volume}{90},
  \bibinfo{pages}{101990}.
\bibitem[{Liu et~al.(2015)Liu, Yang and Yin}]{HaiYang2015}
\bibinfo{author}{Liu, W.}, \bibinfo{author}{Yang, H.}, \bibinfo{author}{Yin,
  Y.}, \bibinfo{year}{2015}.
\newblock \bibinfo{title}{Efficiency of a highway use reservation system for
  morning commute}.
\newblock \bibinfo{journal}{Transportation Research Part C: Emerging
  Technologies} \bibinfo{volume}{56}, \bibinfo{pages}{293--308}.
\bibitem[{Lu et~al.(2023)Lu, Yang, Yang, Zhou and Gao}]{Lu2023}
\bibinfo{author}{Lu, Y.}, \bibinfo{author}{Yang, L.}, \bibinfo{author}{Yang,
  H.}, \bibinfo{author}{Zhou, H.}, \bibinfo{author}{Gao, Z.},
  \bibinfo{year}{2023}.
\newblock \bibinfo{title}{Robust collaborative passenger flow control on a
  congested metro line: A joint optimization with train timetabling}.
\newblock \bibinfo{journal}{Transportation Research Part B: Methodological}
  \bibinfo{volume}{168}, \bibinfo{pages}{27--55}.
\bibitem[{Lu et~al.(2022)Lu, Yang, Yang, Gao, Zhou, Meng and Qi}]{Lu2022}
\bibinfo{author}{Lu, Y.}, \bibinfo{author}{Yang, L.}, \bibinfo{author}{Yang,
  K.}, \bibinfo{author}{Gao, Z.}, \bibinfo{author}{Zhou, H.},
  \bibinfo{author}{Meng, F.}, \bibinfo{author}{Qi, J.}, \bibinfo{year}{2022}.
\newblock \bibinfo{title}{A distributionally robust optimization method for
  passenger flow control strategy and train scheduling on an urban rail transit
  line}.
\newblock \bibinfo{journal}{Engineering} \bibinfo{volume}{12},
  \bibinfo{pages}{202--220}.
\bibitem[{Meng et~al.(2022)Meng, Yang, Shi, Jiang and Gao}]{Meng2022}
\bibinfo{author}{Meng, F.}, \bibinfo{author}{Yang, L.}, \bibinfo{author}{Shi,
  J.}, \bibinfo{author}{Jiang, Z.}, \bibinfo{author}{Gao, Z.},
  \bibinfo{year}{2022}.
\newblock \bibinfo{title}{Collaborative passenger flow control for
  oversaturated metro lines: A stochastic optimization method}.
\newblock \bibinfo{journal}{Transportmetrica A: Transport Science}
  \bibinfo{volume}{18}, \bibinfo{pages}{619--658}.
\bibitem[{{Ministry of Transport of the People's Republic of
  China}(2022)}]{china}
\bibinfo{author}{{Ministry of Transport of the People's Republic of China}},
  \bibinfo{year}{2022}.
\newblock \bibinfo{title}{{The transportation authority optimizes the capacity
  to serve the large passenger flow, so that the city artery keeps surging
  up}}.
\newblock
  \bibinfo{howpublished}{\url{https://www.mot.gov.cn/jiaotongyaowen/202212/t20221228_3730584.html/}}.
\newblock \bibinfo{note}{[Accessed March 29, 2023]}.
\bibitem[{{Motoring}(2024)}]{Singapore}
\bibinfo{author}{{Motoring}}, \bibinfo{year}{2024}.
\newblock \bibinfo{title}{{Electronic Road Pricing}}.
\newblock
  \bibinfo{howpublished}{\url{https://onemotoring.lta.gov.sg/content/onemotoring/home/driving/ERP.html}}.
\newblock \bibinfo{note}{[Accessed April 12, 2024]}.
\bibitem[{{Nederlandse Spoorwegen}(2024)}]{NS}
\bibinfo{author}{{Nederlandse Spoorwegen}}, \bibinfo{year}{2024}.
\newblock \bibinfo{title}{{When can you travel with a discount?}}
\newblock
  \bibinfo{howpublished}{\url{https://www.ns.nl/en/featured/traveling-with-discount/when-can-you-travel-with-a-discount.html}}.
\newblock \bibinfo{note}{[Accessed April 12, 2024]}.
\bibitem[{{People}(2020)}]{BeijingNews}
\bibinfo{author}{{People}}, \bibinfo{year}{2020}.
\newblock \bibinfo{title}{{Reservations can reduce waiting in line due to
  passenger flow control}}.
\newblock
  \bibinfo{howpublished}{\url{http://bj.people.com.cn/n2/2020/0430/c82840-33987970.html}}.
\newblock \bibinfo{note}{[Accessed September 12, 2024]}.
\bibitem[{Polinder et~al.(2022)Polinder, Cacchiani, Schmidt and
  Huisman}]{Polinder2022}
\bibinfo{author}{Polinder, G.J.}, \bibinfo{author}{Cacchiani, V.},
  \bibinfo{author}{Schmidt, M.}, \bibinfo{author}{Huisman, D.},
  \bibinfo{year}{2022}.
\newblock \bibinfo{title}{An iterative heuristic for passenger-centric train
  timetabling with integrated adaption times}.
\newblock \bibinfo{journal}{Computers \& Operations Research}
  \bibinfo{volume}{142}, \bibinfo{pages}{105740}.
\bibitem[{Rahmaniani et~al.(2020)Rahmaniani, Ahmed, Crainic, Gendreau and
  Rei}]{Rahmaniani2020}
\bibinfo{author}{Rahmaniani, R.}, \bibinfo{author}{Ahmed, S.},
  \bibinfo{author}{Crainic, T.G.}, \bibinfo{author}{Gendreau, M.},
  \bibinfo{author}{Rei, W.}, \bibinfo{year}{2020}.
\newblock \bibinfo{title}{The Benders dual decomposition method}.
\newblock \bibinfo{journal}{Operations Research} \bibinfo{volume}{68},
  \bibinfo{pages}{878--895}.
\bibitem[{Robenek et~al.(2018)Robenek, {Sharif Azadeh}, Maknoon, {De Lapparent}
  and Bierlaire}]{shadi2018}
\bibinfo{author}{Robenek, T.}, \bibinfo{author}{{Sharif Azadeh}, S.},
  \bibinfo{author}{Maknoon, Y.}, \bibinfo{author}{{De Lapparent}, M.},
  \bibinfo{author}{Bierlaire, M.}, \bibinfo{year}{2018}.
\newblock \bibinfo{title}{Train timetable design under elastic passenger
  demand}.
\newblock \bibinfo{journal}{Transportation Research Part B: Methodological}
  \bibinfo{volume}{111}, \bibinfo{pages}{19--38}.
\bibitem[{Rothstein(1985)}]{Air2}
\bibinfo{author}{Rothstein, M.}, \bibinfo{year}{1985}.
\newblock \bibinfo{title}{OR and the airline overbooking problem}.
\newblock \bibinfo{journal}{Operations Research} \bibinfo{volume}{33},
  \bibinfo{pages}{237--248}.
\bibitem[{Shi et~al.(2018)Shi, Yang, Yang and Gao}]{Shi2018}
\bibinfo{author}{Shi, J.}, \bibinfo{author}{Yang, L.}, \bibinfo{author}{Yang,
  J.}, \bibinfo{author}{Gao, Z.}, \bibinfo{year}{2018}.
\newblock \bibinfo{title}{Service-oriented train timetabling with collaborative
  passenger flow control on an oversaturated metro line: An integer linear
  optimization approach}.
\newblock \bibinfo{journal}{Transportation Research Part B: Methodological}
  \bibinfo{volume}{110}, \bibinfo{pages}{26--59}.
\bibitem[{{Transport for London}(2024)}] {London}
\bibinfo{author}{{Transport for London}}, \bibinfo{year}{2024}.
\newblock \bibinfo{title}{{Congestion Charge zone}}.
\newblock
  \bibinfo{howpublished}{\url{https://tfl.gov.uk/modes/driving/congestion-charge/congestion-charge-zone}}.
\newblock \bibinfo{note}{[Accessed April 12, 2024]}.
\bibitem[{{United Nations}(2019)}]{un75}
\bibinfo{author}{{United Nations}}, \bibinfo{year}{2019}.
\newblock \bibinfo{title}{{Shifting Demographics: A Visual Guide}}.
\newblock
  \bibinfo{howpublished}{\url{https://www.un.org/en/un75/shifting-demographics}}.
\newblock \bibinfo{note}{[Accessed March 28, 2023]}.

\bibitem[{Wang et~al.(2023)Wang, Jin, Sibul and Wei}]{Jin2023}
\bibinfo{author}{Wang, L.}, \bibinfo{author}{Jin, J.G.},
  \bibinfo{author}{Sibul, G.}, \bibinfo{author}{Wei, Y.}, \bibinfo{year}{2023}.
\newblock \bibinfo{title}{Designing metro network expansion: Deterministic and
  robust optimization models}.
\newblock \bibinfo{journal}{Networks and Spatial Economics}
  \bibinfo{volume}{23}, \bibinfo{pages}{317--347}.
\bibitem[{Xia et~al.(2024)Xia, Ma and {Sharif Azadeh}}]{XIA2024}
\bibinfo{author}{Xia, D.}, \bibinfo{author}{Ma, J.}, \bibinfo{author}{{Sharif
  Azadeh}, S.}, \bibinfo{year}{2024}.
\newblock \bibinfo{title}{Integrated timetabling and vehicle scheduling of an
  intermodal urban transit network: A distributionally robust optimization
  approach}.
\newblock \bibinfo{journal}{Transportation Research Part C: Emerging
  Technologies} \bibinfo{volume}{162}, \bibinfo{pages}{104610}.
\bibitem[{Xia et~al.(2023)Xia, Ma, {Sharif Azadeh} and Zhang}]{Xia2023}
\bibinfo{author}{Xia, D.}, \bibinfo{author}{Ma, J.}, \bibinfo{author}{{Sharif
  Azadeh}, S.}, \bibinfo{author}{Zhang, W.}, \bibinfo{year}{2023}.
\newblock \bibinfo{title}{Data-driven distributionally robust timetabling and
  dynamic-capacity allocation for automated bus systems with modular vehicles}.
\newblock \bibinfo{journal}{Transportation Research Part C: Emerging
  Technologies} \bibinfo{volume}{155}, \bibinfo{pages}{104314}.
\bibitem[{Xia et~al.(2024)Xia, Ma, {Sharif Azadeh} and Zhang}]{bib0099}
\bibinfo{author}{Xia, D.}, \bibinfo{author}{Ma, J.}, \bibinfo{author}{{Sharif
  Azadeh}, S.},
\bibinfo{year}{2024}.
\newblock \bibinfo{title}{Integrated timetabling, vehicle scheduling, and dynamic capacity allocation of modular autonomous vehicles under demand uncertainty}.
\newblock \bibinfo{volume}{arXiv:2410.16409}.
\newblock \bibinfo{howpublished}{\url{https://arxiv.org/abs/2410.16409}}.
\bibitem[{Xiao et~al.(2015)Xiao, Huang and Liu}]{Xiao2015}
\bibinfo{author}{Xiao, L.L.}, \bibinfo{author}{Huang, H.J.},
  \bibinfo{author}{Liu, R.}, \bibinfo{year}{2015}.
\newblock \bibinfo{title}{Congestion behavior and tolls in a bottleneck model
  with stochastic capacity}.
\newblock \bibinfo{journal}{Transportation Science} \bibinfo{volume}{49},
  \bibinfo{pages}{46--65}.
\bibitem[{Yang and Bell(1998)}]{HaiYang1998}
\bibinfo{author}{Yang, H.}, \bibinfo{author}{Bell, M.G.H.},
  \bibinfo{year}{1998}.
\newblock \bibinfo{title}{Models and algorithms for road network design: A
  review and some new developments}.
\newblock \bibinfo{journal}{Transport Reviews} \bibinfo{volume}{18},
  \bibinfo{pages}{257--278}.
\bibitem[{Yang et~al.(2020)Yang, Shao, Wang and Ye}]{HaiYang2020}
\bibinfo{author}{Yang, H.}, \bibinfo{author}{Shao, C.}, \bibinfo{author}{Wang,
  H.}, \bibinfo{author}{Ye, J.}, \bibinfo{year}{2020}.
\newblock \bibinfo{title}{Integrated reward scheme and surge pricing in a
  ridesourcing market}.
\newblock \bibinfo{journal}{Transportation Research Part B: Methodological}
  \bibinfo{volume}{134}, \bibinfo{pages}{126--142}.
\bibitem[{Yang and Wang(2011)}]{YangHai2011}
\bibinfo{author}{Yang, H.}, \bibinfo{author}{Wang, X.}, \bibinfo{year}{2011}.
\newblock \bibinfo{title}{Managing network mobility with tradable credits}.
\newblock \bibinfo{journal}{Transportation Research Part B: Methodological}
  \bibinfo{volume}{45}, \bibinfo{pages}{580--594}.
\bibitem[{Yin et~al.(2021)Yin, D’Ariano, Wang, Yang and Tang}]{YIN2021183}
\bibinfo{author}{Yin, J.}, \bibinfo{author}{D’Ariano, A.},
  \bibinfo{author}{Wang, Y.}, \bibinfo{author}{Yang, L.},
  \bibinfo{author}{Tang, T.}, \bibinfo{year}{2021}.
\newblock \bibinfo{title}{Timetable coordination in a rail transit network with
  time-dependent passenger demand}.
\newblock \bibinfo{journal}{European Journal of Operational Research}
  \bibinfo{volume}{295}, \bibinfo{pages}{183--202}.
\bibitem[{Yin et~al.(2023)Yin, Pu, Yang, D’Ariano and Wang}]{YIN2023}
\bibinfo{author}{Yin, J.}, \bibinfo{author}{Pu, F.}, \bibinfo{author}{Yang,
  L.}, \bibinfo{author}{D’Ariano, A.}, \bibinfo{author}{Wang, Z.},
  \bibinfo{year}{2023}.
\newblock \bibinfo{title}{Integrated optimization of rolling stock allocation
  and train timetables for urban rail transit networks: A Benders decomposition
  approach}.
\newblock \bibinfo{journal}{Transportation Research Part B: Methodological}
  \bibinfo{volume}{176}, \bibinfo{pages}{102815}.
\bibitem[{Yuan et~al.(2022)Yuan, Li, Liu, Yang and Gao}]{Yuan2022}
\bibinfo{author}{Yuan, Y.}, \bibinfo{author}{Li, S.}, \bibinfo{author}{Yang, L.}, \bibinfo{author}{Gao, Z.},
  \bibinfo{year}{2022}.
\newblock \bibinfo{title}{Real-time optimization of train regulation and passenger flow control for urban rail transit network under frequent disturbances}.
\newblock \bibinfo{journal}{Transportation Research Part E: Logistics and Transportation Review} \bibinfo{volume}{168}, \bibinfo{pages}{102942}.
\bibitem[{Yuan et~al.(2023)Yuan, Li, Liu, Yang and Gao}]{YuanPartC}
\bibinfo{author}{Yuan, Y.}, \bibinfo{author}{Li, S.}, \bibinfo{author}{Liu,
  R.}, \bibinfo{author}{Yang, L.}, \bibinfo{author}{Gao, Z.},
  \bibinfo{year}{2023}.
\newblock \bibinfo{title}{Decomposition and approximate dynamic programming
  approach to optimization of train timetable and skip-stop plan for metro
  networks}.
\newblock \bibinfo{journal}{Transportation Research Part C: Emerging
  Technologies} \bibinfo{volume}{157}, \bibinfo{pages}{104393}.



\end{thebibliography}
	\end{document}